\documentclass[12pt,vatola]{article}

\usepackage{graphicx}

\def\c{\centerline}

\def\re#1{\par\hangindent\parindent\indent\llap{#1\enspace}\ignorespaces}

\def\no{\noindent}

\begin{document}

\thispagestyle{empty}

\vskip 45mm

\c{\bf\large AUTOMORPHISM GROUPS OF MAPS, SURFACES}

\vskip 3mm

\c{\bf\large AND SMARANDACHE GEOMETRIES }

\vskip 15mm

\c{\bf  LINFAN MAO}

\vskip 8mm

\hskip 40mm {Institute of Systems Science}

\hskip 40mm {Academy of Mathematics and Systems}

\hskip 40mm {Chinese Academy of Sciences}

\hskip 40mm {Beijing 100080, P.R.China}

\hskip 40mm {maolinfan@163.com}

\newpage

\pagenumbering{roman}

\vskip 25mm

\no{\bf Abstract}:\vskip 25mm

\no A combinatorial map is a connected topological graph
cellularly embedded in a surface. This monograph concentrates on
the automorphism group of a map, which is related to the
automorphism groups of a Klein surface and a Smarandache manifold,
also applied to the enumeration of unrooted maps on orientable and
non-orientable surfaces. A number of results for the automorphism
groups of maps, Klein surfaces and Smarandache manifolds and the
enumeration of unrooted maps underlying a graph on orientable and
non-orientable surfaces are discovered. An elementary
classification for the closed $s$-manifolds is found. Open
problems related the combinatorial maps with the differential
geometry, Riemann geometry and Smarandache geometries are also
presented in this monograph for the further applications of the
combinatorial maps to the classical mathematics.

\vskip 15mm

\no{\bf AMS(2000)} 05C15, 05C25, 20H15, 51D99, 51M05

\newpage

%%%%%%%%%%%%%%%Headings%%%%%%%%%%%%%%%%%%%%%%%%%%%%%%%%%%%%%%%%%
\thispagestyle{empty} \pagestyle{myheadings} \topmargin 5mm
\headheight 8mm \headsep 10mm

\markright {\scriptsize Contents}
%%%%%%%%%%%%%%%%%%%%%%%%%%%%%%%%%%%%%%%%%%%%%%%%%%%%%%%%%%%%%%%%

\vskip 35mm

\no{\large\bf Contents}

\vskip 15mm

\no{\bf Abstract}\dotfill i

\vskip3mm

\no{\bf Chapter $1$ Preliminary}\dotfill 1

\no $\S 1$ Klein surface and $s$-manifolds\dotfill 1

\no $1.1$ Definitions \dotfill 1

\no$1.2$ Classification of Klein surfaces and $s$-manifolds
\dotfill 3

\no $\S 2$ Map and embedding of a graph on surface\dotfill 4

\no $2.1$ Graphs \dotfill 4

\no $2.2$ The embedding of a graph on surfaces\dotfill 5

\no $2.3$ Map and rooted map on surface\dotfill 6

\no $2.4$ Classification maps and embeddings of a graph on a
surfaces\dotfill 8

\no $2.5$ Maps as a combinatorial model of Klein surfaces and
$s$-manifolds \dotfill 10

\no $\S 3$ The semi-arc automorphism group of a graph with
application to maps enumeration\dotfill 11

\no $3.1$ The semi-arc automorphism group of a graph\dotfill 11

\no $3.2$ A scheme for enumerating maps underlying a graph\dotfill
13

\no $\S 4$ A relation among the total embeddings and rooted maps
of a graph on genus\dotfill 16

\no $4.1$ The rooted total map and embedding polynomial of a
graph\dotfill 16

\no $4.2$ The number of rooted maps underlying a graph on
genus\dotfill 20

 \no{\bf Chapter $2$
On the automorphisms of a Klein Surface and a
$s$-manifold}\dotfill 26

\no $\S 1$ An algebraic definition of a voltage map\dotfill 26

\no$1.1$ Coverings of a map\dotfill 26

\no $1.2$ Voltage maps\dotfill 28

\no $\S 2$ Combinatorial conditions for a group being that of a
map\dotfill 29

\no $2.1$ Combinatorial conditions for an automorphism group of a
map\dotfill 30

\no $2.2$ The measures on a map\dotfill 34

\no $\S 3$ A combinatorial refinement of Huriwtz theorem\dotfill
37

\no $\S 4$ The order of an automorphism of a Klein surface\dotfill
42

\no $4.1$ The minimum genus of a fixed-free automorphism\dotfill
42

\no $4.2$ The maximum order of an automorphism of a map\dotfill 45

\no{\bf Chapter $3$ On the Automorphisms of a graph on
Surfaces}\dotfill 48

\no $\S 1$ A necessary and sufficient condition for a group of a
graph being\\ that of a map\dotfill 48

\no $\S 2$ The automorphisms of a complete graph on
surfaces\dotfill 53

\no $\S 3$ The automorphisms of a semi-regular graph on
surfaces\dotfill 61

\no $\S 4$ The automorphisms of one face maps\dotfill 64

 \no{\bf Chapter $4$
Application to the Enumeration of Unrooted Maps and
$s$-Manifolds}\dotfill 68

\no $\S 1$ The enumeration of unrooted complete maps on
surfaces\dotfill 68

\no $\S 2$ The enumeration of a semi-regular graph on
surfaces\dotfill 74

\no $\S 3$ The enumeration of a bouquet on surfaces\dotfill 77

\no $\S 4$ A classification of the closed $s$-manifolds\dotfill 82

\no{\bf Chapter $5$ Open Problems for Combinatorial Maps}\dotfill
88

\no $5.1$ The uniformization theorem for simple connected Riemann
surfaces\dotfill 88

\no $5.2$ The Riemann-Roch theorem \dotfill 89

\no $5.3$ Combinatorial construction of an algebraic curve of
genus \dotfill 89

\no $5.4$ Classification of $s$-manifolds by maps\dotfill 90

\no $5.5$ Gauss mapping among surfaces\dotfill 90

\no $5.6$ The Gauss-Bonnet theorem \dotfill 91

\no $5.7$ Riemann manifolds\dotfill 91

\no{\bf References}\dotfill 92

\newpage

%%%%%%%%%%%%%%%Headings%%%%%%%%%%%%%%%%%%%%%%%%%%%%%%%%%%%%%%%%%
\thispagestyle{empty} \pagestyle{myheadings} \topmargin 5mm
\headheight 8mm \headsep 10mm

\markright {\scriptsize Chapter $1$\quad Preliminary}
%%%%%%%%%%%%%%%%%%%%%%%%%%%%%%%%%%%%%%%%%%%%%%%%%%%%%%%%%%%%%%%%

\vskip 30mm

\pagenumbering{arabic}

\no{\bf\large Chapter $1$\quad Preliminary}

\vskip 10mm

All surfaces  are 2-dimensional compact manifolds without
boundary, graphs are connected, possibly with loops or multiple
edges and groups are finite in the context. For terminology and
notation not defined in this book can be seen in $[33],[34]$ and
$[35]$ for graphs and maps and in $[6],[73]$ for groups.

\vskip 8mm

{\bf\S 1. Klein surface and $s$-manifold}

\vskip 6mm

{\bf $1.1$ Definitions}\vskip 4mm

{$1.1.1$ Definition of a Klein surface}\vskip 3mm

The notion of {\it Klein surface} is introduced by Alling and
Greenleaf [$2$] in $1971$ concerned with real algebraic curves,
correspondence with that of {\it Riemann surface} concerned with
complex algebraic curves. For introducing this concept, it is need
to enlarge analytic functions to those of dianalytic functions
first.

Now let $f: A\longrightarrow  {\cal C}$ be a mapping. Write
$z=x+iy, x,y\in {\cal R}, i=\sqrt {-1}, \overline{z} =x-iy $ and
$f(z)=u(x,y)+iv(x,y)$ $\overline{f(z)}=u(x,y)- iv(x,y)$for certain
functions $u,v:A \longrightarrow {\cal R}.$

A mapping $f: A\longrightarrow  {\cal C}$ is {\it analytic} on $A$
if $\frac{\partial f}{\partial\overline{z}}=0$ ({\it
Cauchy-Riemann equation}) and is {\it antianlytic} on $A$ if
$\frac{\partial f}{\partial z}=0.$

A mapping $f$ is said to be {\it dianalytic} if its restriction to
every connected component of $A$ is analytic or antianalytic.

Now we can formally define a Klein surface.

A {\it Klein surface} is a Hausdorff, connected, topological space
$S$ together with a family $\sum =\{(U_i,\phi_i)$ $ | i \in I \}$
such that the chart $\{U_i |i\in I\}$  is an open covering of $S$,
each map $\phi_i: U_i\longrightarrow A_i$ is a homeomorphism onto
an open subset $A_i$ of ${\cal C}$ or ${\cal C}^+=\{z\in{\cal C}:
Imz\geq 0\}$ and the transition functions

$$
\phi_{ij} =\phi_i\phi_j^{-}: \phi_j (U_i\bigcap U_j
)\longrightarrow  \phi_i (U_i\bigcap U_j ).
$$

\no{are dianalytic.}

The family $\sum$ is said to be a topological {\it atlas} on $S$.

The {\it boundary} of $S$ is the set

\begin{eqnarray*}
\partial S &=& \{x\in S |there \ exists \ i\in I,x\in U_i,\phi_i(x) \in {\cal
R}\\
&and& \phi_i (U_i)\subseteq {\cal C}^+\}.
\end{eqnarray*}

The existence of Klein surfaces is obvious, for example, a Riemann
surface is a Klein surface viewed as an orientable surface with
empty boundary and $\sum$ to be analytic functions. Whence, we
have the following relation:

$$\{Riemann \ Sufaces\}\subset \{Klein \ surfaces\}.$$

The upper half plane $H=\{z\in {\cal C}|Imz > 0\}$ with
$\{(U_1=H,\phi_1=1_H)\}$ and the open unit disc $D=\{z\in {\cal
C}| |z| < 1\}$ with $\{(U_1=d,\phi_1=1_D)\}$ in ${\mathcal C}$ are
two Klein surfaces with empty boundary and analytic structures.

If $k(S), g(S)$ and $\chi (S)$ are the number of connected
components of $\partial S$, the topological genus and the Euler
characteristic of a surface $S$, then we have that

\vskip 3mm

\no{\bf Theorem $1.1.1$}([2])

\[
\chi (S)=\left\{\begin{array}{lr}
2-2g(S)-k(S),& {\rm if} \ S \ {\rm is \ orientable} ,\\
2-g(S)-k(S).& {\rm if} \ S {\rm is \ non-orientable.}
\end{array}
\right.
\]

{ $1.1.2$ Definition of a Smarandache geometry}\vskip 3mm

By the history, we know that classical geometries include the
Euclid geometry, the hyperbolical geometry and the Riemann's
geometry. Each of the later two is proposed by denial the 5th
postulate for parallel lines in the Euclid postulates of geometry.
The {\it Smarandache geometries} is proposed by Smarandache in
1969 ($[61]$), which is a generalization of the classical
geometries, i.e., the Euclid, Lobachevshy-Bolyai-Gauss and
Riemannian geometries may be united altogether in the same space,
by some Smarandache geometries. These last geometries can be
either partially Euclidean and partially Non-Euclidean, or
Non-Euclidean. It seems that the Smarandache geometries are
connected with the {\it Relativity Theory} (because they include
the Riemann geometry in a subspace) and with the {\it Parallel
Universes} (because they combine separate spaces into one space)
too([$32$]).

In $[61]$, Smarandache defined several specific types of
Smarandache geometries, such as the {\it paradoxist geometry}, the
{\it non-geometry}, the {\it counter-projective geometry} and the
{\it anti-geometry}. He also posed a question on the paradoxist
geometry, i.e., {\it find a nice model on manifolds for this
paradoxist geometry and study some of its characteristics}.

An axiom is said {\it smarandachely denied} if in the same space
the axiom behaves differently (i.e., validated and invalided; or
only invalided but in at least two district ways).

A {\it Smarandache geometry} is a geometry which has at least one
smarandachely denied axiom$^{1*}$ \footnotetext[1]{Also see {\it
www.gallup.unm.edu/\~\ samrandache/geometries.htm}}. At present,
the Smarandache manifolds (abbreviated $s$-manifolds) are the
central object discussed in the Smarandache geometries today. More
results for the Smarandache geometries can be seen in the
references $[4],[16]$,$[27]-[28]$, $[32]$ and $[58]-[59]$ etc..

The idea of an $s$-manifold was based on a hyperbolic paper in
$[69]$ and credited to W.Thurston. A more general idea can be
found in $[59]$. According to the survey $[27]$ of H.Iseri, an
{\it $s$-manifold} is combinatorially defined as follows:

{\it An $s$-manifold is any collection ${\mathcal C}(T,n)$ of
these equilateral triangular disks $T_i, 1\leq i\leq n$ satisfying
the following conditions:

$(i)$ Each edge $e$ is the identification of at most two edges
$e_i,e_j$ in two distinct triangular disks $T_i,T_j, 1\leq i,j\leq
n$ and $i\not= j$;

$(ii)$ Each vertex $v$ is the identification of one vertex in each
of five, six or seven distinct triangular disks.}

The vertices are classified by the number of the disks around
them. A vertex around five, six or seven triangular disks is
called an {\it elliptic vertex}, a {\it Euclid vertex} or a {\it
hyperbolic vertex}, respectively.

An $s$-manifold is called closed if the number of triangular disks
is finite and each edge is shared by exactly two triangular disks,
each vertex is completely around by triangular disks. It is
obvious that a closed $s$-manifold is a surface and its Euler
characteristic can be defined by the Theorem $1.1.1$.

\vskip 4mm

{\bf $1.2$ Classification of Klein surfaces and $s$-manifolds}

\vskip 3mm

A {\it morphism} between the Klein surfaces $S$ and $S'$ is a
continuous map $f: S \rightarrow S'$ such that $f(\partial
S)\subseteq \partial S'$ and given $s\in S$, there exist chart
$(U,\phi)$ and $(V,\psi)$ at $s$ and $f(s)$ respectively, and an
analytic function $F:\phi (U)\rightarrow {\mathcal C}$ such that

$$\psi (f(s))=\Phi (F(\phi (s))),$$

\no{where, $\Phi:{\mathcal C}\rightarrow {\mathcal C}^+: x+iy
\rightarrow x+i|y|$ is a continuous map.}

An {\it automorphism} of a Klein surface $S$ is an $1-1$ morphism
$f:S \rightarrow S$. It has been known that for a given Klein
surface $S$, the set ${\rm Aut}S$ of automorphisms of $S$ forms a
group with respect to the composition operation and ${\rm Aut}H=
PGL(2,R)$.

Let $\Gamma$ be a discrete subgroup of Aut$H$. We say that
$\Gamma$ is a {\it non-euclidean crystallographic group}( shortly
$NEC$ group) if the quotient $H/\Gamma$ is compact.

More results can be seen in $[11]$. Typical results for
automorphisms of a Klein surface $S$ are as follows.

\vskip 4mm

\no{\bf Theorem $1.1.2$}([11]) {\it Let $S$ be a compact Klein
surface, $g=g(S)$ and $k=k(S)$, then

($i$) there exists an $NEC$ group $\Gamma$ such that ${\rm
Aut}S\cong N_{\Omega}(\Gamma)/\Gamma$, where $\Omega = {\rm
Aut}H$.

($ii$) if $S$ satisfies the condition $2g+k\geq 3$ if $S$ is
orientable and $g+k\geq 3$ if $S$ is non-orientable, then ${\rm
Aut}S$ is finite.}

\vskip 3mm

Similarly, two $s$-manifolds ${\mathcal C}_1(T,n)$ and  ${\mathcal
C}_2(T,n)$ are called to be {\it isomorphic} if there is an $1-1$
mapping $\tau: {\mathcal C}_1(T,n)\rightarrow {\mathcal C}_2(T,n)$
such that for $\forall T_1,T_2\in {\mathcal C}_1(T,n)$,

$$\tau (T_1\bigcap T_2)=\tau (T_1)\bigcap\tau (T_2).$$

\no{If ${\mathcal C}_1(T,n)={\mathcal C}_1(T,n)={\mathcal
C}(T,n)$, $\tau$ is called an {\it automorphism} of the
$s$-manifold ${\mathcal C}(T,n)$. All automorphisms of an
$s$-manifold form a group under the composition operation, called
the automorphism group of an $s$-manifold ${\mathcal C}(T,n)$,
denoted by ${\rm Aut}{\mathcal C}(T,n)$.}

\vskip 6mm

{\bf\S 2. Map and embedding of a graph on surface}

\vskip 4mm

{\bf $2.1$ Graphs}

\vskip 3mm

Combinatorially, a {\it graph} $\Gamma$ is a $2$-tuple $(V,E)$
consists of a finite non-empty set $V$ of vertices together with a
set $E$ of unordered pairs of vertices, i.e., $E\subseteq V\times
V$($[22], [35], [70]$). Often denoted by $V(\Gamma )$, $E(\Gamma
)$ the vertex set and edge set of the graph $\Gamma$.

The cardinal numbers of $|V|$ and $|E|$ are called the {\it order}
and the {\it size} of the graph $\Gamma$.

We can also obtain a representation of a graph $\Gamma$
representing a vertex $u$ by a point $p(u)$, $p(u)\not=p(v)$ if
$u\not= v$ and an edge $(u,v)$ by a curve connecting the points
$p(u)$ and $p(v)$ on the plane.

For example, the graph in the Fig. $1.1$

\includegraphics[bb=5 10 200 180]{sg1.eps}\vskip 2mm

\c{\bf Fig $1.1$} \vskip 3mm

\no{is a graph $\Gamma =(V,E)$ with $V=\{u,v,w,x\}$ and

$$E=\{(u,u),(v,v),(w,w), (x,x),(u,v),(v,w),(w,x),(x,u)\}.$$

A {\it walk} of a graph $\Gamma$ is an alternating sequence of
vertices and edges $u_1, e_1,u_2,e_2,$ $\cdots,e_n, u_{n_1}$ with
$e_i=(u_i,u_{i+1})$ for $1\leq i\leq n$. The number $n$ is the
length of the walk. If $u_1=u_{n+1}$, the walk is said to be {\it
closed}, and {\it open} otherwise. For example,
$ue_1ve_2we_6we_3xe_3we_2v$ is a walk in the graph of the Fig.
$1.1$. A walk is called a {\it trail} if all its edges are
distinct and a {\it path} if all the vertices are distinct. A
closed path is said to be a circuit.

A graph $\Gamma$ is {\it connected} if there is paths connecting
any two vertices in this graph and is {\it simple} if any
$2$-tuple $(u,v)\in E(\Gamma )\subseteq V(\Gamma )\times V(\Gamma
)$ appears once at most and $u\not= v$.

Let $\Gamma$ be a graph. For $\forall u\in V(\Gamma )$, the
neighborhood $N_{\Gamma}^v(u)$ is defined by
$N_{\Gamma}^v(u)=\{v|(u,v)\ {\rm or} \ (v,u)\in E(\Gamma )\}$. Its
cardinal $|N_{\Gamma}^v(u)|$ is called the valency of the vertex
$u$ in the graph $\Gamma$, denoted by $\rho_{\Gamma}(u)$. By the
enumeration of edges, we know the following result

$$\sum\limits{u\in V(\Gamma )}\rho_{\Gamma}(u)= 2|E(\Gamma )|.$$

\vskip 6mm

{\bf $2.2$ The embedding of a graph on surfaces}

\vskip 4mm

A map on a surface $S$ is a kind of partition $S$ which enables us
to obtain homeomorphisms of 2-cells $\{(x,y)|x^2+y^2 \ < \ 1,
x,y\in {\mathcal R}\}$ if we remove from $S$ all the curves used
to partite $S$. There is a classical result for the partition of a
surface gotten by T.Rad\'{o} in 1925.

\vskip 4mm

\no{\bf Theorem $1.2.1$}([52]) {\it For any compact surface $S$,
there exist a triangulation $\{{\mathcal T}_i,i\geq 1\}$ on $S$.}

\vskip 3mm

This theorem is fundamental for the topological graph theory,
which enables us to discussion a surface combinatorially.

For any connected graph $\Gamma = (V(\Gamma),E(\Gamma))$ and a
surface $S$, an {\it embedding} of the graph $\Gamma$ in the
surface $S$ is geometrical defined to be a continuous $1-1$
mapping $\tau: \Gamma\rightarrow S$. The image $\tau (\Gamma)$ is
contained in the $1$-skeleton of a triangulation of the surface
$S$.  If each component in $S - \tau (\Gamma)$ homeomorphic to an
open disk, then the embedding is said a $2$-cell embedding, where,
components in  $S - \tau (\Gamma)$ are called faces. All
embeddings considered in this book are $2$-cell embeddings.

Let $\Gamma$ be a graph. For $v\in V(\Gamma)$, denote by
$N_{\Gamma}^e(v)=\{e_1,e_2,\cdots, e_{\rho (v)}\}$ all the edges
incident with the vertex $v$. A permutation on $e_1,e_2,\cdots,
e_{\rho (v)}$ is said a {\it pure rotation}. All pure rotations
incident a vertex $v$ is denoted by $\varrho (v)$. A {\it pure
rotation system} of the graph $\Gamma$ is defined to be

$$\rho (\Gamma) = \{\varrho (v)| v\in V(\Gamma)\}$$

\no{and all pure rotation systems of the graph $\Gamma$ is denoted
by $\varrho (\Gamma)$.}

The first characteristic for embedding of a graph on orientable
surfaces is found by Heffter in 1891 and formulated by Edmonds in
1962, states as follows.

\vskip 4mm

\no{\bf Theorem $1.2.2$}([17]) {\it Every pure rotation system for
a graph $\Gamma$ induces a unique embedding of $\Gamma$ into an
orientable surface. Conversely, every embedding of a graph
$\Gamma$ into an orientable surface induces a unique pure rotation
system of $\Gamma$. }

\vskip 3mm

According to this theorem, we know that the number of orientable
embeddings of a graph $\Gamma$ is $\prod_{v\in V(\Gamma)}(\rho
(v)-1)!$.

The characteristic for embedding of a graph on locally orientable
surface is used by Ringel in the 1950s and gave a formal proof by
Stahl in $1978 ([22][62])$.

From the topological theory, embedded vertex and face can be
viewed as disk, and an embedded edge can be viewed as an $1$-band
which is defined as a topological space $B$ together with a
homeomorphism $h:I\times I \rightarrow B$, where $I=[0,1]$, the
unit interval.

Define a rotation system $\rho^L (\Gamma)$ to be a pair $({\cal
J},\lambda)$, where ${\cal J}$  is a pure rotation system of
$\Gamma$, and $\lambda : E(\Gamma) \rightarrow Z_2$. The edge with
$\lambda (e)= 0$ or $\lambda (e)=1$ is called {\it type $0$} or
{\it type $1$} edge, respectively.  The {\it rotation system} of a
graph $\Gamma$ are defined by

$$\varrho^L (\Gamma) = \{({\cal
J},\lambda)| {\cal J}\in \varrho (\Gamma),\lambda : E(\Gamma)
\rightarrow Z_2 \}.$$

Then we know that

\vskip 4mm

\no{\bf Theorem $1.2.3$}([22][62]) {\it Every rotation system on a
graph $\Gamma$ defines a unique locally orientable embedding of
$\Gamma\rightarrow S$. Conversely, every embedding of a graph
$\Gamma\rightarrow S$ defines a rotation system for $\Gamma$. }

\vskip 3mm

For any embedding of the graph $\Gamma$, there is a spanning tree
$T$ such that every edge on this tree is type $0 ([43])$. Whence
the number of embeddings of a graph $\Gamma$ on locally orientable
surfaces is

$$2^{\beta (\Gamma)}\prod_{v\in V(\Gamma)}(\rho
(v)-1)!$$

\no{and the number of embeddings of $\Gamma$ on the non-orientable
surfaces is}

$$(2^{\beta (\Gamma)}-1)\prod_{v\in V(\Gamma)}(\rho
(v)-1)!.$$

The following result is the {\it Euler-Poincar\'{e}} formula for
an embedding of a graph on surface.

\vskip 4mm

\no{\bf Theorem $1.2.4$} {\it If a graph $\Gamma$ can be embedded
into a surface $S$, then

$$\nu(\Gamma)-\varepsilon(\Gamma)+\phi (\Gamma) = \chi (S),$$

\no{where, $\nu(\Gamma), \varepsilon(\Gamma)$ and $\phi (\Gamma) $
are the order, size and the number of faces of the graph $\Gamma$,
and $\chi (S)$ is the Euler characteristic of the surface $S$:}}

\[
\chi (S)=\left\{\begin{array}{lr}
2-2p, & if \ S \ is \  orientable,\\
2-q, &  if \ S \ is \  non-orientable.
\end{array}
\right.
\]

\vskip 5mm

{\bf $2.3.$ Map and rooted map on surface}

\vskip 3mm

In 1973, Tutte gave an algebraic representation for an embedding
of a graph on locally orientable surface $([66]$, which transfer a
geometrical partition of a surface to a kind of permutation in
algebra.

According to the summary in $[33]$, a {\it{map}} $M = ({\cal X}
_{\alpha,\beta},\cal{P})$ is defined  to be a basic permutation
$\cal{P}$, i.e, for any $x\in {\cal X}_{\alpha,\beta}$, no integer
$k$ exists such that  ${\cal{P}}^{k}x = \alpha x$, acting on
${\cal X} _{\alpha,\beta}$, the disjoint union of {\it
quadricells} $Kx$ of $x\in  X$ (the base set), where
$K=\{1,\alpha,\beta,\alpha\beta \}$ is the {\it Klein group},
satisfying the following two conditions:\vskip 3mm

{\it ($Ci$) $\alpha{\cal{P}}={\cal{P}}^{-1}\alpha$;\vskip 2mm

($Cii$) the group $\Psi_{J}=<\alpha,\beta,\cal{P}>$ is transitive
on ${\cal X}_{\alpha,\beta}$.}\vskip 2mm

For a given map $M=({\mathcal X}_{\alpha ,\beta},{\mathcal P})$,
it can be shown that $M^* = ({\mathcal X}_{\beta ,\alpha
},{\mathcal P}\alpha\beta)$ is also a map, call it the {\it dual}
of the map $M$. The vertices of $M$ are defined as the pairs of
conjugatcy orbits of ${\mathcal P}$ action on ${\mathcal
X}_{\alpha,\beta}$ by the condition $(Ci)$ and edges the orbits of
$K$ on ${\mathcal X}_{\alpha,\beta}$, for example,$\forall x\in
{\mathcal X}_{\alpha,\beta}$, $\{x,\alpha x,\beta x,\alpha\beta
x\}$ is an edge of the map $M$. Define the faces of $M$ to be the
vertices in the dual map $M^*$. Then the Euler characteristic
$\chi (M)$ of the map $M$ is

$$\chi (M)= \nu (M)-\varepsilon (M)+\phi (M)$$

\no{where,$\nu (M), \varepsilon (M), \phi (M)$ are the number of
vertices, edges and faces of the map $M$, respectively.}

For example, the graph $K_4$ on the tours with one face length $4$
and another $8$ , shown in the Fig. $1.2$, can be algebraic
represented as follows:\vskip 3mm

A map $({\mathcal X}_{\alpha,\beta},\mathcal{P})$ with ${\mathcal
X}_{\alpha,\beta}= \{x,y,z,u,v,w,\alpha x,\alpha y, \alpha
z,\alpha u,\alpha v,\alpha w, \beta x,\beta y,\beta z,$ $\beta u,
 \beta v,\beta w,\alpha\beta x,
\alpha \beta y,\alpha \beta z,\alpha \beta u,\alpha \beta
v,\alpha\beta w \}$ and

\begin{eqnarray*}
{\mathcal P} &=& (x,y,z)(\alpha \beta x,u,w)(\alpha \beta z,\alpha
\beta u,v)
(\alpha \beta y,\alpha \beta v,\alpha \beta w)\\
&\times& (\alpha x,\alpha z,\alpha y)(\beta x,\alpha w,\alpha
u)(\beta z,\alpha v,\beta u)(\beta y,\beta w,\beta v)
\end{eqnarray*}

\no{The four vertices of this map are $ \{(x,y,z), (\alpha
x,\alpha z,\alpha y)\}$, $\{(\alpha \beta x,u,w),(\beta x,\alpha
w,\alpha u)\}$, $\{(\alpha \beta z,\alpha \beta u,v),(\beta
z,\alpha v,\beta u)\}$ and $\{(\alpha \beta y,\alpha \beta
v,\alpha \beta w),(\beta y,\beta w,\beta v)\}$ and six edges are
$\{e,\alpha e,\beta e,\alpha\beta e\}$, where, $e\in
\{x,y,z,u,v,w\}$. The Euler characteristic $\chi (M)$ is $\chi
(M)=4-6+2=0$.}

\includegraphics[bb=5 10 200 200]{7mg1.eps}

\vskip 2mm \c{\bf Fig $1.2$} \vskip 2mm

Geometrically, an embedding $M$ of a graph $\Gamma$ on a surface
is a map and has an algebraic representation. The graph $\Gamma$
is said the {\it underlying graph} of the map $M$ and denoted by
$\Gamma =\Gamma (M)$. For determining a given map $({\mathcal
X}_{\alpha,\beta},\mathcal{P})$ is orientable or not, the
following condition is needed.\vskip 3mm

{\it ($Ciii$) If the group $\Psi_I=<\alpha\beta ,{\mathcal P}>$ is
transitive on ${\mathcal X}_{\alpha,\beta}$, then $M$ is
non-orientable. Otherwise, orientable.}\vskip 2mm

It can be shown that the number of orbits of the group
$\Psi_I=<\alpha\beta ,{\mathcal P}>$ in the Fig.$1.1$ action on
${\mathcal X}_{\alpha,\beta}= \{x,y,z,u,v,w,\alpha x,\alpha y,$ $
\alpha z,\alpha u,\alpha v,\alpha w, \beta x,\beta y,\beta z,\beta
u,$ $ \beta v,\beta w,\alpha\beta x, \alpha \beta y,\alpha \beta
z,\alpha \beta u,\alpha \beta v,\alpha\beta w \}$ is $2$. Whence,
it is an orientable map and the genus of the surface is $1$.
Therefore, the algebraic representation is correspondent with its
geometrical mean.

A {\it rooted map} $M^x$ is a map $M$ such that one quadricell
$x\in{\mathcal X}_{\alpha ,\beta}$ is marked beforehand, which is
introduced by Tutte in the enumeration of planar maps. The
importance of the root is to destroy the symmetry in a map. That
is the reason why we can enumerate rooted maps on surfaces by
combinatorial approaches.

\vskip 5mm

{\bf 2.4. Classification maps and embeddings of a graph on
surfaces}

\vskip 3mm { $2.4.1$ Equivalent embeddings of a graph}\vskip 2mm

From references, such as, $[22],[70],$ etc., two embeddings
$({\cal J}_1,\lambda_1),({\cal J}_2,\lambda_2)$ of $\Gamma$ on an
orientable surface $S$ are called {\it equivalent} if there exists
an orientation-preserving homeomorphism $\tau$ of the surface $S$
such that $\tau : {\cal J}_1 \rightarrow {\cal J}_2$, and
$\tau\lambda =\lambda\tau$. If $({\cal J}_1,\lambda_1)=({\cal
J}_2,\lambda_2)=({\cal J},\lambda )$, then an
orientation-preserving homeomorphism mapping $({\cal
J}_1,\lambda_1)$ to $({\cal J}_2,\lambda_2)$ is called an
automorphism of the embedding $({\cal J},\lambda )$. Certainly,
all automorphisms of an embedding form a group, denoted by ${\rm
Aut}({\cal J},\lambda )$.

Enumerating the non-equivalent orientable embeddings of a complete
graph and a complete bipartite graph are considered by Biggs,
White, Mull and Rieper et al in $[6],[54]-[55]$. Their approach is
generalized in the following Section $2.3.2$ for enumerating
non-equivalent embeddings of a given graph on locally orientable
surface in the view of maps on surfaces.

\vskip 3mm

{ $2.4.2$ Isomorphism of maps}

\vskip 2mm

Two maps $M_{1} = ({\cal X}_{\alpha,\beta}^{1},{\cal{P}}_{1})$ and
$M_{2} = ({\cal X}_{\alpha,\beta}^{2},{\cal{P}}_{2})$ are said to
be  {\it isomorphic} if there exists a bijection $\xi$

$$\xi: {\cal X}_{\alpha,\beta}^{1} \longrightarrow
{\cal X}_{\alpha,\beta}^{2}$$

\no{such that for $\forall{x}\in{{\cal X}_{\alpha,\beta}^{1}}$,}

$$\xi\alpha(x)=\alpha\xi(x),
\xi\beta(x)=\beta\xi(x)\quad {\rm and}\quad
\xi{\cal{P}}_{1}(x)={\cal{P}}_{2}\xi(x).$$

\no{Call $\xi$ an {\it isomorphism} between $M_1$ and $M_2$. If
$M_{1}=M_{2}=M$, then an isomorphism between $M_{1}$ and $M_{2}$
is called an {\it automorphism} of $M$. All automorphisms of a map
$M$ form a group, called the {\it automorphism group} of $M$ and
denoted by ${\rm Aut}M$.}

Similarly, two rooted maps $M^{x}_{1},$ $M^{y}_{2}$ are said to be
{\it isomorphic} if there is an isomorphism $\theta$ between them
such that $\theta(x)= y$. All automorphisms of a rooted map $M^r$
also form a group, denoted by ${\rm Aut}M^{r}$. It has been known
that ${\rm Aut}M^{r}$ is trivial $([33])$.

Using isomorphisms between maps, an alternative approach for
determining equivalent embeddings of a graph on locally orientable
surfaces can be gotten, which has been used in $[43],[49]-[50]$
for determining the number of non-equivalent embeddings of a
complete graph, a semi-regular graph and a Cayley graph $\Gamma
=Cay(G:S)$ with ${Aut}\Gamma \cong R(G)\times H$, is defined as
follows.

For a given map $M$ underlying a graph $\Gamma$, it is obvious
that ${\rm Aut}M|_{\Gamma}\leq {\rm Aut}\Gamma$. We extend the
action $\forall g\in {\rm Aut}\Gamma$ on $V(\Gamma)$ to ${\mathcal
X}_{\alpha ,\beta}$, where $X=E(\Gamma)$, as follows:

{\it $\forall x\in {\mathcal X}_{\alpha ,\beta}$, if $x^g=y$, then
define $(\alpha x)^g=\alpha y, (\beta x)^g=\beta y$ and
$(\alpha\beta x)^g=\alpha\beta y$.}

Two maps (embeddings) $M_1,M_2$ with a given underlying graph
$\Gamma$ are equivalent if there exists an isomorphism $\zeta$
between them induced by an element $\xi$. Call $\zeta$ an
equivalence between $M_1,M_2$. Certainly, on an orientable
surface, an equivalence preserve the orientation on this surface.

\vskip 3mm

\no{\bf Theorem $1.2.5$} {\it Let $M=({\mathcal X}_{\alpha
,\beta},{\mathcal P})$ be a map with an underlying graph $\Gamma$,
$\forall g\in {Aut}\Gamma$. Then the extend action of $g$ on
${\mathcal X}_{\alpha ,\beta}$ with $X=E(\Gamma)$ is an
automorphism of map $M$ iff \ $\forall v\in V(M)$, $g$ preserves
the cyclic order of $v$. }

\vskip 2mm

{\it Proof} Assume that $\zeta\in {Aut}M$ is induced by the extend
action of an automorphism $g$ in $\Gamma$, $u,v\in V(M)$ and
$u^g=v$. Not loss of the generality, we assume that

$$u=(x_1,x_2,\cdots,x_{\rho (u)})(\alpha x_{\rho (u)} ,\cdots, \alpha x_2,\alpha x_1)$$

$$v=(y_1,y_2,\cdots,y_{\rho (v)})(\alpha y_{\rho (v)} ,\cdots, \alpha y_2,\alpha y_1).$$

Without loss of generality , we can assume that

$$(x_1,x_2,\cdots,x_{\rho (u)})^g = (y_1,y_2,\cdots,y_{\rho (v)}),$$

\no{that is,}

$$(g(x_1),g(x_2),\cdots,g(x_{\rho (u)})) = (y_1,y_2,\cdots,y_{\rho (v)}).$$

\no{Whence, $g$ preserves the cyclic order of vertices in the map
$M$.}

On the other hand, if the extend action of $g\in Aut\Gamma$ on
${\mathcal X}_{\alpha ,\beta}$ preserves the cyclic order of each
vertex in $M$, i.e., $\forall u\in V(\Gamma),\exists v\in
V(\Gamma)$ such that $u^g=v$. Assume that

$${\mathcal P}=\prod_{u\in V(M)}u.$$

\no{Then}

$${\mathcal P}^g=\prod_{u\in V(M)}u^g = \prod_{v\in V(M)}v={\mathcal P}.$$

\no{Therefore, the extend action of $g$ on ${\mathcal X}_{\alpha
,\beta}$ is an automorphism of the map $M$. \quad\quad $\natural$}

\vskip 5mm

{\bf $2.5$ Maps as a combinatorial model of Klein surfaces and
$s$-manifolds}

\vskip 3mm

{$2.5.1$ The model of Klein surfaces}

\vskip 2mm

Given a complex algebraic curve, it is a very important problem to
determine its birational automorphisms. For curve ${\mathcal C}$
of genus $g\geq 2$, Schwarz proved that ${\rm Aut}({\mathcal C})$
is finite in 1879 and Hurwitz proved $|Aut({\mathcal C})|\leq
84(g-1)$(see $[18]$ ). As observed by Riemann, groups of
birational automorphisms of complex algebraic curves are the same
as the automorphism groups of the compact Riemann surfaces. The
latter can be combinatorially dealt with the approach of maps.

\vskip 3mm

\no{\bf Theorem $1.2.6$}([8][29]) {\it If $M$ is an orientable map
of genus $p$, then ${\rm Aut}M$ is isomorphic to a group of
conformal transformations of a Riemann surface of genus $p$.}

\vskip 2mm

According to the Theorem $1.1.2$, the automorphism group of a
Klein surface has the same form as a Riemann surface. Similar to
the proof of the Theorem $5.6$ in $[29]$, we can get a result
similar to the Theorem $1.2.6$ for Klein surfaces.

\vskip 3mm

\no{\bf Theorem $1.2.7$} {\it If $M$ is a locally orientable map
of genus $q$, then ${\rm Aut}M$ is isomorphic to a group of
conformal transformations of a Klein surface of genus $q$.}

\vskip 2mm

{\it Proof} By a result in $[8]$, ${\rm Aut}M\cong N_T(A)/A$,
Where $T = <a,b,c|a^2=b^2=c^2=(ba)^2=(ac)^m=(cb)^n=1>$, $A\leq T$
and $T$ can be realized by an automorphism group of a tessellation
on the upper plane, $A$ a $NEC$ subgroup. According to the Theorem
$1.1.2$, The underlying surface $S$ of $M$ has $S= H/A$ with
$\Omega ={\rm Aut}H = PGL(2,{\mathcal R})$ being the automorphism
group of the upper half plane $H$. Since $T\leq \Omega$, we know
that ${\rm Aut}M\cong N_T(A)/A \leq N_{\Omega}(A) /A$, isomorphic
to a group of conformal transformations of the Klein surface $S=
H/G$. \quad\quad $\natural$

\vskip 3mm

{$2.5.2$ The model of closed $s$-manifolds}

\vskip 2mm

For a closed $s$-manifold ${\mathcal C}(T,n)$, we can define a map
$M$ by $V(M)=\{{\rm the \ vertices \ in \ }$ $ {\mathcal
C}(T,n)\},$ $E(M)=\{{\rm the \ edges \ in \ } {\mathcal C}(T,n)\}
$  and $F(M) =\{T, T\in {\mathcal C}(T,n) \}.$ Then, $M$ is a
triangular map with vertex valency $\in\{5,6,7\}$. On the other
hand, if $M$ is a triangular map on surface with vertex valency
$\in\{5,6,7\}$, we can define ${\mathcal C}(T,\phi (M))$ by

$${\mathcal C}(T,\phi (M))= \{f | f\in F(M)\}.$$

\no{Then, ${\mathcal C}(T,\phi (M))$ is an $s$-manifold.
Therefore, we get the following result. } \vskip 3mm

\no{\bf Theorem $1.2.8$} {\it Let $\widehat{{\mathcal C}(T,n)}$,
${\mathcal M}(T,n)$ and ${\mathcal M}^*(T,n)$ be the set of
$s$-manifolds with $n$ triangular disks, triangular maps with $n$
faces and vertex valency $\in\{5,6,7\}$ and cubic maps of order
$n$ with face valency $\in\{5,6,7\}$. Then

$(i)$ There is a bijection between ${\mathcal M}(T,n)$ and
$\widehat{{\mathcal C}(T,n)}$;

$(ii)$ There is also a bijection between ${\mathcal M}^*(T,n)$ and
$\widehat{{\mathcal C}(T,n)}$.}

\vskip 6mm

{\bf\S 3. The semi-arc automorphism group of a graph with
application to maps enumeration}

\vskip 4mm

{\bf $3.1$ The semi-arc automorphism group of a graph}

\vskip 3mm

Let $\Gamma$ be a graph  with vertex set $V(\Gamma)$ and edge set
$E(\Gamma)$. By the definition, an {\it automorphism} of $\Gamma$
on $V(\Gamma)\bigcup E(\Gamma)$ is an $1-1$ mapping $(\xi,\eta)$
on $\Gamma$ such that

$$\xi : V(\Gamma)\rightarrow V(\Gamma), \ \  \eta : E(\Gamma)\rightarrow E(\Gamma),$$

\no{satisfying that for any incident elements $e,f$,
$(\xi,\eta)(e)$ and $(\xi,\eta)(f)$ are also incident. Certainly,
all automorphisms of a graph $\Gamma$ form a group, denoted by
${\rm Aut}\Gamma$.}

Now an edge $e=uv\in E(\Gamma)$ can be divided into two {\it
semi-arcs} $e_u, e_v$. Call $u$ the {\it root vertex} in the
semi-arc $e_u$. Two semi-arc $e_u,f_v$ are said {\it incident} if
$u=v$ or $e=f$. The set of all semi-arcs of a graph $\Gamma$ is
denoted by $X_{\frac{1}{2}}(\Gamma)$. A semi-arc automorphism of a
graph, first appeared in $[43]$ and then applied to the
enumeration rooted maps on surfaces underlying a graph $\Gamma$ in
$[46]$, is defined as follows.

\vskip 3mm

\no{\bf Definition $1.3.1$}  {\it An $1-1$ mapping $\xi$ on
$X_{\frac{1}{2}}(\Gamma)$ such that $\forall e_u,f_v\in
X_{\frac{1}{2}}(\Gamma)$, $\xi (e_u)$ and $\xi (f_v)$ are incident
if $e_u$ and $f_v$ are incident, is called a semi-arc automorphism
of the graph $\Gamma$.}

\vskip 2mm

All semi-arc automorphisms of a graph also form a group under the
composition operation,  denoted by ${\rm
Aut}_{\frac{1}{2}}\Gamma$, which is more important for the
enumeration of maps on surfaces and  also important for determine
the conformal transformations on a Klein surface. The following
table lists semi-arc automorphism groups of some well-known
graphs, which give us some useful information for the semi-arc
automorphism groups, for example, ${\rm Aut}_{\frac{1}{2}}K_n=
S_n$ but ${\rm Aut}_{\frac{1}{2}}B_n= S_n[S_2]\not={\rm Aut}B_n.$

\vskip 3mm
\begin{center}
\begin{tabular}{|c|c|c|} \hline
$\Gamma$ & $Aut_{\frac{1}{2}}\Gamma $ &  $ order$  \\ \hline
$K_n$ & $S_n$ & $n!$  \\
$K_{m,n}(m\not=n)$ & $S_m\times S_n$ & $m!n!$\\
\ \ \ $K_{n,n}$ \ \ \  & \ \ \ $S_2[S_n]$ \ \ \ & \ \ \ $2n!^2$ \ \ \ \\
$B_n$ & $S_n[S_2]$ & $2^nn!$ \\
${Dp}_n$ & $S_2\times S_n$ & $2n!$ \\
${Dp}_n^{k,l} (k\not= l)$ & $S_2[S_k]\times S_n \times S_2[S_l]$ & $2^{k+l}n!k!l!$ \\
${Dp}_n^{k,k}$ & $S_2\times S_n \times (S_2[S_k])^2$ &
$2^{2k+1}n!k!^2$ \\ \hline
\end{tabular}
\end{center}
\vskip 2mm \c{\ table $3.1$} \vskip 3mm

\no{Here, $Dp_n$ is a dipole graph with 2 vertices, $n$ multiple
edges and $Dp_n^{k,l}$ is a generalized dipole graph with 2
vertices, $n$ multiple edges, and one vertex with $k$ bouquets and
another, $l$ bouquets. Comparing the semi-arc automorphism groups
in the second column with automorphism groups of graphs in the
first column in table $3.1$, it is easy to note that the semi-arc
automorphism group are the same as the automorphism group in the
first two cases. In fact, it is so by the following Theorem
$1.3.1$.}

 For $\forall g\in {\rm
Aut}\Gamma$, there is also an induced action $g|^{\frac{1}{2}}$ on
$X_{\frac{1}{2}}(\Gamma)$, $g: X_{\frac{1}{2}}(\Gamma)\rightarrow
X_{\frac{1}{2}}(\Gamma)$, as follows:

$$\forall e_u\in X_{\frac{1}{2}}(\Gamma), g(e_u)=(g(e)_{g(u)}. $$

\no{All induced action of the elements in ${\rm Aut}\Gamma$ on
$X_{\frac{1}{2}}(\Gamma)$ is denoted by ${\rm
Aut}\Gamma|^{\frac{1}{2}}$. Notice that }

$${\rm Aut}\Gamma \cong {\rm Aut}\Gamma|^{\frac{1}{2}}.$$

We have the following result.

\vskip 4mm

\no{\bf Theorem $1.3.1$} {\it For a graph $\Gamma$ without loops,}

$${\rm Aut}_{\frac{1}{2}}\Gamma = {\rm Aut}\Gamma|^{\frac{1}{2}}.$$

\vskip 3mm

{\it Proof} By the definition, we only need to prove that for
$\forall\xi_{\frac{1}{2}}\in {\rm Aut}_{\frac{1}{2}}\Gamma $, $\xi
= \xi_{\frac{1}{2}} |_{\Gamma}\in {\rm Aut}\Gamma$ and
$\xi_{\frac{1}{2}}=\xi |^{\frac{1}{2}}$. In fact, for any $\forall
e_u,f_x\in X_{\frac{1}{2}}(\Gamma)$, where, $e=uv\in E(\Gamma)$
and $f=xy\in E(\Gamma)$, if

$$\xi_{\frac{1}{2}}(e_u)=f_x,$$

\no{then by the definition, we know that}

$$\xi_{\frac{1}{2}}(e_v)=f_y .$$

\no{Whence, $\xi_{\frac{1}{2}}(e)=f$. That is,
$\xi_{\frac{1}{2}}|_{\Gamma}\in {\rm Aut}\Gamma$.}

Now since there is not a loop in $\Gamma$, we know that
$\xi_{\frac{1}{2}}|_{\Gamma} ={\bf id}_{\Gamma}$ if and only if
$\xi_{\frac{1}{2}}={\bf id}_{\Gamma}$. Therefore,
$\xi_{\frac{1}{2}}$ is induced by $\xi_{\frac{1}{2}}|_{\Gamma}$ on
$X_{\frac{1}{2}}(\Gamma)$, that is,

$${\rm Aut}_{\frac{1}{2}}\Gamma = {\rm Aut}\Gamma|^{\frac{1}{2}}.\quad\natural$$

Notice that for a given graph $\Gamma$, $X_{\frac{1}{2}}(\Gamma) =
{\cal X}_{\beta},$ if we equal $e_u$ to $e$ and $e_v$ to $\beta e$
for an edge $e=uv\in E(\Gamma)$.

For a given map $M=({\mathcal X}_{\alpha ,\beta},{\mathcal P})$
underlying a graph $\Gamma$, we have known that ${\rm
Aut}M|_{\Gamma}\leq {\rm Aut}\Gamma$, which made us to extend the
action of an automorphism $g$ of the graph $\Gamma$ on ${\mathcal
X}_{\alpha ,\beta}$ with $X=E(\Gamma )$ to get automorphisms of a
map induced by automorphisms of its underlying graph. More detail,
we can get the following result.

\vskip 4mm

\no{\bf Theorem $1.3.2$} {\it Two maps $M_1=({\mathcal X}_{\alpha
,\beta},{\mathcal P}_1)$ and $M_2=({\mathcal X}_{\alpha
,\beta},{\mathcal P}_2)$ underlying a graph $\Gamma$ are

($i$) equivalent iff there is an element $\zeta\in {\rm
Aut}_{\frac{1}{2}}\Gamma$ such that ${\mathcal
P}_1^{\zeta}={\mathcal P}_2$ and

($ii$)isomorphic iff there is an element $\zeta\in {\rm
Aut}_{\frac{1}{2}}\Gamma$ such that ${\mathcal
P}_1^{\zeta}={\mathcal P}_2$ or ${\mathcal P}_1^{\zeta}={\mathcal
P}_2^{-1}$. }

\vskip 3mm

{\it Proof} By the definition of equivalence between maps, if
$\kappa$ is an equivalence between $M_1$ and $M_2$, then $\kappa$
is an isomorphism between $M_1$ and $M_2$ induced by an
automorphism $\iota\in {\rm Aut}\Gamma$. Notice that

$${\rm Aut}\Gamma \cong {\rm Aut}\Gamma|^{\frac{1}{2}}
\leq {\rm Aut}_{\frac{1}{2}}\Gamma .$$

\no{Whence, we know that $\iota\in {\rm
Aut}_{\frac{1}{2}}\Gamma$.}

Now if there is a $\zeta\in {\rm Aut}_{\frac{1}{2}}\Gamma$ such
that ${\mathcal P}_1^{\zeta}={\mathcal P}_2$, then $\forall e_x\in
X_{\frac{1}{2}}(\Gamma)$, $\zeta (e_x)=\zeta (e)_{\zeta (x)}$. Now
assume that $e=(x,y)\in E(\Gamma )$, then by our convention, we
know that if $e_x=e\in {\mathcal X}_{\alpha ,\beta}$, then
$e_y=\beta e$. Now by the definition of an automorphism on the
semi-arc set $X_{\frac{1}{2}}(\Gamma)$, if $\zeta (e_x)=f_u$,
where $f=(u,v)$, then there must be $\zeta (e_y)=f_v$. Notice that
$X_{\frac{1}{2}}(\Gamma) = {\cal X}_{\beta}$. We know that $\zeta
(e_y)=\zeta (\beta e)=\beta f=f_v$. We can also extend the action
of $\zeta$ on  $X_{\frac{1}{2}}(\Gamma)$ to ${\mathcal X}_{\alpha
,\beta}$ by $\zeta (\alpha e)=\alpha\zeta (e)$. Whence, we know
that $\forall e\in {\mathcal X}_{\alpha ,\beta}$,

$$\alpha\zeta (e)=\zeta\alpha (e), \ \beta\zeta (e)=\zeta\beta (e) \ {\rm
and} \ {\mathcal P}_1^{\zeta }(e)= {\mathcal P}_2(e).$$

\no{Therefore, the extend action of $\zeta$ on  ${\mathcal
X}_{\alpha ,\beta}$ is an isomorphism between the map $M_1$ and
$M_2$. Whence, $\zeta$ is an equivalence between the map $M_1$ and
$M_2$. So the assertion in ($i$) is true. }

For the assertion in ($ii$), if there is an element $\zeta\in {\rm
Aut}_{\frac{1}{2}}\Gamma$ such that ${\mathcal
P}_1^{\zeta}={\mathcal P}_2$, then the map $M_1$ is isomorphic to
$M_2$. If ${\mathcal P}_1^{\zeta}={\mathcal P}_2^{-1}$, then
${\mathcal P}_1^{\zeta\alpha}={\mathcal P}_2$. The map $M_1$ is
also isomorphic to $M_2$. This is the sufficiency of ($ii$).

By the definition of an isomorphism $\xi$ between maps $M_1$ and
$M_2$, we know that $\forall x\in {\mathcal X}_{\alpha ,\beta}$,

$$\alpha\xi(x)=\xi\alpha (x), \ \beta\xi (x)=\xi\beta (x) \ {\rm
and} \ {\mathcal P}_1^{\xi }(x)= {\mathcal P}_2(x).$$

 By the convention, the condition

$$\beta\xi (x)=\xi\beta (x) \ {\rm and} \ {\mathcal P}_1^{\xi }(x)= {\mathcal P}_2(x).$$

\no{is just the condition of an automorphism $\xi$ or $\alpha\xi$
on $X_{\frac{1}{2}}(\Gamma)$. Whence, the assertion in ($ii$) is
also true.\quad\quad $\natural$}

\vskip 6mm

{\bf $3.2$ A scheme for enumerating maps underlying a graph}

\vskip 4mm

For a given graph $\Gamma$, denoted by ${\mathcal E }^{O}(\Gamma),
{\mathcal E}^{N}(\Gamma)$ and ${\mathcal E}^{L}(\Gamma)$ the sets
of embeddings of $\Gamma$ on the orientable surfaces, on the
non-orientable surfaces and on the locally orientable surfaces,
respectively. For determining the number of non-equivalent
embeddings of a graph on surfaces and non-isomorphic unrooted maps
underlying a graph, another form of the Theorem $1.3.2$ using the
group action idea is need, which is stated as follows.

\vskip 4mm

\no{\bf Theorem $1.3.3$} {\it For two maps $M_1=({\mathcal
X}_{\alpha ,\beta},{\mathcal P}_1)$ and $M_2=({\mathcal X}_{\alpha
,\beta},{\mathcal P}_2)$ underlying a graph $\Gamma$, then

($i$) $M_1,M_2$ are equivalent iff $M_1,M_2$ are in one orbits of
${\rm Aut}_{\frac{1}{2}}\Gamma$ action on
$X_{\frac{1}{2}}(\Gamma)$;

($ii$) $M_1,M_2$ are isomorphic iff $M_1,M_2$ are in one orbits of
${\rm Aut}_{\frac{1}{2}}\Gamma\times <\alpha >$ action on
${\mathcal X}_{\alpha ,\beta}$. }

\vskip 3mm

Now we can established a scheme for enumerating the number of
non-isomorphic unrooted maps and non-equivalent embeddings in a
given set of embeddings of a graph on surfaces by using the {\it
Burnside Lemma} as the following.

\vskip 4mm

\no{\bf Theorem $1.3.4$} {\it For a given graph $\Gamma$, let
${\mathcal E} \subset {\mathcal E}^{L}(\Gamma)$, then the numbers
$n({\mathcal E}, \Gamma)$ and $\eta ({\mathcal E}, \Gamma)$ of
non-isomorphic unrooted maps and non-equivalent embeddings in
${\mathcal E}$ are respective}

$$
n({\mathcal E},\Gamma)=\frac{1}{2|{\rm
Aut}_{\frac{1}{2}}\Gamma|}\sum\limits_{g \in {\rm
Aut_{\frac{1}{2}}\Gamma}} |\Phi_1 (g)|,
$$

\no{\it where, $\Phi_1(g)=\{ {\mathcal P}| {\mathcal P}\in
{\mathcal E}$ and ${\mathcal P}^g= {\mathcal P} \ {\rm or} \
{\mathcal P}^{g\alpha}= {\mathcal P} \}$ and}

$$
\eta ({\mathcal E},\Gamma)=\frac{1}{|{\rm
Aut}_{\frac{1}{2}}\Gamma|}\sum\limits_{g \in {\rm
Aut_{\frac{1}{2}}\Gamma}} |\Phi_2 (g)|,
$$

\no{\it where, $\Phi_2 (g)=\{ {\mathcal P}| {\mathcal P}\in
{\mathcal E}$ and ${\mathcal P}^g= {\mathcal P}\}$.}

\vskip 3mm

{\it Proof} \quad Define the group ${\mathcal H}={\rm
Aut}_{\frac{1}{2}}\Gamma\times <\alpha >$. Then by the Burnside
Lemma and the Theorem $1.3.3$, we get that

$$
n({\mathcal E},\Gamma)=\frac{1}{|\mathcal H|}\sum\limits_{g \in
{\mathcal H }} |\Phi_1 (g)|,
$$

\no{where, $\Phi_1 (g)=\{ {\mathcal P}| {\mathcal P}\in {\mathcal
E}$ and ${\mathcal P}^g= {\mathcal P}  \}$. Now $|{\mathcal H}|=
2|{\rm Aut}_{\frac{1}{2}}\Gamma|$}. Notice that if ${\mathcal
P}^g={\mathcal P}$, then ${\mathcal P}^{g\alpha } \not={\mathcal
P}$, and if ${\mathcal P}^{g\alpha } ={\mathcal P}$, then
${\mathcal P}^g \not={\mathcal P}$. Whence, $\Phi_1
(g)\bigcap\Phi_1 (g\alpha )=\emptyset $. We have that

$$
n({\cal E},\Gamma)=\frac{1}{2|\rm
Aut_{\frac{1}{2}}\Gamma|}\sum\limits_{g \in {\rm
Aut_{\frac{1}{2}}\Gamma}} |\Phi_1 (g)|,
$$

\no{where, $\Phi_1(g)=\{ {\mathcal P}| {\mathcal P}\in {\mathcal
E}$ and ${\mathcal P}^g= {\mathcal P} \ {\rm or} \ {\mathcal
P}^{g\alpha}= {\mathcal P}  \}$.}

A similar proof enables us to obtain that

$$
\eta ({\mathcal E},\Gamma)=\frac{1}{|{\rm
Aut}_{\frac{1}{2}}\Gamma|}\sum\limits_{g \in {\rm
Aut_{\frac{1}{2}}\Gamma}} |\Phi_2 (g)|,
$$

\no{\it where, $\Phi_2 (g)=\{ {\mathcal P}| {\mathcal P}\in
{\mathcal E}$ and ${\mathcal P}^g= {\mathcal P}\}$. \quad\quad
$\natural$}

From the Theorem $1.3.4$, we get the following results.

\vskip 4mm

\no{{\bf Corollary $1.3.1$} {\it The numbers
$n^{O}(\Gamma),n^{N}(\Gamma)$ and $n^{L}(\Gamma)$ of
non-isomorphic unrooted orientable maps ,non-orientable maps and
locally orientable maps underlying a graph $\Gamma$ are
respective}}

$$
n^{O}(\Gamma)=\frac{1}{2|\rm
Aut_{\frac{1}{2}}\Gamma|}\sum\limits_{g \in {\rm
Aut_{\frac{1}{2}}\Gamma}} |\Phi_1^{O} (g)|; \quad\quad (1.3.1)
$$

$$
n^{N}(\Gamma)=\frac{1}{2|\rm
Aut_{\frac{1}{2}}\Gamma|}\sum\limits_{g \in {\rm
Aut_{\frac{1}{2}}\Gamma}} |\Phi_1^{N} (g)|; \quad\quad (1.3.2)
$$

$$
n^{L}(\Gamma)=\frac{1}{2|\rm
Aut_{\frac{1}{2}}\Gamma|}\sum\limits_{g \in {\rm
Aut_{\frac{1}{2}}\Gamma}} |\Phi_1^{L} (g)|, \quad\quad (1.3.3)
$$

\no{\it where, $\Phi_1^{O}(g)=\{ {\mathcal P}| {\mathcal P}\in
{\mathcal E}^{O}(\Gamma)$ and ${\mathcal P}^g= {\mathcal P}$ or
${\mathcal P}^{g\alpha}= {\mathcal P}  \}$, $\Phi_1^{N}(g)=\{
{\mathcal P}| {\mathcal P}\in {\mathcal E}^{N}(\Gamma)$ and
${\mathcal P}^g= {\mathcal P}$ or ${\mathcal P}^{g\alpha}=
{\mathcal P} \}$, $\Phi_1^{L}(g)=\{ {\mathcal P}| {\mathcal P}\in
{\mathcal E}^{L}(\Gamma)$ and ${\mathcal P}^g= {\mathcal P}$ or
${\mathcal P}^{g\alpha}= {\mathcal P} \}$.}

\vskip 4mm

\no{{\bf Corollary $1.3.2$} {\it The numbers
$\eta^{O}(\Gamma),\eta^{N}(\Gamma)$ and $\eta^{L}(\Gamma)$ of
non-equivalent embeddings of a graph $\Gamma$ on orientable
,non-orientable  and locally orientable surfaces are respective}}

$$
\eta^{O}(\Gamma)=\frac{1}{|\rm
Aut_{\frac{1}{2}}\Gamma|}\sum\limits_{g \in {\rm
Aut_{\frac{1}{2}}\Gamma}} |\Phi_2^{O} (g)|; \quad\quad (1.3.4)
$$

$$
\eta^{N}(\Gamma)=\frac{1}{|\rm
Aut_{\frac{1}{2}}\Gamma|}\sum\limits_{g \in {\rm
Aut_{\frac{1}{2}}\Gamma}} |\Phi_2^{N} (g)|; \quad\quad (1.3.5)
$$

$$
\eta^{L}(\Gamma)=\frac{1}{|\rm
Aut_{\frac{1}{2}}\Gamma|}\sum\limits_{g \in {\rm
Aut_{\frac{1}{2}}\Gamma}} |\Phi_2^{L} (g)|, \quad\quad (1.3.6)
$$

\no{\it where, $\Phi_2^{O}(g)=\{ {\mathcal P}| {\mathcal P}\in
{\mathcal E}^{O}(\Gamma)$ and ${\mathcal P}^g= {\mathcal P} \}$,
$\Phi_2^{N}(g)=\{ {\mathcal P}| {\mathcal P}\in {\mathcal
E}^{N}(\Gamma)$ and ${\mathcal P}^g= {\mathcal P} \}$,
$\Phi_2^{L}(g)=\{ {\mathcal P}| {\mathcal P}\in {\mathcal
E}^{L}(\Gamma)$ and ${\mathcal P}^g= {\mathcal P} \}$.}

\vskip 3mm

For a simple graph $\Gamma$, since ${\rm Aut}_{\frac{1}{2}}\Gamma
={\rm Aut}\Gamma$ by the Theorem $1.3.1$, the formula $(1.3.4)$ is
just the scheme used for counting the non-equivalent embeddings of
a complete graph, a complete bipartite graph in the references
$[6],[54]-[55],[70]$. For an {\it asymmetric graph} $\Gamma$, that
is, ${\rm Aut}_{\frac{1}{2}}\Gamma = {\bf id}_{\Gamma}$, we get
the numbers of
 non-isomorphic maps underlying a graph $\Gamma$ and the numbers of
 non-equivalent embeddings of the graph $\Gamma$ by the
Corollary $1.3.1$ and $1.3.2$ as follows.

\vskip 4mm \no{\bf Theorem $1.3.5$}  {\it The numbers $n^O(\Gamma
)$, $n^N(\Gamma )$ and $n^L(\Gamma )$ of non-isomorphic unrooted
maps on orientable, non-orientable surface or locally surface with
an asymmetric underlying graph $\Gamma$ are respective }

$$n^O(\Gamma) = \frac{\prod\limits_{v\in V(\Gamma)} (\rho (v)-1)!}{2},$$

$$n^L( \Gamma) = 2^{\beta (\Gamma)-1}\prod\limits_{v\in V(\Gamma)} (\rho (v)-1)!$$

\no{\it and}

$$n^N(\Gamma) = (2^{\beta (\Gamma)-1}-\frac{1}{2})\prod\limits_{v\in V(\Gamma)}
(\rho (v)-1)!,$$

\no{\it where, $\beta (\Gamma)$ is the Betti number of the graph
$\Gamma$.

The numbers $\eta^O(\Gamma )$, $\eta^N(\Gamma )$ and
$\eta^L(\Gamma )$ of non-equivalent embeddings of an asymmetric
underlying graph $\Gamma$ are respective}

$$\eta^O(\Gamma) = \prod\limits_{v\in V(\Gamma)} (\rho (v)-1)!,$$

$$\eta^L( \Gamma) = 2^{\beta (\Gamma)}\prod\limits_{v\in V(\Gamma)} (\rho (v)-1)!$$

\no{\it and}

$$\eta^N(\Gamma) = (2^{\beta (\Gamma)}-1)\prod\limits_{v\in V(\Gamma)} (\rho (v)-1)!.$$

\vskip 15mm

{\bf\S 4. A relation among the total embeddings and rooted maps of
a graph on genus}

\vskip 8mm

{\bf $4.1$ The rooted total map and embedding polynomial of a
graph}

\vskip 6mm

For a given graph $\Gamma$ with maximum valency $\geq 3$, assume
that $r_i(\Gamma), \widetilde{r_i}(\Gamma), i\geq 0$ are
respectively the numbers of rooted maps with an underlying graph
$\Gamma$ on the orientable surface with genus $\gamma
(\Gamma)+i-1$ and on the non-orientable surface with genus
$\widetilde{\gamma} (\Gamma)+i-1$,  where $\gamma (\Gamma)$ and
$\widetilde{\gamma}(\Gamma)$ denote the minimum orientable genus
and minimum non-orientable genus of the graph $\Gamma$,
respectively. Define its {\it rooted orientable map polynomial}
$r[\Gamma ](x)$ , {\it rooted non-orientable map polynomial}
$\widetilde{r}[\Gamma ](x)$ and {\it rooted total map polynomial}
$R[\Gamma](x)$ on genus by:\vskip 3mm

$$
r[\Gamma] (x) =\sum\limits_{i\geq 0} r_{i}(\Gamma) x^{i},
$$

$$
\widetilde{r}[\Gamma] (x) =\sum\limits_{i\geq 0}
\widetilde{r_{i}}(\Gamma) x^{i}
$$

\no{and}

$$
R[\Gamma](x)=\sum\limits_{i\geq 0} r_{i}(\Gamma)
x^{i}+\sum\limits_{i\geq 1} \widetilde{r_{i}}(\Gamma) x^{-i}.
$$

The total number of orientable embeddings of $\Gamma$ is
$\prod\limits_{d\in D(\Gamma)} (d-1)!$ and non-orientable
embeddings is $(2^{\beta (\Gamma)}-1)\prod\limits_{d\in D(\Gamma)}
(d-1)!$, where $D(\Gamma)$ is its valency sequence. Now let
$g_i(\Gamma)$ and $\widetilde{g_{i}}(\Gamma)$, $i\geq 0$
respectively be the number of embeddings of $\Gamma$ on the
orientable surface with genus $\gamma(\Gamma)+i-1$ and on the
non-orientable surface with genus $\widetilde{\gamma}
(\Gamma)+i-1$. The {\it orientable genus polynomial} $g[\Gamma]
(x) $ , {\it non-orientable genus polynomial}
$\widetilde{g}[\Gamma] (x)$ and {\it total genus polynomial}
$G[\Gamma](x)$ of $\Gamma$ are defined by\vskip 3mm

$$
g[\Gamma] (x) = \sum\limits_{i\geq 0} g_{i}(\Gamma) x^{i},
$$

$$
\widetilde{g}[\Gamma] (x) = \sum\limits_{i\geq 0}
\widetilde{g_{i}}(\Gamma) x^{i}
$$

\no{and}

$$
G[\Gamma](x)=\sum\limits_{i\geq 0} g_{i}(\Gamma)
x^{i}+\sum\limits_{i\geq 1} \widetilde{g_{i}}(\Gamma) x^{-i}.
$$

The orientable genus polynomial $g[\Gamma](x)$ is introduced by
Gross and Furst in $[23]$, and in $[19],[23]-[24]$, the orientable
genus polynomials of a necklace, a bouquet, a closed-end ladder
and a cobblestone are determined. The total genus polynomial is
introduced by Chern et al. in $[13]$, and in $[31]$, recurrence
relations for the total genus polynomial of a bouquet and a dipole
are found. The rooted orientable map polynomial is introduced in
$[43]-[44],[47]$ and the rooted non-orientable map polynomial in
$[48]$. All the polynomials $r[\Gamma] (x),\widetilde{r}[\Gamma]
(x),R[\Gamma](x)$ and $g[\Gamma] (x),\widetilde{g}[\Gamma]
(x),G[\Gamma](x)$ are finite by the properties of embeddings of
$\Gamma$ on surfaces.

Now we establish relations of $r[\Gamma](x)$ with $g[\Gamma](x)$,
$\widetilde{r}[\Gamma](x)$ with $\widetilde{g}[\Gamma](x)$ and
$R[\Gamma](x)$ with $G[\Gamma](x)$ as follows.

\vskip 4mm \no{{\bf Lemma $1.4.1([25][45])$} {\it For a given map
$M$, the number of non-isomorphic rooted maps by rooting on $M$
is}}

$$
\frac{4\varepsilon (M)}{|{\rm Aut}M|},
$$

\no{\it where $\varepsilon (M)$ is the number of edges in $M$.}

\vskip 4mm

\no{{\bf Theorem $1.4.1$ }\quad {\it For a given graph $\Gamma$,}}

$$
|\rm{Aut_{\frac{1}{2}}\Gamma}|r[\Gamma] (x) = 2\varepsilon(\Gamma)
g[\Gamma] (x),
$$

$$
|\rm{Aut_{\frac{1}{2}}\Gamma}|\widetilde{r}[\Gamma] (x) =
2\varepsilon(\Gamma) \widetilde{g}[\Gamma] (x)
$$

\no{and}

$$
|\rm{Aut_{\frac{1}{2}}\Gamma}|R[\Gamma] (x) = 2\varepsilon(\Gamma)
G[\Gamma] (x),
$$

\no{\it where ${\rm Aut}_{\frac{1}{2}}\Gamma$ and $\varepsilon
(\Gamma)$ denote the semi-arc automorphism group and the size of
$\Gamma$, respectively.}

\vskip 3mm

{\it Proof}\quad  For an integer $k$, denotes by ${\cal
M}_{k}(\Gamma,S)$ all the non-isomorphic unrooted maps on an
orientable surface $S$ with genus $\gamma (\Gamma)+k-1$. According
to the Lemma $1.4.1$, we know that

\begin{eqnarray*}
r_{k}(\Gamma) &=& \sum\limits_{M\in{{\cal M}_{k}(\Gamma,S)}}
\frac{4 \varepsilon (M)}{|{\rm Aut}M|}\\
&=& \frac{4\varepsilon (\Gamma)}{|{\rm
Aut_{\frac{1}{2}}}\Gamma\times < \alpha >|} \sum\limits_{M\in
{{\cal M}_{k}(\Gamma,S)}} \frac{|{\rm
Aut_{\frac{1}{2}}}\Gamma\times < \alpha >|}{|{\rm Aut}M|} .
\end{eqnarray*}

Since $|{\rm Aut_{\frac{1}{2}}}\Gamma\times < \alpha >|= |({\rm
Aut_{\frac{1}{2}}}\Gamma\times <\alpha >)_{M}||M^{{\rm
Aut_{\frac{1}{2}}}\Gamma\times <\alpha >}|$ and $|({\rm
Aut_{\frac{1}{2}}}\Gamma\times\\ <\alpha >)_{M}| = |{\rm Aut}M|$,
we have that

\begin{eqnarray*}
r_{k}(\Gamma) &=& \frac{4\varepsilon (\Gamma)}{|{\rm
Aut_{\frac{1}{2}}}\Gamma\times < \alpha >|}
\sum\limits_{M\in {{\cal M}_{k}(\Gamma,S)}}|M^{{\rm Aut_{\frac{1}{2}}}\Gamma\times <\alpha >}|\\
&=& \frac{2\varepsilon (\Gamma)
g_{k}(\Gamma)}{|\rm{Aut_{\frac{1}{2}}\Gamma}|}.
\end{eqnarray*}

\no{Therefore, we get that}

\begin{eqnarray*}
|\rm{Aut_{\frac{1}{2}}\Gamma}|r[\Gamma](x) &=& |\rm{Aut_{\frac{1}{2}}\Gamma}|\sum\limits_{i\geq 0} r_{i}(\Gamma)x^{i}\\
&=& \sum\limits_{i\geq 0} |\rm{Aut_{\frac{1}{2}}\Gamma}| r_{i}(\Gamma)x^{i}\\
&=& \sum\limits_{i\geq 0} 2\varepsilon (\Gamma) g_{i}(\Gamma)x^{i}
= 2\varepsilon (\Gamma) g[\Gamma](x).
\end{eqnarray*}

Similarly, let $\widetilde{{\cal M}_{k}} (\Gamma,\widetilde{S})$
be all the non-isomorphic unrooted maps on an non-orientable
surface $\widetilde{S}$ with genus $\widetilde{\gamma
}(\Gamma)+k-1$. Similar to the proof for orientable case, we can
get that

\begin{eqnarray*}
\widetilde{r_{k}}(\Gamma)&=& \frac{4\varepsilon (\Gamma)}{|{\rm
Aut_{\frac{1}{2}}}\Gamma\times < \alpha >|} \sum\limits_{M\in
{\widetilde{{\cal M}_{k}}(\Gamma,\widetilde{S})}} \frac{|{\rm
Aut_{\frac{1}{2}}}\Gamma\times < \alpha >|}{|{\rm Aut}M|}\\
&=& \frac{4\varepsilon (\Gamma)}{|{\rm
Aut_{\frac{1}{2}}}\Gamma\times < \alpha >|} \sum\limits_{M\in
{\widetilde{{\cal M}_{k}}(\Gamma,\widetilde{S})}}
|M^{{\rm Aut_{\frac{1}{2}}}\Gamma\times <\alpha >}|\\
&=& \frac{2\varepsilon (\Gamma)
\widetilde{g_{k}}(\Gamma)}{|\rm{Aut_{\frac{1}{2}}\Gamma}|}.
\end{eqnarray*}

\no{Therefore, we also get that}

\begin{eqnarray*}
|\rm{Aut_{\frac{1}{2}}\Gamma}|\widetilde{r}[\Gamma](x)
&=& \sum\limits_{i\geq 0} |\rm{Aut_{\frac{1}{2}}\Gamma}| \widetilde{r_{i}}(\Gamma)x^{i}\\
&=& \sum\limits_{i\geq 0} 2\varepsilon (\Gamma)
\widetilde{g_{i}}(\Gamma)x^{i}= 2\varepsilon (\Gamma)
\widetilde{g}[\Gamma](x).
\end{eqnarray*}

\no{Notice that}

$$
R[\Gamma](x)=\sum\limits_{i\geq 0} r_{i}(\Gamma)
x^{i}+\sum\limits_{i\geq 1} \widetilde{r_{i}}(\Gamma) x^{-i}
$$

\no{and}

$$
G[\Gamma](x)=\sum\limits_{i\geq 0} g_{i}(\Gamma)
x^{i}+\sum\limits_{i\geq 1} \widetilde{g_{i}}(\Gamma) x^{-i}.
$$

\no{By the previous discussion, we know that for $k\geq 0$,}

$$
r_{k}(\Gamma)=\frac{2\varepsilon (\Gamma)
g_{k}(\Gamma)}{|\rm{Aut_{\frac{1}{2}}\Gamma}|} \ {\rm and} \
\widetilde{r_{k}}(\Gamma)=\frac{2\varepsilon (\Gamma)
\widetilde{g_{k}}(\Gamma)}{|\rm{Aut_{\frac{1}{2}}\Gamma}|}.
$$

\no{Therefore, we get that}

\begin{eqnarray*}
|\rm{Aut_{\frac{1}{2}}\Gamma}|R[\Gamma](x) &=&
|\rm{Aut_{\frac{1}{2}}\Gamma}|( \sum\limits_{i\geq 0}
r_{i}(\Gamma)
x^{i}+\sum\limits_{i\geq 1} \widetilde{r_{i}}(\Gamma) x^{-i})\\
&=& \sum\limits_{i\geq 0} |\rm{Aut_{\frac{1}{2}}\Gamma}|
r_{i}(\Gamma)x^{i}+\sum\limits_{i\geq 1}
|\rm{Aut_{\frac{1}{2}}\Gamma}|
\widetilde{r_{i}}(\Gamma)x^{-i}\\
&=& \sum\limits_{i\geq 0} 2\varepsilon (\Gamma)
g_{i}(\Gamma)x^{i}+\sum\limits_{i\geq 0} 2\varepsilon (\Gamma)
\widetilde{g_{i}}(\Gamma)x^{-i}= 2\varepsilon
(\Gamma)G[\Gamma](x).
\end{eqnarray*}

\no{This completes the proof. \quad\quad\quad\quad $\natural$}

\vskip 4mm

\no{\bf Corollary $1.4.1$} {\it Let be $\Gamma$ a graph and $s\geq
0$ be an integer. If $r_s(\Gamma)$ and $g_s(\Gamma)$ are the
numbers of rooted maps underlying the graph $\Gamma$ and
embeddings of $\Gamma$ on a locally orientable surface of genus
$s$, respectively, then}

$$|{\rm Aut}_{\frac{1}{2}}\Gamma|r_s(\Gamma) \ = \ 2\varepsilon (\Gamma)g_s(\Gamma).$$

\vskip 6mm

{\bf $4.2$ The number of rooted maps underlying a graph on genus}

\vskip 4mm

The Corollary $1.4.1$ in the previous section can be used to find
the implicit relations among $r[\Gamma] (x)$,
$\widetilde{r}[\Gamma] (x)$ or $R[\Gamma](x)$ if the implicit
relations among $g[\Gamma] (x)$,  $\widetilde{g}[\Gamma] (x)$ or
$G[\Gamma](x)$ are known, and vice via.

Denote the variable vector $(x_1,x_2,\cdots)$ by $\b{x}$,

$$
\b{r}(\Gamma)=(\cdots,\widetilde{r}_2(\Gamma),\widetilde{r}_1(\Gamma),
r_0(\Gamma),r_1(\Gamma),r_2(\Gamma),\cdots),
$$

$$
\b{g}(\Gamma)=(\cdots,\widetilde{g}_2(\Gamma),\widetilde{g}_1(\Gamma),
g_0(\Gamma),g_1(\Gamma),g_2(\Gamma),\cdots).
$$

\no{The $\b{r}(\Gamma)$ and $\b{g}(\Gamma)$ are called the {\it
rooted map sequence} and the {\it embedding sequence} of the graph
$\Gamma$.

Define a function $F(\b{x},\b{y})$ to be {\it $y$-linear} if it
can be represented as the following form}

$$
F(\b{x},\b{y})= f(x_1,x_2,\cdots) \ + \
h(x_1,x_2,\cdots)\sum\limits_{i\in I}y_i \ + \
l(x_1,x_2,\cdots)\sum\limits_{\Lambda\in {\bf\mathcal O}}\Lambda
(\b{y}),
$$

\no{where, $I$ denotes a subset of index and ${\bf\mathcal O}$ a
set of linear operators. Notice that $f(x_1,x_2,\cdots
)=F(\b{x},\b{0})$, where $\b{0}=(0,0,\cdots )$. We have the
following general result.}

\vskip 4mm

\no{\bf Theorem $1.4.2$} {\it Let ${\mathcal G}$ be a graph family
and ${\mathcal H}\subseteq {\mathcal G}$.  If their embedding
sequences $\b{g}(\Gamma), \Gamma\in {\mathcal H},$ satisfies the
equation

$$F_{\mathcal H}(\b{x},\b{g}(\Gamma))=0,\quad\quad\quad (4.1)$$

\no{then its rooted map sequences $\b{r}(\Gamma), \Gamma\in
{\mathcal H}$ satisfies the equation}

$$F_{\mathcal H}(\b{x},\frac{|{\rm Aut}_{\frac{1}{2}}\Gamma|}{2\varepsilon (\Gamma)}\b{r}(\Gamma))=0,$$

\no{and vice via, if the rooted map sequences $\b{r}(\Gamma),
\Gamma\in {\mathcal H}$ satisfies the equation}

$$F_{\mathcal H}(\b{x},\b{r}(\Gamma))=0,\quad\quad\quad (4.2)$$

\no{then its embedding sequences $\b{g}(\Gamma), \Gamma\in
{\mathcal H}$ satisfies the equation}

$$F_{\mathcal H}(\b{x},\frac{2\varepsilon (\Gamma)}{|{\rm Aut}_{\frac{1}{2}}\Gamma|}\b{g}(\Gamma))=0.$$

\no{Even more, assume the function $F(\b{x},\b{y})$ is $y$-linear
and $\frac{2\varepsilon (\Gamma)}{|{\rm
Aut}_{\frac{1}{2}}\Gamma|},\Gamma\in {\mathcal H}$ is a constant.
If the embedding sequences $\b{g}(\Gamma),\Gamma\in {\mathcal H}$
satisfies the equation $(4.1)$, then}

$$F_{\mathcal H}^{\diamond}(\b{x},\b{r}(\Gamma))=0,$$

\no{where $F_{\mathcal
H}^{\diamond}(\b{x},\b{y})=F(\b{x},\b{y})+(\frac{2\varepsilon
(\Gamma)}{|{\rm Aut}_{\frac{1}{2}}\Gamma|}-1)F(\b{x},\b{0})$ and
vice via, if its rooted map sequence $\b{g}(\Gamma),\Gamma\in
{\mathcal H}$ satisfies the equation $(4.2)$, then}

$$F_{\mathcal H}^{\star}(\b{x}, \b{g}(\Gamma))=0.$$

\no{where $F_{\mathcal H}^{\star}=F(\b{x},\b{y})+(\frac{|{\rm
Aut}_{\frac{1}{2}}\Gamma|}{2\varepsilon
(\Gamma)}-1)F(\b{x},\b{0})$.}}

\vskip 3mm

{\it Proof} According to the Corollary $1.4.1$ in this chapter,
for any integer $s\geq o$ and $\Gamma\in {\mathcal H}$, we know
that

$$|{\rm Aut}_{\frac{1}{2}}\Gamma|r_s(\Gamma) \ = \ 2\varepsilon (\Gamma)g_s(\Gamma).$$

\no{Whence,}

$$
r_s(\Gamma) \ = \ \frac{2\varepsilon (\Gamma)}{|{\rm
Aut}_{\frac{1}{2}}\Gamma|} \ g_s(\Gamma)
$$

\no{and}

$$
g_s(\Gamma) \ = \ \frac{|{\rm
Aut}_{\frac{1}{2}}\Gamma|}{2\varepsilon (\Gamma)} r_s(\Gamma).
$$

Therefore, if the embedding sequences $\b{g}(\Gamma),
\Gamma\in{\mathcal H}$ satisfies the equation $(4.1)$, then

$$F_{\mathcal H}(\b{x},\frac{|{\rm Aut}_{\frac{1}{2}}\Gamma|}{2\varepsilon (\Gamma)}\b{r}(\Gamma))=0,$$

\no{and vice via, if the rooted map sequences $\b{r}(\Gamma),
\Gamma\in {\mathcal H}$ satisfies the equation $(4.2)$, then }

$$F_{\mathcal H}(\b{x},\frac{2\varepsilon (\Gamma)}{|{\rm Aut}_{\frac{1}{2}}\Gamma|}\b{g}(\Gamma))=0.$$

Now assume that $F_{\mathcal H}(\b{x},\b{y})$ is a $y$-linear
function and has the following form

$$
F_{\mathcal H}(\b{x},\b{y})= f(x_1,x_2,\cdots) \ + \
h(x_1,x_2,\cdots)\sum\limits_{i\in I}y_i \ + \
l(x_1,x_2,\cdots)\sum\limits_{\Lambda\in {\bf\mathcal O}}\Lambda
(\b{y}),
$$

\no{where ${\bf\mathcal O}$ is a set of linear operators. If
$F_{\mathcal H}(\b{x},\b{g}(\Gamma))=0$, that is}

$$
f(x_1,x_2,\cdots) \ + \ h(x_1,x_2,\cdots)\sum\limits_{i\in
I,\Gamma\in {\mathcal H}}g_i(\Gamma) \ + \
l(x_1,x_2,\cdots)\sum\limits_{\Lambda\in {\bf\mathcal O},\Gamma\in
{\mathcal H}}\Lambda (\b{g}(\Gamma))=0,
$$

\no{we get that}

\begin{eqnarray*}
f(x_1,x_2,\cdots) \ &+& \ h(x_1,x_2,\cdots)\sum\limits_{i\in
I,\Gamma\in {\mathcal H}}\frac{|{\rm
Aut}_{\frac{1}{2}}\Gamma|}{2\varepsilon (\Gamma)} r_i(\Gamma) \\
&+& \ l(x_1,x_2,\cdots)\sum\limits_{\Lambda\in {\bf\mathcal
O},\Gamma\in {\mathcal H}}\Lambda (\frac{|{\rm
Aut}_{\frac{1}{2}}\Gamma|}{2\varepsilon (\Gamma)}\b{r}(\Gamma))=0.
\end{eqnarray*}

Since $\Lambda\in{\bf\mathcal O}$ is a linear operator and
$\frac{2\varepsilon (\Gamma)}{|{\rm
Aut}_{\frac{1}{2}}\Gamma|},\Gamma\in {\mathcal H}$ is a constant,
we also have

\begin{eqnarray*}
f(x_1,x_2,\cdots) \ &+& \ \frac{|{\rm
Aut}_{\frac{1}{2}}\Gamma|}{2\varepsilon
(\Gamma)}h(x_1,x_2,\cdots)\sum\limits_{i\in I,\Gamma\in {\mathcal
H}} r_i(\Gamma) \\
&+& \ \frac{|{\rm Aut}_{\frac{1}{2}}\Gamma|}{2\varepsilon
(\Gamma)}l(x_1,x_2,\cdots)\sum\limits_{\Lambda\in {\bf\mathcal
O},\Gamma\in {\mathcal H}}\Lambda (\b{r}(\Gamma))=0,
\end{eqnarray*}

\no{that is,}

$$
\frac{2\varepsilon (\Gamma)}{|{\rm
Aut}_{\frac{1}{2}}\Gamma|}f(x_1,x_2,\cdots) \ + \
h(x_1,x_2,\cdots)\sum\limits_{i\in I,\Gamma\in {\mathcal H}}
r_i(\Gamma) \ + \ l(x_1,x_2,\cdots)\sum\limits_{\Lambda\in
{\bf\mathcal O},\Gamma\in {\mathcal H}}\Lambda (\b{r}(\Gamma))=0.
$$

\no{Therefore, we get that}

$$F_{\mathcal H}^{\diamond}( \b{x},\b{r}(\Gamma))=0.$$

Similarly, if

$$F_{\mathcal H}(\b{x},\b{r}(\Gamma))=0,$$

\no{then we can also get that}

$$F_{\mathcal H}^{\star}(\b{x}, \b{g}(\Gamma))=0.$$

\no{This completes the proof.\quad\quad\quad $\natural$}

\vskip 4mm

\no{\bf Corollary $1.4.2$} {\it Let ${\mathcal G}$ be a graph
family and ${\mathcal H}\subseteq {\mathcal G}$. If the embedding
sequence $\b{g}(\Gamma)$ of a graph $\Gamma\in {\mathcal G}$
satisfies a recursive relation

$$\sum\limits_{i\in J,\Gamma\in {\mathcal H}}a(i, \Gamma)g_i(\Gamma)=0,$$

\no{where $J$ is the set of index, then the rooted map sequence
$\b{r}(\Gamma)$ satisfies a recursive relation}

$$\sum\limits_{i\in J,\Gamma\in {\mathcal H}}
\frac{a(i, \Gamma)|{\rm Aut}_{\frac{1}{2}}\Gamma|}{2\varepsilon
(\Gamma)}r_i(\Gamma)=0,$$

\no{and vice via.} }

\vskip 3mm

A typical example of the Corollary $1.4.2$ is the graph family
bouquets $B_n, n\geq 1$. Notice that in $[24]$, the following
recursive relation for the number $g_m(n)$ of embeddings of a
bouquet $B_n$ on an orientable surface with genus $m$ for $n\geq
2$ was found.

\begin{eqnarray*}
(n+1)g_m(n) &=& 4(2n-1)(2n-3)(n-1)^2(n-2)g_{m-1}(n-2)\\
&+& 4(2n-1)(n-1)g_m(n-1).
\end{eqnarray*}

\no{and with the boundary conditions}
\vskip 2mm

$g_m(n)=0$ if $m\leq 0$ or $n\leq 0;$

$g_0(0)=g_0(1)=1$ and $g_m(0)=g_m(1)=0$ for $m\geq 0;$

$g_0(2)=4, g_1(2)=2, g_m(2)=0$ for $m\geq 1.$

Since $|{\rm Aut}_{\frac{1}{2}}B_n|=2^nn!$, we get the following
recursive relation for the number $r_m(n)$ of rooted maps on an
orientable surface of genus $m$ underlying a graph $B_n$ by the
Corollary $1.4.2$

\begin{eqnarray*}
(n^2-1)(n-2)r_m(n) &=& (2n-1)(2n-3)(n-1)^2(n-2)r_{m-1}(n-2)\\
&+& 2(2n-1)(n-1)(n-2)r_m(n-1),
\end{eqnarray*}

\no{and with the boundary conditions}

 $r_m(n)=0$ if $m\leq
0$ or $n\leq 0;$

$r_0(0)=r_0(1)=1$ and $r_m(0)=r_m(1)=0$ for $m\geq 0;$

$r_0(2)=2, r_1(2)=1, g_m(2)=0$ for $m\geq 1.$

\vskip 4mm

\no{\bf Corollary $1.4.3$} {\it Let ${\mathcal G}$ be a graph
family and ${\mathcal H}\subseteq {\mathcal G}$. If the embedding
sequences $\b{g}(\Gamma), \Gamma\in {\mathcal G}$ satisfies an
operator equation

$$\sum\limits_{\Lambda\in {\bf\mathcal O},
\Gamma\in {\mathcal H}}\Lambda (\b{g}(\Gamma))=0,$$

\no{where ${\mathcal O}$ denotes a set of linear operators, then
the rooted map sequences $\b{r}(\Gamma), \Gamma\in {\mathcal H}$
satisfies an operator equation}

$$\sum\limits_{\Lambda\in {\bf\mathcal O},
\Gamma\in {\mathcal H}}\Lambda (\frac{|{\rm
Aut}_{\frac{1}{2}}\Gamma|}{2\varepsilon
(\Gamma)}\b{r}(\Gamma))=0,$$

\no{and vice via.} }

\vskip 3mm

Let $\theta =(\theta_1,\theta_2,\cdots,\theta_k)\vdash 2n$, i.e.,
$\sum\limits_{j=1}^k\theta_j=2n$ with positive integers
$\theta_j$. Kwak and Shim introduced three linear operators
$\Gamma ,\Theta$ and $\Delta$ to find the total genus polynomial
of the bouquets $B_n, n\geq 1$ in $[31]$, which are defined as
follows.

Denotes by $z_{\theta}$ and $z_{\theta}^{-1}= 1/z_{\theta}$ the
multivariate monomials $\prod\limits_{i=1}^k z_{\theta_i}$ and $
1/\prod\limits_{i=1}^k z_{\theta_i}$, where $\theta
=(\theta_1,\theta_2,\cdots,\theta_k)\vdash 2n$. Then the  linear
operators $\Gamma ,\Theta$ and $\Delta$ are defined by

$$
\Gamma (z_{\theta}^{\pm
1})=\sum\limits_{j=1}^k\sum\limits_{l=0}^{\theta_j} \theta_j
\{(\frac{z_{1+l}z_{\theta_{j+1-l}}}{z_{\theta_j}})z_{\theta}
\}^{\pm 1}
$$

$$
\Theta (z_{\theta}^{\pm
1})=\sum\limits_{j=1}^k(\theta_j^2+\theta_j)
(\frac{z_{\theta_j+2}z_{\theta}}{z_{\theta_j}})^{-1}
$$

\no{and}

$$
\Delta (z_{\theta}^{\pm 1})=\sum\limits_{1\leq i< j\leq k}
2\theta_i\theta_j[\{(\frac{z_{\theta_j+\theta_i+2}}
{z_{\theta_j}z_{\theta_i}})z_{\theta}\}^{\pm
1}+\{(\frac{z_{\theta_j+\theta_i+2}}
{z_{\theta_j}z_{\theta_i}})z_{\theta}\}^{-1}]
$$

Denotes by $\hat{i}[B_n](z_j)$ the sum of all monomial
$z_{\theta}$ or $1/z_{\theta}$ taken over all embeddings of $B_n$
into an orientable or non-orientable surface, that is

$$
\hat{i}[B_n](z_j)=\sum\limits_{\theta\vdash
2n}i_{\theta}(B_n)z_{\theta}+\sum\limits_{\theta\vdash
2n}\tilde{i}_{\theta}(B_n)z_{\theta}^{-1},
$$

\no{where, $i_{\theta}(B_n)$ and $\tilde{i}_{\theta}(B_n)$ denote
the number of embeddings of $B_n$ into orientable and
non-orientable surface of region type $\theta$. They proved in
$[31]$ that}

$$\hat{i}[B_{n+1}](z_j)=(\Gamma+\Theta+\Delta)\hat{i}[B_n](z_j)
=(\Gamma+\Theta+\Delta)^n(\frac{1}{z_2}+z_1^2).$$

\no{and}

$$G[B_{n+1}](x)=(\Gamma+\Theta+\Delta)^n(\frac{1}{z_2}+z_1^2)|_{z_j=x \ for \ j\geq 1 \
and \ (C*)}.$$

\no{Where, $(C*)$ denotes the condition}\vskip 2mm

{\it $(C*)$: replacing the power $1+n-2i$ of x by $i$ if $i\geq 0$
and $-(1+n+i)$ by $-i$ if $i\leq 0$.}\vskip 2mm

\no{Since}

$$\frac{|{\rm Aut}_{\frac{1}{2}}B_n|}{2\varepsilon
(B_n)}=\frac{2^nn!}{2n}=2^{n-1}(n-1)!$$

\no{and $\Gamma ,\Theta ,\Delta$ are linear, by the Corollary
$1.4.3$ we know that}

\begin{eqnarray*}
R[B_{n+1}](x) &=&
\frac{(\Gamma+\Theta+\Delta)\hat{i}[B_n](z_j)}{2^nn!}|_{z_j=x \
for \ j\geq
1 \ and \ (C*)}\\
&=& \frac{(\Gamma+\Theta+\Delta)^n(\frac{1}{z_2}+z_1^2)}
{\prod\limits_{k=1}^n 2^kk!}|_{z_j=x \ for \ j\geq 1 \ and \
(C*)}.
\end{eqnarray*}

\no{For example, calculation shows that}

$$R[B_1](x)= x+\frac{1}{x};$$

$$R[B_2](x)=2+x+\frac{5}{x}+\frac{4}{x^2};$$

$$R[B_3](x)=\frac{41}{x^3}+\frac{42}{x^2}+\frac{22}{x}+5+10x;$$

\no{and}

$$R[B_4](x)=\frac{488}{x^4}+\frac{690}{x^3}+\frac{304}{x^2}+\frac{93}{x}+14+70x+21x^2.$$

\newpage

%%%%%%%%%%%%%%%Headings%%%%%%%%%%%%%%%%%%%%%%%%%%%%%%%%%%%%%%%%%
\thispagestyle{empty} \pagestyle{myheadings} \topmargin 5mm
\headheight 8mm \headsep 10mm

\markright {\scriptsize Chapter $2$\quad On the Automorphisms of a
Klein Surface and a $s$-Manifold}
%%%%%%%%%%%%%%%%%%%%%%%%%%%%%%%%%%%%%%%%%%%%%%%%%%%%%%%%%%%%%%%%

\vskip 35mm

\no{\bf\large Chapter $2$\quad On the Automorphisms of a Klein
Surface and a $s$-Manifold}

\vskip 10mm

Many papers concerned the automorphisms of a Klein surface, such
as,$[1],[15],$ $[26],[38]$ for a Riemann surface by using Fuchsian
group and $[9]-[10],[21]$ for a Klein surface without boundary by
using $NEC$ groups. Since maps is a natural model for the Klein
surfaces, an even more efficient approach is, perhaps, by using
the combinatorial map theory. Establishing some classical results
again and finding their combinatorial refinement are the central
topics in this chapter.

\vskip 8mm

{\bf\S 1. An algebraic definition of a voltage map}

\vskip 6mm

{\bf $1.1$ Coverings of a map }

\vskip 5mm

For two maps $\widetilde{M}=(\widetilde{{\mathcal X}}_{\alpha
,\beta},\widetilde{{\mathcal P}})$ and $M=({\mathcal X}_{\alpha
,\beta},{\mathcal P})$, call the map $\widetilde{M}$ covering the
map $M$ if there is a mapping $\pi:\widetilde{{\mathcal
X}}_{\alpha ,\beta}\rightarrow {\mathcal X}_{\alpha ,\beta}$ such
that $\forall x\in \widetilde{{\mathcal X}_{\alpha ,\beta}}$,

$$\alpha \pi (x) =\pi\alpha (x), \ \beta\pi (x) =\pi\beta (x) \
and \ \pi \widetilde {\mathcal P}(x)={\mathcal P}\pi(x).$$

\no{The mapping $\pi$ is called a {\it covering mapping}. For
$\forall x\in {\mathcal X}_{\alpha ,\beta}$, define the {\it
quadricell set} $\pi^{-1}(x)$ by}

$$\pi^{-1}(x) = \{\widetilde{x}|\widetilde{x}\in (\widetilde{{\mathcal X}}_{\alpha
,\beta} \ {\rm and} \ \pi (\widetilde{x})=x \}.$$

We have the following result.

\vskip 4mm

\no{\bf Lemma $2.1.1$} {\it Let $\pi:\widetilde{{\mathcal
X}}_{\alpha ,\beta}\rightarrow {\mathcal X}_{\alpha ,\beta}$ be a
covering mapping. Then for any two quadricells $x_1,x_2\in
{\mathcal X}_{\alpha ,\beta}$,

($i$) $|\pi^{-1}(x_1)| = |\pi^{-1}(x_2)|$.

($ii$) If $x_1\not=x_2$, then
$\pi^{-1}(x_1)\bigcap\pi^{-1}(x_2)=\emptyset$.}

\vskip 3mm

{\it Proof} $(i)$ By the definition of a map, for $x_1,x_2\in
{\mathcal X}_{\alpha ,\beta}$, there exists an element $\sigma\in
\Psi_J=<\alpha ,\beta ,{\mathcal P}>$ such that $x_2=\sigma
(x_1)$.

Since $\pi$ is an covering mapping from $\widetilde{M}$ to $M$, it
is commutative with $\alpha ,\beta$ and ${\mathcal P}$. Whence,
$\pi$ is also commutative with $\sigma$. Therefore,

$$\pi^{-1}(x_2)=\pi^{-1}(\sigma (x_1))=\sigma(\pi^{-1}(x_1)).$$

Notice that $\sigma\in\Psi_J$ is an $1-1$ mapping on ${\mathcal
X}_{\alpha ,\beta}$. Hence,  $|\pi^{-1}(x_1)| = |\pi^{-1}(x_2)|$.

($ii$) If $x_1\not=x_2$  and there exists an element $y\in
\pi^{-1}(x_1)\bigcap\pi^{-1}(x_2)$, then there must be $x_1=\pi
(y)=x_2$. Contradicts the assumption. \quad\quad $\natural$

The relation of a covering mapping with an isomorphism is in the
following theorem.

\vskip 4mm

\no{\bf Theorem $2.1.1$} {\it Let $\pi:\widetilde{{\mathcal
X}}_{\alpha ,\beta}\rightarrow {\mathcal X}_{\alpha ,\beta}$ be a
covering mapping. Then $\pi$ is an isomorphism  iff $\pi$ is an
$1-1$ mapping.}

\vskip 3mm

{\it Proof} If $\pi$ is an isomorphism between the maps
$\widetilde{M}=(\widetilde{{\mathcal X}}_{\alpha
,\beta},\widetilde{{\mathcal P}})$ and $M=({\mathcal X}_{\alpha
,\beta},{\mathcal P})$, then it must be an $1-1$ mapping by the
definition, and vice via.\quad\quad $\natural$

A covering mapping $\pi$ from $\widetilde{M}$ to $M$ naturally
induces a mapping $\pi^*$ by the following condition:

$$\forall x\in {\mathcal X}_{\alpha
,\beta},g\in {\rm Aut}\widetilde{M}, \pi^*: g\rightarrow \pi
g\pi^{-1} (x) .$$

We have the following result.

\vskip 4mm

\no{\bf Theorem $2.1.2$} {\it If $\pi:\widetilde{{\mathcal
X}}_{\alpha ,\beta}\rightarrow {\mathcal X}_{\alpha ,\beta}$ is a
covering mapping, then the induced mapping $\pi^*$ is a
homomorphism from ${\rm Aut}\widetilde{M}$ to ${\rm Aut}M$.}

\vskip 3mm

{\it Proof} First, we prove that for $\forall g\in {\rm
Aut}\widetilde{M}$ and $ x\in {\mathcal X}_{\alpha ,\beta}$,
$\pi^*(g)\in {\rm Aut M}$.

Notice that for $\forall g\in {\rm Aut}\widetilde{M}$ and $ x\in
{\mathcal X}_{\alpha ,\beta}$,

$$\pi g\pi^{-1}(x)= \pi (g\pi^{-1}(x))\in {\mathcal X}_{\alpha ,\beta}$$

\no{and $\forall x_1,x_2\in {\mathcal X}_{\alpha ,\beta}$, if
$x_1\not=x_2$, then $\pi g\pi^{-1}(x_1)\not= \pi g\pi^{-1}(x_2)$.
Otherwise, assume that}

$$\pi g\pi^{-1}(x_1) = \pi g\pi^{-1}(x_2) = x_0\in {\mathcal X}_{\alpha ,\beta},$$

\no{then we have that $x_1= \pi g^{-1}\pi^{-1}(x_0)=x_2$.
Contradicts to the assumption.}

By the definition, for $x\in {\mathcal X}_{\alpha ,\beta}$ we get
that

$$\pi^*\alpha (x)=\pi g\pi^{-1}\alpha (x)=\pi g\alpha\pi^{-1} (x)=\pi\alpha g\pi^{-1} (x)
=\alpha\pi g\pi^{-1} (x)=\alpha\pi^* (x),$$

$$\pi^*\beta (x)=\pi g\pi^{-1}\beta (x)=\pi g\beta\pi^{-1} (x)=\pi\beta g\pi^{-1} (x)
=\beta\pi g\pi^{-1} (x)=\beta\pi^* (x).$$

Notice that $\pi (\widetilde {\mathcal P})= {\mathcal P}$. We get
that

$$\pi^*{\mathcal P}(x)=\pi g\pi^{-1}{\mathcal P}(x)
= \pi g \widetilde {\mathcal P}\pi^{-1}(x) = \pi \widetilde
{\mathcal P}g\pi^{-1}(x) = {\mathcal P} \pi g\pi^{-1}(x)=
{\mathcal P}\pi^* (x).$$

Therefore, we get that $\pi g\pi^{-1}\in {\rm Aut}M$, i.e.,
$\pi^*: {\rm Aut}\widetilde{M}\rightarrow {\rm Aut}M$.

Now we prove that $\pi^*$ is a homomorphism from ${\rm
Aut}\widetilde{M}$ to ${\rm Aut}M$. In fact, for $\forall
g_1,g_2\in {\rm Aut}\widetilde{M}$, we have that

$$\pi^*(g_1g_2)=\pi (g_1g_2)\pi^{-1}= (\pi g_1\pi^{-1})(\pi g_2\pi^{-1})=\pi^*(g_1)\pi^*(g_2).$$

\no{Whence,  $\pi^*: {\rm Aut}\widetilde{M}\rightarrow {\rm Aut}M$
is a homomorphism. \quad\quad $\natural$}

\vskip 6mm

{\bf $1.2$ Voltage maps}

\vskip 4mm

For creating a homomorphism between Klein surfaces, voltage maps
are extensively used, which is introduced by Gustin in 1963 and
extensively used by Youngs in 1960s for proving the Heawood map
coloring theorem and generalized by Gross in 1974 ($[22]$). Now it
already become a powerful approach for getting regular maps on a
surface, see $[5],[7],[56]-[57],[65]]$, especially, $[56]-[57]$.
It often appears as an embedded voltage graph in references.
Notice that by using the voltage graph theory, the $2$-factorable
graphs are enumerated in $[51]$. Now we give a purely algebraic
definition for voltage maps, not using geometrical intuition and
establish its theory in this section and the next section again.

\vskip 4mm

\no{\bf Definition $2.1.1$} {\it Let $M=({\mathcal X}_{\alpha
,\beta},{\mathcal P})$ be a map and $G$ a finite group. Call a
pair $(M,\vartheta )$ a voltage map with group $G$ if $\vartheta
:{\mathcal X}_{\alpha ,\beta}\rightarrow G$, satisfying the
following condition:

$(i)$ $\forall x\in {\mathcal X}_{\alpha ,\beta}, \vartheta
(\alpha x)=\vartheta (x), \vartheta (\alpha\beta x)=\vartheta
(\beta x)= \vartheta^{-1}(x);$

$(ii)$ $\forall F=(x,y,\cdots,z)(\beta z,\cdots ,\beta y, \beta
x)\in F(M)$, the face set of $M$, $\vartheta (F)=\vartheta
(x)\vartheta (y)\cdots \vartheta (z)$ and $<\vartheta (F)| F\in
F(u), u\in V(M)> = G$, where, $F(u)$ denotes all the faces
incident with the vertex $u$.}

\vskip 3mm

For a given voltage graph $(M,\vartheta )$, define

$${\mathcal X}_{\alpha^{\vartheta} ,\beta^{\vartheta}}={\mathcal X}_{\alpha ,\beta}\times G$$

$${\mathcal P}^{\vartheta}=
\prod\limits_{(x,y,\cdots,z)(\alpha z,\cdots,\alpha y,\alpha x)\in
V(M)}\prod\limits_{g\in G}(x_g,y_g,\cdots ,z_g)(\alpha z_g,\cdots
,\alpha y_g,\alpha x_g),$$

\no{and}

$$\alpha^{\vartheta} = \alpha$$

$$ \beta^{\vartheta}= \prod\limits_{x\in{\mathcal X}_{\alpha
,\beta},g\in G}(x_g,(\beta x)_{g\vartheta (x)}).$$

\no{where, we use $u_g$ denoting an element $(u,g)\in {\mathcal
X}_{\alpha ,\beta}\times G$. }

It can be shown that $M^{\vartheta}= ({\mathcal
X}_{\alpha^{\vartheta} ,\beta^{\vartheta}},{\mathcal
P}^{\vartheta} )$ also satisfying the conditions of a map with the
same orientation as the map $M$. Hence, we can define the lifting
map of a voltage map as follows.

\vskip 4mm

\no{\bf Definition $2.1.2$} {\it For a voltage map $(M,\vartheta
)$ with group $G$, the map $M^{\vartheta}=({\mathcal X}_{\alpha
,\beta}^{\vartheta},{\mathcal P}^{\vartheta})$ is called its
lifting map.}

For a vertex $v=(C)(\alpha C\alpha^{-1})\in V(M)$, denote by
$\{C\}$ the quadricells in the cycle $C$. The following numerical
result is obvious by the definition of a lifting map.

\vskip 4mm

\no{\bf Lemma $2.1.2$} {\it The numbers of vertices and edges in
the lifting map $M^{\vartheta}$ are respective}

$$\nu (M^{\vartheta})= \nu (M)|G| \ \ and \ \ \varepsilon (M^{\vartheta}) = \varepsilon (M)|G|$$

\vskip 4mm

\no{\bf Lemma $2.1.3$} {\it Let $F=(C^*)(\alpha C^*\alpha^{-1})$
be a face in the map $M$. Then there are $|G|/o(F)$ faces in the
lifting map $M^{\vartheta}$ with length $|F|o(F)$ lifted from the
face $F$, where $o(F)$ denotes the order of $\prod\limits_{x\in
\{C\}}\vartheta (x)$ in the group $G$.}

\vskip 3mm

{\it Proof} Let $F=(u,v\cdots ,w)(\beta w,\cdots,\beta v,\beta u
)$ be a face in the map $M$ and $k$ is the length of $F$. Then, by
the definition, for $\forall g\in G$, the conjugate cycles

\begin{eqnarray*}
(C^*)^{\vartheta}&=& (u_g,v_{g\vartheta (u)},\cdots, u_{g\vartheta
(F) },v_{g\vartheta (F)\vartheta
(u)},\cdots,w_{g\vartheta (F)^2},\cdots , w_{g\vartheta^{o(F)-1} (F)})\\
& {\ }& \beta(u_g,v_{g\vartheta (u)},\cdots, u_{g\vartheta (F)
},v_{g\vartheta (F)\vartheta (u)},\cdots,w_{g\vartheta
(F)^2},\cdots , w_{g\vartheta^{o(F)-1} (F)})^{-1}\beta^{-1}.
\end{eqnarray*}

\no{is a face in $M^{\vartheta}$ with length $k o(F)$. Therefore,
there are $|G|/o(F)$ faces in the lifting map $M^{\vartheta}$.
altogether.\quad\quad $\natural$.}

Therefore,we get that

\vskip 4mm

\no{\bf Theorem $2.1.3$} {\it The Euler characteristic $\chi
(M^{\vartheta})$ of the lifting map $M^{\vartheta}$ of the voltage
map $(M,G)$ is}

$$
\chi (M^{\vartheta})=|G|(\chi (M)+\sum\limits_{m\in {\mathcal
O}(F(M))}(-1+\frac{1}{m})),
$$

\no{\it where ${\mathcal O}(F(M))$ denotes the order $o(F)$ set
of the faces in $M$.}

\vskip 3mm

{\it Proof} According to the Lemma $2.1.2$ and $2.1.3$, the
lifting map $M^{\vartheta}$ has $|G|\nu (M)$ vertices,
$|G|\varepsilon (M)$ edges and $|G|\sum\limits_{m\in {\mathcal
O}(F(M))}\frac{1}{m}$ faces. Therefore, we get that

\begin{eqnarray*}
\chi (M^{\vartheta}) &=& \nu (M^{\vartheta})-\varepsilon
(M^{\vartheta})+\phi (M^{\vartheta})\\
&=& |G|\nu (M)-|G|\varepsilon (M)+ |G|\sum\limits_{m\in {\mathcal
O}(F(M))}\frac{1}{m}\\
&=& |G|(\chi (M)-\phi (M)+\sum\limits_{m\in {\mathcal
O}(F(M))}\frac{1}{m})\\
&=& |G|(\chi (M)+\sum\limits_{m\in {\mathcal
O}(F(M))}(-1+\frac{1}{m})). \quad\quad \natural
\end{eqnarray*}

\vskip 8mm

{\bf \S $2.$ Combinatorial conditions for a group being that of a
map}

\vskip 5mm

Locally characterizing that an automorphism of a voltage map is
that of its lifting is well-done in the references $[40]-[41]$.
Among them, a typical result is the following:\vskip 3mm

 {\it An
automorphism $\zeta$ of a map $M$ with voltage assignment
$\vartheta\rightarrow G$ is an automorphism of its lifting map
$M^{\vartheta}$ if for each face $F$ with $\vartheta (F)={\bf
1_G}$, $\vartheta (\zeta (F))={\bf 1_G}$.}\vskip 2mm

\no{Since the central topic in this chapter is found what a finite
group is an automorphism group of a map, i.e., a global question,
the idea used in the references $[40]-[41]$ are not applicable.
New approach should be used.}

\vskip 4mm

{\bf $2.1$ Combinatorial conditions for an automorphism group of a
map}

\vskip 3mm

First, we characterize an automorphism group of a map.

A permutation group $G$ action on $\Omega$ is called {\it
fixed-free} if $G_x={\bf 1}_G$ for $\forall x\in\Omega$. We have
the following.

\vskip 4mm

\no{\bf Lemma $2.2.1$} {\it Any automorphism group $G$ of a map
$M=({\mathcal X}_{\alpha ,\beta},{\mathcal P})$ is fixed-free on
${\mathcal X}_{\alpha ,\beta}$.}

\vskip 3mm

{\it Proof} For $\forall x\in {\mathcal X}_{\alpha ,\beta}$,since
$G\preceq {\rm Aut}M$, we get that $G_x\preceq ({\rm Aut}M)_x$.
Notice that $({\rm Aut}M)_x={\bf 1}_G$. Whence, we know that
$G_x={\bf 1}_G$, i.e., $G$ is fixed-free. \quad\quad $\natural$

Notice that the automorphism group of a lifting map has a obvious
subgroup, determined by the following lemma.

\vskip 4mm

\no{\bf Lemma $2.2.2$} {\it Let $M^{\vartheta}$ be a lifting map
by the voltage assignment $\vartheta: {\mathcal X}_{\alpha
,\beta}\rightarrow G$. Then $G$ is isomorphic to a fixed-free
subgroup of ${\rm Aut}M^{\vartheta}$ on $V(M^{\vartheta})$.}

\vskip 3mm

{\it Proof} For $\forall g\in G$, we prove that the induced action
$g^*:{\mathcal X}_{\alpha^{\vartheta}
,\beta^{\vartheta}}\rightarrow {\mathcal X}_{\alpha^{\vartheta}
,\beta^{\vartheta}}$ by $g^*: x_h\rightarrow x_{gh}$ is an
automorphism of the map $M^{\vartheta}$.

In fact, $g^*$ is a mapping on ${\mathcal X}_{\alpha^{\vartheta}
,\beta^{\vartheta}}$ and for $\forall x_u\in {\mathcal
X}_{\alpha^{\vartheta} ,\beta^{\vartheta}}$, we get
$g^*:x_{g^{-1}u}\rightarrow x_u$.

Now if for $x_h,y_f\in {\mathcal X}_{\alpha^{\vartheta}
,\beta^{\vartheta}}, x_h\not=y_f$, we have that
$g^*(x_h)=g^*(y_f)$, that is, $x_{gh}=y_{gf}$, by the definition,
we must have that $x=y$ and $gh=gf$, i.e., $h=f$. Whence,
$x_h=y_f$, contradicts to the assumption. Therefore, $g^*$ is
$1-1$ on ${\mathcal X}_{\alpha^{\vartheta} ,\beta^{\vartheta}}$.

We prove that for $x_u\in{\mathcal X}_{\alpha^{\vartheta}
,\beta^{\vartheta}}$, $g^*$ ia commutative with
$\alpha^{\vartheta},\beta^{\vartheta}$ and ${\mathcal
P}^{\vartheta}$. Notice that

$$g^*\alpha^{\vartheta}{x_u}=g^*(\alpha x)_u=(\alpha x)_{gu}=\alpha x_{gu}=\alpha g^*(x_u);$$

$$g^*\beta^{\vartheta}(x_u)=g^*(\beta x)_{u\vartheta (x)}=(\beta x)_{gu\vartheta (x)}
=\beta x_{gu\vartheta
(x)}=\beta^{\vartheta}(x_{gu})=\beta^{\vartheta}g^*(x_u)$$

\no{and}

\begin{eqnarray*}
&g^*{\mathcal P}^{\vartheta}(x_u)& \\
&=&
g^*\prod\limits_{(x,y,\cdots,z)(\alpha z,\cdots,\alpha y,\alpha
x)\in V(M)}\prod\limits_{u\in G}(x_u,y_u,\cdots ,z_u)(\alpha
z_u,\cdots ,\alpha y_u,\alpha x_u)(x_u)\\
&=& g^*y_u=y_{gu}\\
&=& \prod\limits_{(x,y,\cdots,z)(\alpha z,\cdots,\alpha y,\alpha
x)\in V(M)}\prod\limits_{gu\in G}(x_{gu},y_{gu},\cdots
,z_{gu})(\alpha
z_{gu},\cdots ,\alpha y_{gu},\alpha x_{gu})(x_{gu})\\
&=& {\mathcal P}^{\vartheta}(x_{gu})={\mathcal
P}^{\vartheta}g^*(x_u).
\end{eqnarray*}

Therefore, $g^*$ is an automorphism of the lifting map
$M^{\vartheta}$.

To see $g^*$ is fixed-free on $V(M)$, choose $\forall
u=(x_h,y_h,\cdots,z_h)(\alpha z_h,\cdots ,\alpha y_h,\alpha
x_h)\in V(M), h\in G$. If $g^*(u)=u$, i.e.,

$$
(x_{gh},y_{gh},\cdots,z_{gh})(\alpha z_{gh},\cdots ,\alpha
y_{gh},\alpha x_{gh})=(x_h,y_h,\cdots,z_h)(\alpha z_h,\cdots
,\alpha y_h,\alpha x_h).
$$

Assume that $x_{gh}=w_h$, where $w_h\in\{x_h,y_h,\cdots,z_h,
\alpha x_h,\alpha y_h,\cdots ,\alpha z_h\}$. By the definition,
there must be that $x=w$ and $gh=h$. Therefore, $g={\bf 1}_G$,
i.e., $\forall g\in G$, $g^*$ is fixed-free on $V(M)$.

Now define $\tau: g^*\rightarrow g$. Then $\tau$ is an isomorphism
between the action of elements in $G$ on ${\mathcal
X}_{\alpha^{\vartheta} ,\beta^{\vartheta}}$ and the group $G$
itself.\quad\quad $\natural$

According to the Lemma $2.2.1$, given a map $M$ and a group
$G\preceq {\rm Aut}M$, we can define a {\it quotient map}
$M/G=({\mathcal X}_{\alpha ,\beta}/G,{\mathcal P}/G)$ as follows.

$$
{\mathcal X}_{\alpha ,\beta}/G =\{x^G|x\in{\mathcal X}_{\alpha
,\beta}\},
$$

\no{where $x^G$ denotes an orbit of $G$ action on ${\mathcal
X}_{\alpha ,\beta}$ containing $x$ and}

$${\mathcal P}/G=\prod\limits_{(x,y,\cdots ,z)
(\alpha z ,\cdots ,\alpha y, \alpha x)\in V(M)}(x^G,y^G,\cdots
)(\cdots ,\alpha y^G,\alpha x^G),$$

\no{since $G$ action on ${\mathcal X}_{\alpha ,\beta}$ is
fixed-free.}

Notice that the map $M$ may be not a regular covering of its
quotient map $M/G$. We have the following theorem characterizing a
fixed-free automorphism group of a map on $V(M)$.

\vskip 4mm

\no{\bf Theorem $2.2.1$} {\it An finite group $G$ is a fixed-free
automorphism group of a map $M=({\mathcal X}_{\alpha
,\beta},{\mathcal P})$ on $V(M)$ iff there is a voltage map
$(M/G,G)$ with an assignment $\vartheta:{\mathcal X}_{\alpha
,\beta}/G\rightarrow G $ such that $M\cong (M/G)^{\vartheta}$.}

\vskip 3mm

{\it Proof} The necessity of the condition is already proved in
the Lemma $2.2.2$. We only need to prove its sufficiency.

Denote by $\pi: M\rightarrow M/G$ the quotient map from $M$ to
$M/G$. For each element of  $\pi^{-1}(x^G)$, we give it a label.
Choose $x\in\pi^{-1}(x^G)$. Assign its label $l:x\rightarrow
x_{{\bf 1}_G}$, i.e., $l(x)=x_{{\bf 1}_G}$. Since the group $G$
acting on ${\mathcal X}_{\alpha ,\beta}$ is fixed-free, if
$u\in\pi^{-1}(x^G)$ and $u=g(x), g\in G$, we label $u$ with
$l(u)=x_g$. Whence, each element in $\pi^{-1}(x^G)$ is labelled by
a unique element in $G$.

Now we assign voltages on the quotient map $M/G=({\mathcal
X}_{\alpha ,\beta}/G,{\mathcal P}/G)$. If $\beta
x=y,y\in\pi^{-1}(y^G)$ and the label of $y$ is $l(y)=y_h^*, h\in
G$, where, $l(y^*)={\bf 1}_G$, then we assign a voltage $h$ on
$x^G$,i.e., $\vartheta (x^G)=h$. We should prove this kind of
voltage assignment is well-done, which means that we must prove
that for $\forall v\in\pi^{-1}(x^G)$ with $l(v)=j, j\in G$, the
label of $\beta v$ is $l(\beta v)=jh$. In fact, by the previous
labelling approach, we know that the label of $\beta v$ is

\begin{eqnarray*}
l(\beta v)&=& l(\beta gx)=l(g\beta x)\\
&=& l(gy)=l(ghy^*)=gh.
\end{eqnarray*}

Denote by $M^l$ the labelled map $M$ on each element in ${\mathcal
X}_{\alpha ,\beta}$. Whence, $M^l\cong M$. By the previous voltage
assignment, we also know that $M^l$ is a lifting of the quotient
map $M/G$ with the voltage assignment $\vartheta:{\mathcal
X}_{\alpha ,\beta}/G\rightarrow G$. Therefore,

$$M \ \cong \ (M/G)^{\vartheta}.$$

\no{This completes the proof. \quad\quad $\natural$}

According to the Theorem $2.2.1$, we get the following result for
a group to be an automorphism group of a map.

\vskip 4mm

\no{\bf Theorem $2.2.2$} {\it If a group $G, G\preceq {\rm Aut}M$,
is fixed-free on $V(M)$, then}

$$
|G|(\chi (M/G)+\sum\limits_{m\in {\mathcal
O}(F(M/G))}(-1+\frac{1}{m}))=\chi (M).
$$

\vskip 3mm

{\it Proof} By the Theorem $2.2.1$, we know that there is a
voltage assignment $\vartheta$ on the quotient map $M/G$ such that

$$M \ \cong \ (M/G)^{\vartheta}.$$

Applying the Theorem $2.1.3$, we know the Euler characteristic of
the map $M$ is

$$
\chi (M)=|G|(\chi (M/G)+\sum\limits_{m\in {\mathcal
O}(F(M/G))}(-1+\frac{1}{m})). \quad\quad \natural
$$

Theorem $2.2.2$ has some useful corollaries for determining the
automorphism group of a map.

\vskip 4mm

\no{\bf Corollary $2.2.1$} {\it If $M$ is an orientable map of
genus $p$, $G\preceq {\rm Aut}M$ is fixed-free on $V(M)$ and the
quotient map $M/G$ with genus $\gamma$, then }

$$
|G|=\frac{2p-2}{2\gamma-2+\sum\limits_{m\in {\mathcal
O}(F(M/G))}(1-\frac{1}{m}))}.
$$

\no{\it Particularly, if $M/G$ is planar, then}

$$
|G|=\frac{2p-2}{-2+\sum\limits_{m\in {\mathcal
O}(F(M/G))}(1-\frac{1}{m}))}.
$$

\vskip 3mm

\no{\bf Corollary $2.2.2$} {\it If $M$ is a non-orientable map of
genus $q$, $G\preceq {\rm Aut}M$ is fixed-free on $V(M)$ and the
quotient map $M/G$ with genus $\delta$, then }

$$
|G|=\frac{q-2}{\delta -2+\sum\limits_{m\in {\mathcal
O}(F(M/G))}(1-\frac{1}{m}))}.
$$

\no{\it Particularly, if $M/G$ is projective planar, then}

$$
|G|=\frac{q-2}{-1+\sum\limits_{m\in {\mathcal
O}(F(M/G))}(1-\frac{1}{m}))}.
$$

\vskip 3mm

By applying the Theorem $2.2.1$, we can also calculate the Euler
characteristic of the quotient map, which enables us to get the
following result for a group being that of a map.

\vskip 4mm

\no{\bf Theorem $2.2.3$} {\it If a group $G, G\preceq {\rm Aut}M$,
then}

$$\chi (M)+\sum\limits_{g\in G, g\not={\bf 1}_G}(|\Phi_v (g)|+|\Phi_f (g)|)=|G|\chi(M/G),$$

\no{where, $\Phi_v(g)=\{v|v\in V(M), v^g=v\}$ and
$\Phi_f(g)=\{f|f\in F(M), f^g=f\}$, and if $G$ is fixed-free on
$V(M)$, then}

$$\chi (M)+\sum\limits_{g\in G, g\not={\bf 1}_G}|\Phi_f (g)|=|G|\chi(M/G).$$

\vskip 3mm

{\it Proof} By the definition of a quotient map, we know that

$$\phi_v (M/G)=orb_v(G)=\frac{1}{|G|}\sum\limits_{g\in G}|\Phi_v (g)|$$

\no{and}

$$\phi_f (M/G)=orb_f(G)=\frac{1}{|G|}\sum\limits_{g\in G}|\Phi_f (g)|,$$

\no{by applying the Burnside Lemma. Since $G$ is fixed-free on
${\mathcal X}_{\alpha ,\beta}$ by the Lemma $2.1$, we also know
that}

$$\varepsilon (M/G)=\frac{\varepsilon (M)}{|G|}.$$

Applying the Euler-Poincar\'{e} formula for the quotient map
$M/G$, we get that

$$\frac{\sum\limits_{g\in G}|\Phi_v (g)|}{|G|}
-\frac{\varepsilon (M)}{|G|}+\frac{\sum\limits_{g\in G}|\Phi_f
(g)|}{|G|}=\chi (M/G).$$

\no{Whence, we have}

$$\sum\limits_{g\in G}|\Phi_v (g)|-\varepsilon (M)+\sum\limits_{g\in G}|\Phi_f (g)|=|G|\chi (M/G).$$

\no{Notice that $\nu (M)=|\Phi_v ({\bf 1}_G)|$, $\phi (M)=|\Phi_f
({\bf 1}_G)|$ and $\nu (M)-\varepsilon (M)+\phi (M)=\chi (M)$. We
know that}

$$\chi (M)+\sum\limits_{g\in G, g\not={\bf 1}_G}(|\Phi_v (g)|+|\Phi_f (g)|)=|G|\chi(M/G).$$

Now if $G$ is fixed-free on $V(M)$, by the Theorem $2.1$, there is
a voltage assignment $\vartheta$ on the quotient map $M/G$ such
that $M\cong (M/G)^{\vartheta}$. According to the Lemma $2.1.2$,
we know that

$$\nu (M/G)=\frac{\nu (M)}{|G|}.$$

\no{Whence, $\sum\limits_{g\in G}|\Phi_v (g)|=\nu (M)$ and
$\sum\limits_{g\in G, g\not={\bf 1}_G}(|\Phi_v (g)|=0$. Therefore,
we get that}

$$\chi (M)+\sum\limits_{g\in G, g\not={\bf 1}_G}|\Phi_f (g)|=|G|\chi(M/G).
\quad\quad \natural$$

Consider the properties of the group $G$ on $F(M)$, we get the
following interesting results.

\vskip 4mm

\no{\bf Corollary $2.2.3$} {\it If a finite group $G, G\preceq
{\rm Aut}M$ is fixed-free on $V(M)$ and transitive on $F(M)$, for
example, $M$ is regular and $G={\rm Aut}M$, then $M/G$ is an one
face map and}

$$\chi (M)=|G|(\chi (M/G)-1)+\phi (M)$$

\vskip 3mm

Particularly, for an one face map, we know that

\vskip 4mm

\no{\bf Corollary $2.2.4$} {\it For an one face map $M$, if $G$,
$G\preceq {\rm Aut}M$ is fixed-free on $V(M)$, then}

$$\chi (M)-1=|G|(\chi (M/G)-1),$$

\no{\it and $|G|$, especially, $|{\rm Aut}M|$ is an integer factor
of $\chi (M)-1$.}

\vskip 3mm

\no{\bf Remark $2.2.1$} For an one face planar map, i.e., the
plane tree, the only fixed-free automorphism group on its vertices
is the trivial group by the Corollary $2.4$.

\vskip 6mm

{\bf $2.2$ The measures on a map}

\vskip 4mm

On the classical geometry, its central question is to determine
the measures on the objects, such as the distance, angle, area,
volume, curvature, $\dots$. For maps being a combinatorial model
of Klein surfaces, we also wish to introduce various measures on a
map and enlarge its application filed to other branch of
mathematics..

{ $2.2.1$ The angle on a map}

\vskip 3mm

For a map $M=({\mathcal X}_{\alpha ,\beta},{\mathcal P})$, $x\in
{\mathcal X}_{\alpha ,\beta}$, the permutation pair
$\{(x,{\mathcal P}x),(\alpha x,{\mathcal P}^{-1}\alpha x)\}$ is
called an {\it angle} incident with $x$, which is introduced by
Tutte in $[66]$. We prove in this section that any automorphism of
a map is a conformal mapping and affirm the Theorem $1.2.7$ in
Chapter $1$ again.

Define an {\it angle transformation} $\Theta$ of angles of a map
$M=({\mathcal X}_{\alpha ,\beta},{\mathcal P})$ as follows.

$$\Theta =\prod\limits_{x\in {\mathcal X}_{\alpha ,\beta}}(x,{\mathcal P}x).$$

Then we have

\vskip 4mm

\no{\bf Theorem $2.2.4$} {\it Any automorphism of a map
$M=({\mathcal X}_{\alpha ,\beta},{\mathcal P})$ is conformal.}

\vskip 3mm

{\it Proof} By the definition, for $\forall g\in {\rm Aut}M$, we
know that

$$\alpha g=g\alpha, \ \beta g=g\beta \ {\rm and} \ {\mathcal P}g=g{\mathcal P}.$$

Therefore, for $\forall x\in {\mathcal X}_{\alpha ,\beta}$, we
have

$$\Theta g(x)=(g(x),{\mathcal P}g(x))$$

\no{and}

$$g\Theta (x)=g(x,{\mathcal P}x)=(g(x),{\mathcal P}g(x)).$$

Whence, we get that for $\forall x\in {\mathcal X}_{\alpha
,\beta}$, $\Theta g(x)=g\Theta (x).$ Therefore, we get that
$\Theta g=g\Theta$,i.e., $g\Theta g^{-1}=\Theta$.

Since for $\forall x\in {\mathcal X}_{\alpha ,\beta}$, $g\Theta
g^{-1}(x)=(g(x),{\mathcal P}g(x))$ and $\Theta (x)=(x,{\mathcal
P}(x))$, we have that

$$(g(x),{\mathcal P}g(x))=(x,{\mathcal
P}(x)).$$

\no{That is, $g$ is a conformal mapping. \quad\quad $\natural$}

\vskip 4mm

{ $2.2.2$ The non-Euclid area on a map}

\vskip 3mm

For a given voltage map $(M,G)$, its {\it non-Euclid area} $\mu
(M,G)$ is defined by

$$
\mu (M,G)=2\pi(-\chi (M)+\sum\limits_{m\in {\mathcal
O}(F(M))}(-1+\frac{1}{m})).
$$

\no{Particularly, since any map $M$ can be viewed as a voltage map
$(M,{\bf 1}_G)$, we get the non-Euclid area of a map $M$}

$$\mu (M)=\mu (M,{\bf 1}_G)=-2\pi\chi (M).$$

Notice that the area of a map is only dependent on the genus of
the surface. We know the following result.

\vskip 4mm

\no{\bf Theorem $2.2.5$} {\it Two maps on one surface have the
same non-Euclid area.}

\vskip 3mm

By the non-Euclid area, we get the {\it Riemann-Hurwitz formula}
in Klein surface theory for a map in the following result.

\vskip 4mm

\no{\bf Theorem $2.2.6$} {\it If $G\preceq {\rm Aut}M$ is
fixed-free on $V(M)$, then}

$$|G|=\frac{\mu (M)}{\mu (M/G,\vartheta)},$$

\no{\it where $\vartheta$ is constructed in the proof of the
Theorem $2.2.1$.}

\vskip 3mm

{\it Proof} According to the Theorem $2.2.2$, we know that

\begin{eqnarray*}
|G|&=& \frac{-\chi (M)}{-\chi (M)+\sum\limits_{m\in {\mathcal
O}(F(M))}(-1+\frac{1}{m})}\\
&=& \frac{-2\pi\chi (M)}{2\pi(-\chi (M)+\sum\limits_{m\in
{\mathcal O}(F(M))}(-1+\frac{1}{m}))}=\frac{\mu (M)}{\mu
(M/G,\vartheta)}. \quad\quad \natural
\end{eqnarray*}

As an interesting result, we can obtain the same result for the
non-Euclid area of a triangle as the classical differential
geometry.

\vskip 4mm

\no{\bf Theorem $2.2.7([42])$} {\it The non-Euclid area $\mu
(\Delta )$ of a triangle $\Delta$ on a surface ${\mathcal S}$ with
internal angles $\eta ,\theta ,\sigma$ is}

$$\mu (\Delta )=\eta +\theta +\sigma -\pi.$$

\vskip 3mm

{\it Proof} According to the Theorem $1.2.1$ and $2.2.5$, we can
assume there is a triangulation $M$ with internal angles $\eta
,\theta ,\sigma$ on ${\mathcal S}$ and with an equal non-Euclid
area on each triangular disk. Then

\begin{eqnarray*}
\phi (M)\mu (\Delta)&=& \mu (M)=-2\pi\chi (M)\\
&=& -2\pi (\nu (M)-\varepsilon (M)+\phi (M)).
\end{eqnarray*}

Since $M$ is a triangulation, we know that

$$2\varepsilon (M)=3\phi (M).$$

Notice that the sum of all the angles in the triangles on the
surface ${\mathcal S}$ is $2\pi\nu (M)$, we get that

\begin{eqnarray*}
\phi (M)\mu (\Delta)&=& -2\pi (\nu (M)-\varepsilon (M)+\phi
(M))=(2\nu (M)-\phi (M))\pi\\
&=& \sum\limits_{i=1}^{\phi (M)}[(\eta +\theta +\sigma)-\pi]=\phi
(M)(\eta +\theta +\sigma-\pi).
\end{eqnarray*}

Whence, we get that

$$\mu (\Delta )=\eta +\theta +\sigma -\pi.\quad\quad \natural$$

\vskip 10mm

{\bf\S 3. A combinatorial refinement of Huriwtz theorem}

\vskip 8mm

In 1893, Hurwitz obtained a famous result for the
orientation-preserving automorphism group ${\rm Aut}^+{\mathcal
S}$ of a Riemann surface ${\mathcal S}$($[11][18][22]$):\vskip 3mm

{\it For a Riemann surface ${\mathcal S}$ of genus $g({\mathcal
S})\geq 2$, ${\rm Aut}^+{\mathcal S}\leq 84(g({\mathcal
S})-1)$.}\vskip 2mm

We have known that the maps are the combinatorial model for Klein
surfaces, especially, the Riemann surfaces. What is its
combinatorial counterpart? What we can say for the automorphisms
of a map?

For a given graph $\Gamma$, define a {\it graphical property} $P$
to be its a kind of subgraphs, such as, regular subgraphs,
circuits, trees, stars, wheels, $\cdots$. Let $M=({\mathcal
X}_{\alpha ,\beta},{\mathcal P})$ be a map. Call a subset $A$ of
${\mathcal X}_{\alpha ,\beta}$ has the graphical property $P$ if
the underlying graph of $A$ has property $P$. Denote by ${\mathcal
A}(P, M)$ the set of all the $A$ subset with property $P$ in the
map $M$.

For refinement the Huriwtz theorem, we get a general combinatorial
result in the following.

\vskip 4mm

\no{\bf Theorem $2.3.1$}\quad  {\it Let $M=({\mathcal X}_{\alpha
,\beta},{\mathcal P})$ be a map. Then for $\forall G\preceq {\rm
Aut}M$,}

$$[|v^G|| v\in V(M)]\ | \ |G|$$

\no{\it and}

$$|G| \ | \ |A||{\mathcal A}(P,M)|,$$

\no{where£¬$[a,b,\cdots]$ denotes least common multiple of
$a,b,\cdots$.} \vskip 3mm

{\it Proof}\quad According to a well-known result in the
permutation group theory, for $\forall v\in V(M)$, we know
$|G|=|G_v||v^G|$. Therefore, $|v^G|\ | \ |G|$. Whence,

$$[|v^G|| v\in V(M)]\ | \ |G|.$$

The group $G$ is fixed-free action on ${\mathcal X}_{\alpha
,\beta}$, i.e., $\forall x\in {\mathcal X}_{\alpha ,\beta}$, we
have $|G_x|=1$ (see also $[28]$).

Now we consider the action of the automorphism group $G$ on
${\mathcal A}(P,M)$. Notice that if $A\in {\mathcal A}(P, M)$,
then then $\forall g\in G$,$A^g)\in {\mathcal A}(P,M)$, i.e.,
$A^G\subseteq {\mathcal A}(P, M)$. That is, the action of $G$ on
${\mathcal A}(P, M)$ is closed. Whence, we can classify the
elements in ${\mathcal A}(P, M)$ by $G$. For $\forall x,y\in
{\mathcal A}(P,M)$, define $x\sim y$ if and only if there is an
element $g, g\in G$ such that $x^g=y$.

Since $|G_x|=1$, i.e., $|x^G|=|G|$, we know that each orbit of $G$
action on ${\mathcal X}_{\alpha ,\beta}$ has a same length $|G|$.
By the previous discussion, the action of $G$ on ${\mathcal A}(P,
M)$ is closed, therefore, the length of each orbit of $G$ action
on ${\mathcal A}(P, M)$ is also $|G|$. Notice that there are
$|A||{\mathcal A}(P, M)|$ quadricells in ${\mathcal A}(P,M)$. We
get that

$$|G| \ | \ |A||{\mathcal A}(P, M)|.$$

\no{This completes the proof. \quad\quad $\natural$}

Choose property $P$ to be tours with each edge appearing at most
$2$ in the map $M$. Then we get the following results by the
Theorem $2.3.1$.

\vskip 4mm

\no{\bf Corollary $2.3.1$} {\it Let ${\mathcal Tr}_2$ be the set
of tours with each edge appearing 2 times. Then for $\forall
G\preceq {\rm Aut}M$,

$$|G| \ | \ (l|{\mathcal Tr}_2|, \ l=|T|=\frac{|T|}{2}\geq 1, \ T\in {\mathcal Tr}_2,).$$

\no{Let ${\mathcal Tr}_1$ be the set of tours without repeat
edges. Then}

$$|G| \ | \ (2l|{\mathcal Tr}_1|,l=|T|=\frac{|T|}{2}\geq 1, \ T\in {\mathcal Tr}_1,).$$

Particularly, denote by $\phi (i,j)$ the number of faces in $M$
with facial length $i$ and singular edges $j$, then

$$|G| \ | \ ((2i-j)\phi (i,j),i,j\geq 1),$$

\no{where,$(a,b,\cdots)$ denotes the greatest common divisor of
$a,b,\cdots$.}}

\vskip 4mm

\no{\bf Corollary $2.3.2$} {\it Let ${\mathcal T}$ be the set of
trees in the map $M$. Then for $\forall G\preceq {\rm Aut}M$,

$$|G| \ | \ (2lt_l,l\geq 1),$$

\no{where $t_l$ denotes the number of trees with $l$ edges.}}

\vskip 4mm

\no{\bf Corollary $2.3.3$} {\it Let $v_i$ be the number of
vertices with valence $i$. Then for $\forall G\preceq {\rm
Aut}M$,}

$$|G| \ | \ (2iv_i,i\geq 1).$$

Theorem $2.3.1$ is a combinatorial refinement of the Hurwitz
theorem. Applying it, we can get the automorphism group of a map
as follows.

\vskip 4mm

\no{\bf Theorem $2.3.2$}  {\it Let $M$ be an orientable map of
genus $g(M)\geq 2$. Then for $\forall G\preceq {\rm Aut}^+M$,

$$|G| \ \leq  \ 84(g(M)-1)$$

\no{and for $\forall G\preceq {\rm Aut}M$,}

$$|G|\leq 168(g(M)-1).$$}

\vskip 3mm

{\it Proof}  Define the average vertex valence $\overline{\nu
(M)}$ and the average face valence $\overline{\phi (M)}$ of a map
$M$ by

$$
\overline{\nu (M)}=\frac{1}{\nu (M)}\sum\limits_{i\geq 1}i\nu_i ,
$$

$$
\overline{\phi (M)}=\frac{1}{\phi (M)}\sum\limits_{j\geq 1}j\phi_j
,
$$

\no{where,$\nu (M)$,$\phi (M)$,$\phi (M)$ and $\phi_j$ denote the
number of vertices, faces, vertices of valence $i$ and faces of
valence $j$, respectively. }

Then we know that $\overline{\nu (M)}\nu(M)=\overline{\phi
(M)}\phi (M)=2\varepsilon (M)$. Whence, $\nu(M)=\frac{2\varepsilon
(M)}{\overline{\nu (M)}}$ and $\phi (M)=\frac{2\varepsilon
(M)}{\overline{\phi (M)}}$. According to the Euler formula, we
have that

$$\nu (M)-\varepsilon (M)+\phi (M)=2-2g(M),$$

\no{where,$\varepsilon (M),g(M)$ denote the number of edges and
genus of the map $M$. We get that}

$$
\varepsilon (M)=\frac{2(g(M)-1)}{(1-\frac{2}{\overline{\nu
(M)}}-\frac{2}{\overline{\phi (M)}})}.
$$

Choose the integers $k=\lceil\overline{\nu (M)}\rceil$ and
$l=\lceil\overline{\phi (M)}\rceil$. We have that

$$
\varepsilon (M)\leq \frac{2(g(M)-1)}{(1-\frac{2}{k}-\frac{2}{l})}.
$$

Because $1-\frac{2}{k}-\frac{2}{l} \
> \ 0$, So $k\geq 3, l \ > \ \frac{2k}{k-2}$. Calculation shows that
the minimum value of $1-\frac{2}{k}-\frac{2}{l}$ is $\frac{1}{21}$
and attains the minimum value if and only if $(k,l)=(3,7)$ or
$(7,3)$. Therefore,

$$
\varepsilon (M\leq 42(g(M)-1)).
$$

According to the Theorem $2.3.1$ and its corollaries, we know that
$|G|\leq 4\varepsilon (M)$ and if $G$ is orientation-preserving,
then $|G|\leq 2\varepsilon (M)$. Whence,

$$
|G| \ \leq \ 168(g(M)-1))
$$

\no{and if $G$ is orientation-preserving, then}

$$
|G| \  \leq  \ 84(g(M)-1)),
$$

\no{with equality hold if and only if $G={\rm AutM},(k,l)=(3,7)$
or $(7,3)$¡£\quad\quad $\natural$}

For the automorphism of a Riemann surface, we have

\vskip 4mm

\no{\bf Corollary $2.3.4$} {\it For any Riemann surface ${\mathcal
S}$ of genus $g\geq 2$,}

$$4g({\mathcal S})+2\leq |{\rm Aut}^+ {\mathcal S}|\leq 84(g({\mathcal S})-1)$$

\no{\it and}

$$8g({\mathcal S})+4\leq |{\rm Aut}{\mathcal S}|\leq 168(g({\mathcal S})-1).$$

\vskip 3mm

{\it Proof} By the Theorem $1.2.6$ and $2.3.2$, we know the upper
bound for $|{\rm Aut}{\mathcal S}|$ and $|{\rm Aut}^+ {\mathcal
S}|$. Now we prove the lower bound. We construct a regular map
$M_k=({\mathcal X}_k,{\mathcal P}_k)$ on a Riemann surface of
genus $g\geq 2$ as follows, where $k=2g+1$.

$${\mathcal X}_k=\{x_1,x_2,\cdots,x_k,\alpha x_1,\alpha x_2,\cdots,\alpha x_k,
\beta x_1,\beta x_2,\cdots,\beta x_k,\alpha\beta x_1,\alpha\beta
x_2,\cdots,\alpha\beta x_k\}$$

$$
{\mathcal P}_k= (x_1,x_2,\cdots,x_k,\alpha\beta x_1,\alpha\beta
x_2,\cdots,\alpha\beta x_k)(\beta x_k,\cdots, \beta x_2,\beta
x_1,\alpha x_k,\cdots,\alpha x_2,\alpha x_1).
$$

It can be shown that $M_k$ is a regular map, and its
orientation-preserving automorphism group ${\rm Aut}^+M_k=
<{\mathcal P}_k>$. Direct calculation shows that if $k\equiv 0(mod
2)$, $M_k$ has $2$ faces, and if $k\equiv 1$, $M_k$ is an one face
map. Therefore, according to the Theorem $1.2.6$, we get that

$$|{\rm Aut}^+ {\mathcal S}|\geq 2\varepsilon (M_k)\geq 4g+2 ,$$

\no{and}

$$|{\rm Aut}{\mathcal S}|\geq 4\varepsilon (M_k)\geq 8g+4. \quad\quad \natural$$

For the non-orientable case, we can also get the bound for the
automorphism group of a map.

\vskip 4mm

\no{\bf Theorem $2.3.3$} {\it Let $M$ be a non-orientable map of
genus $g'(M)\geq 3$. Then for $\forall G\preceq {\rm Aut}^+M$,

$$|G| \ \leq  \ 42(g'(M)-2)$$

\no{and for $\forall G\preceq {\rm Aut}M$,}

$$|G|\leq 84(g'(M)-2),$$

\no{with the equality hold iff $M$ is a regular map with vertex
valence $3$ and face valence $7$ or vice via.}}

\vskip 3mm

{\it Proof} Similar to the proof of the Theorem $2.3.2$, we can
also get that

$$
\varepsilon (M\leq 21(g'(M)-2))
$$

\no{and with equality hold if and only if $G={\rm AutM}$ and $M$
is a regular map with vertex valence $3$, face valence $7$ or vice
via. According to the Corollary $2.3.3$, we get that}

$$|G|\leq 4\varepsilon (M)$$

\no{and if G is orientation-preserving, then}

$$|G|\leq 2\varepsilon (M).$$

Whence, for $\forall G\preceq {\rm Aut}^+M$,

$$|G| \ \leq  \ 42(g'(M)-2)$$

\no{and for $\forall G\preceq {\rm Aut}M$,}

$$|G|\leq 84(g'(M)-2),$$

\no{with the equality hold iff $M$ is a regular map with vertex
valence $3$ and face valence $7$ or vice via. \quad\quad
$\natural$}

Similar to the Hurwtiz theorem for a Riemann surface, we can get
the upper bound for a Klein surface underlying a non-orientable
surface.

\vskip 4mm

\no{\bf Corollary $2.3.5$}  {\it For any Klein surface ${\mathcal
K}$ underlying a non-orientable surface of genus $q\geq 3$,}

$$|{\rm Aut}^+ {\mathcal K}|\leq 42(q-2)$$

\no{\it and}

$$|{\rm Aut}{\mathcal K}|\leq 84(q-2).$$

According to the  Theorem $1.2.8$, similar to the proof of the
Theorem $2.3.2$ and $2.3.3$, we get the following result for the
automorphisms of an $s$-manifold as follows.

\vskip 4mm

\no{\bf Theorem $2.3.4$} {\it Let ${\mathcal C}(T,n)$ be a closed
$s$-manifold with negative Euler characteristic. Then $|{\rm
Aut}{\mathcal C}(T,n)|\leq 6n$ and

$$|{\rm Aut}{\mathcal C}(T,n)|\leq -21\chi ({\mathcal C}(T,n)),$$

\no{with equality hold only if ${\mathcal C}(T,n)$ is hyperbolic,
where $\chi ({\mathcal C}(T,n))$ denotes the genus of ${\mathcal
C}(T,n)$. }}

\vskip 3mm

{\it Proof} \ The inequality $|{\rm Aut}{\mathcal C}(T,n)|\leq 6n$
is known by the Corollary $2.3.1$. Similar to the proof of the
Theorem $2.3.2$, we know that

$$\varepsilon ({\mathcal C}(T,n) )= \frac{-\chi ({\mathcal C}(T,n))}
{\frac{1}{3}-\frac{2}{k}},$$

\no{where $k=\frac{1}{\nu ({\mathcal C}(T,n) )}\sum\limits_{i\geq
1}i\nu_i \leq 7$ and with the equality holds only if $k=7$, i.e.,
${\mathcal C}(T,n)$ is hyperbolic. \quad\quad $\natural$}

\vskip 10mm

{\bf\S 4. The order of an automorphism of a Klein surface}

\vskip 8mm

Harvey $[26]$ in 1966, Singerman $[60]$ in 1971 and Bujalance
$[9]$ in 1983 considered the order of an automorphism of a Riemann
surface of genus $p\geq 2$ and a compact non-orientable Klein
surface without boundary of genus $q\geq 3$. Their approach is by
using the Fuchsian groups and $NEC$ groups for Klein surfaces. The
central idea is by applying the {\it Riemann-Hurwitz equation},
stated as follows:

{\it Let $G$ be an $NEC$ graph and $G'$ be a subgroup of $G$ with
finite index. Then}

$$
\frac{\mu (G')}{\mu (G)}=[G:G'],
$$

\no{\it where, $\mu (G)$ denotes the non-Euclid area of the group
$G$, which is defined  as if

$$\sigma =(g;\pm ;[m_1,\cdots ,m_r];\{(n_{11,\cdots
,n_{1s_1}}),\cdots,(n_{k1},\cdots ,n_{ks})\})$$

\no{is the signature of the group $G$, then }}

$$
\mu (G)=2\pi [\eta
g+k-2+\sum\limits_{i=1}^r(1-1/m_i)+
1/2\sum\limits_{i=1}^k\sum\limits_{j=1}^{s_i}(1-1/n_{ij})],
$$

\no{\it where, $\eta =2$ if $sign (\sigma)=+$ and $\eta =1$
otherwise.}

Notice that we have introduced the conception of non-Euclid area
for the voltage maps and have gotten the Riemann-Hurwitz equation
in the Theorem $2.2.6$ for a fixed-free on $V(M)$ group.
Similarly, we can find the minimum genus of a map, fixed-free on
its vertex set by the voltage assignment on its quotient map and
the maximum order of an automorphism of a map.

\vskip 6mm {\bf $4.1$ The minimum genus of a fixed-free
automorphism}

\vskip 4mm

\no{\bf Lemma $2.4.1$} {Let $N=\prod\limits_{i=1}^kp_i^{r_i},
p_1<p_2<\cdots <p_k$ be the arithmetic decomposition of the
integer $N$ and $m_i\geq 1, m_i|N$ for $i=1,2,\cdots ,k$. Then for
any integer $s\geq1$,}

$$\sum\limits_{i=1}^s(1-\frac{1}{m_i})\geq 2(1-\frac{1}{p_1})\lfloor\frac{s}{2}\rfloor.$$

\vskip 3mm

{\it Proof} If $s\equiv 0(mod 2)$, it is obvious that

$$\sum\limits_{i=1}^s(1-\frac{1}{m_i})\geq\sum\limits_{i=1}^s(1-\frac{1}{p_1}) \geq (1-\frac{1}{p_1})s.$$

Now assume that $s\equiv 1(mod 2)$ and there are $m_{i_j}\not=
p_1, j=1,2,\cdots ,l$. If the assertion is not true, we must have
that

$$(1-\frac{1}{p_1})(l-1) \ > \sum\limits_{j=1}^l(1-\frac{1}{m_{i_j}})\geq (1-\frac{1}{p_2})l.$$

\no{Whence, we get that}

$$(1-\frac{1}{p_1})l \ > (1-\frac{1}{p_2})l+1-\frac{1}{p_1} \ >(1-\frac{1}{p_1})l.$$

\no{A contradiction. Therefore, we have that}

$$\sum\limits_{i=1}^s(1-\frac{1}{m_i})\geq 2(1-\frac{1}{p_1})\lfloor\frac{s}{2}\rfloor. \quad\quad \natural$$

\vskip 4mm

\no{\bf Lemma $2.4.2$} {\it For a map $M=({\mathcal X}_{\alpha
,\beta},{\mathcal P})$ with $\phi (M)$ faces and
$N=\prod\limits_{i=1}^kp_i^{r_i}, p_1<p_2<\cdots <p_k$, the
arithmetic decomposition of the integer $N$, there exists a
voltage assignment $\vartheta :{\mathcal X}_{\alpha
,\beta}\rightarrow Z_N$ such that for $\forall F\in F(M)$,
$o(F)=p_1$ if $\phi (M)\equiv 0(mod 2)$ or there exists a face
$F_0\in F(M)$, $o(F)=p_1$ for $\forall F\in F(M)\setminus
\{F_0\}$, but $o(F_0)=1$. }

\vskip 3mm

{\it Proof} Assume that $f_1,f_2,\cdots ,f_n,$ where,$ n=\phi
(M)$, are the $n$ faces of the map $M$. By the definition of
voltage assignment, if $x,\beta x$ or $x, \alpha\beta x$ appear on
one face $f_i, 1\leq i\leq n$ altogether, then they contribute to
$\vartheta (f_i)$ only with $\vartheta (x)\vartheta^{-1}(x)={\bf
1}_{Z_N}$. Whence, not loss of generality, we only need to
consider the voltage $x_{ij}$ on the common boundary among the
face $f_i$ and $f_j$ for $1\leq i,j\leq n$. Then the voltage
assignment on the $n$ faces are

$$\vartheta (f_1)=x_{12}x_{13}\cdots x_{1n},$$

$$\vartheta (f_2)=x_{21}x_{23}\cdots x_{2n},$$

$$\cdots\cdots\cdots\cdots\cdots\cdots$$

$$\vartheta (f_n)=x_{n1}x_{n2}\cdots x_{n(n-1)}.$$

We wish to find an assignment on $M$ which can enables us to get
as many faces as possible with the voltage of order $p_1$. Not
loss of generality, we can choose $\vartheta^{p_1}(f_1)={\bf
1}_{Z_N}$ in the first. To make $\vartheta^{p_1}(f_2)={\bf
1}_{Z_N}$, choose $x_{23}=x_{13}^{-1},\cdots ,x_{2n}=x_{1n}^{-1}$.
If we have gotten $\vartheta^{p_1}(f_i)={\bf 1}_{Z_N}$ and $i < n$
if $n\equiv 0(mod 2)$ or $i \ <n-1$ if $n\equiv 1(mod 2)$, we can
choose that

$$x_{(i+1)(i+2)}=x_{i(i+2)}^{-1},x_{(i+1)(i+3)}=x_{i(i+3)}^{-1},\cdots
,x_{(i+1)n}=x_{in}^{-1},$$

\no{which also make $\vartheta^{p_1}(f_{i+1})={\bf 1}_{Z_N}$.}

Now if $n\equiv 0(mod 2)$, this voltage assignment makes each face
$f_i, 1\leq i\leq n$ satisfying that $\vartheta^{p_1}(f_i)={\bf
1}_{Z_N}$. But if $n\equiv 1(mod 2)$, it only makes
$\vartheta^{p_1}(f_i)={\bf 1}_{Z_N}$ for $1\leq i\leq n-1$, but
$\vartheta (f_n)={\bf 1}_{Z_N}$. This completes the proof.
\quad\quad $\natural$

Now we can prove a result for the minimum genus of a fixed-free
automorphism of a map.

\vskip 4mm

\no{\bf Theorem $2.4.1$} {\it Let $M=({\mathcal X}_{\alpha
,\beta},{\mathcal P})$ be a map and $N=p_1^{r_1}\cdots p_k^{r_k},
p_1 \ < p_2 \ <\cdots \ < p_k$, be the arithmetic decomposition of
the integer $N$. Then for any voltage assignment $\vartheta
:{\mathcal X}_{\alpha ,\beta}\rightarrow Z_N$,

$(i)$ if $M$ is orientable, the minimum genus $g_{min}$ of the
lifting map $M^{\vartheta}$ which admits an automorphism of order
$N$, fixed-free on $V(M^{\vartheta})$, is

$$g_{min}=1+N\{g(M)-1+(1-\sum\limits_{m\in {{\mathcal O}(F(M))}}\frac{1}{p_1})
\lfloor\frac{\phi (M)}{2}\rfloor\}.$$

$(ii)$ if $M$ is non-orientable, the minimum genus $g'_{min}$ of
the lifting map $M^{\vartheta}$ which admits an automorphism of
order $N$, fixed-free on $V(M^{\vartheta})$, is}

$$g'_{min}= 2+N\{g(M)-2+2(1-\frac{1}{p_1})\lfloor\frac{\phi (M)}{2}\rfloor\}.\quad\quad \natural$$

\vskip 3mm

{\it Proof} $(i)$ According to the Theorem $2.2.1$, we know that

$$2-2g(M^{\vartheta})=N\{(2-2g(M))+\sum\limits_{m\in {{\mathcal O}(F(M))}}(-1+\frac{1}{m})\}.$$

\no{Whence,}

$$2g(M^{\vartheta})=2+N\{2g(M)-2+\sum\limits_{m\in {{\mathcal O}(F(M))}}(1-\frac{1}{m})\}.$$

Applying the Lemma $2.4.1$ and $2.4.2$, we get that

$$g_{min}=1+N\{g(M)-1+(1-\frac{1}{p_1})\lfloor\frac{\phi
(M)}{2}\rfloor\}$$.

($ii$) Similarly, by the Theorem $2.2.1$, we know that

$$2-g(M^{\vartheta})=N\{(2-g(M))+\sum\limits_{m\in {{\mathcal O}(F(M))}}(-1+\frac{1}{m})\}.$$

\no{Whence,}

$$g(M^{\vartheta})=2+N\{g(M)-2+\sum\limits_{m\in {{\mathcal O}(F(M))}}(1-\frac{1}{m})\}.$$

Applying the Lemma $2.4.1$ and $2.4.2$, we get that

$$g'_{min}=2+N\{g(M)-2+2(1-\frac{1}{p_1})\lfloor\frac{\phi (M)}{2}\rfloor\}.\quad\quad \natural$$

\vskip 6mm

{\bf $4.2$ The maximum order of an automorphism of a map}

\vskip 4mm

For the maximum order of an automorphism of a map, we have the
following result.

\vskip 4mm

\no{\bf Theorem $2.4.2$} {\it The maximum order $N_{max}$ of an
 automorphism $g$ of an orientable map $M$ of genus$\geq 2$ is

$$N_{max} \ \leq 2g(M)+1$$

\no{and the maximum order $N'_{max}$ of a non-orientable map of
genus$\geq 3$ is}

$$N'_{max} \ \leq g(M)+1,$$

\no{where $g(M)$ is the genus of the map $M$.}}

\vskip 3mm

{\it Proof} According to the Theorem $2.2.3$, denote $G=<g>$, we
get that

$$\chi (M)+\sum\limits_{g\in G, g\not={\bf 1}_G}(|\Phi_v (g)|+|\Phi_f (g)|)=|G|\chi(M/G),$$

\no{where, $\Phi_f(g)=\{F| F\in F(M), F^g=F\}$ and $\Phi_v(g)=\{
v|v\in V(M), v^g=v\}$.  If $g\not={\bf 1}_{G}$, direct calculation
shows that $\Phi_f(g )=\Phi_f(g^2 )$ and $\Phi_v(g )=\Phi_v(g^2
)$. Whence,}

$$\sum\limits_{g\in G, g\not={\bf 1}_G}|\Phi_v (g)|=(|G|-1)|\Phi_v (g)|$$

\no{and}

$$\sum\limits_{g\in G, g\not={\bf 1}_G}|\Phi_f (g)|=(|G|-1)|\Phi_f (g)|.$$

\no{Therefore, we get that}

$$\chi (M)+(|G|-1)|\Phi_v (g)|+(|G|-1)|\Phi_f (g)|=|G|\chi(M/G)$$

\no{Whence, we have that}

$$\chi (M)-(|\Phi_v (g)|+|\Phi_f (g)|)=|G|(\chi(M/G)-(|\Phi_v (g)|+|\Phi_f (g)|)).$$

If $\chi(M/G)-(|\Phi_v (g)|+|\Phi_f (g)|)=0$,i.e.,
$\chi(M/G)=|\Phi_v (g)|+|\Phi_f (g)|\geq 0$, then we get that
$g(M) \ \leq 1$ if $M$ is orientable or $g(M) \ \leq 2$ if $M$ is
non-orientable. Contradicts to the assumption. Therefore,
$\chi(M/G)-(|\Phi_v (g)|+|\Phi_f (g)|)\not=0$. Whence, we get that

$$
|G|=\frac{\chi (M)-(|\Phi_v (g)|+|\Phi_f (g)|)}{\chi(M/G)-(|\Phi_v
(g)|+|\Phi_f (g)|)}=H(v,f;g).
$$

\no{Notice that $|G|, \chi (M)-(|\Phi_v (g)|+|\Phi_f (g)|)$ and
$\chi(M/G)-(|\Phi_v (g)|+|\Phi_f (g)|)$ are integers. We know that
the function $H(v,f;g)$ takes its maximum value at
$\chi(M/G)-(|\Phi_v (g)|+|\Phi_f (g)|)=-1$ since $\chi (M) \ \leq
-1$. But in this case, we get that }

$$|G|=|\Phi_v
(g)|+|\Phi_f (g)|-\chi (M)=1+\chi (M/G)-\chi (M).$$

We divide the discussion into to cases.\vskip 3mm

{\bf Case $1$}\quad\quad $M$ is orientable.\vskip 2mm

Since $\chi(M/G)+1=(|\Phi_v (g)|+|\Phi_f (g)|) \geq 0$, we know
that $\chi(M/G) \geq -1$. Whence, $\chi (M/G)= 0$ or $2$.
Therefore, we have that

$$|G|=1+\chi (M/G)-\chi (M)\leq 3-\chi (M)=2g(M)+1.$$

\no{That is, $N_{max} \ \leq 2g(M)+1$.}\vskip 3mm

{\bf Case $2$} \quad\quad $M$ is non-orientable.\vskip 2mm

In this case, since $\chi(M/G) \geq -1$, we know that $\chi (M/G)=
-1, 0,1$ or $2$. Whence, we have that

$$|G|=1+\chi (M/G)-\chi (M)\leq 3-\chi (M)=g(M)+1.$$

\no{This completes the proof. \quad\quad $\natural$}

According to this theorem, we get the following result for the
order of an automorphism of a Klein surface without boundary by
the Theorem $1.2.7$, which is even more better than the results
already known.

\vskip 4mm

\no{\bf Corollary $2.4.1$} {\it The maximum order of an
automorphism of a Riemann surface of genus$\geq 2$ is $2g(M)+1$,
and the maximum order of an automorphism of a non-orientable Klein
surface of genus$\geq 3$ without boundary is $g(M)+1$.}

\vskip 3mm

The maximum order of an automorphism of a map can be also
determined by its underlying graph, which is stated as follows.

\vskip 4mm

\no{\bf Theorem $2.4.3$} {\it Let $M$ be a map underlying the
graph $G$ and $o_{max}(M,g), o_{max}(G,g)$ be the maximum order of
orientation-preserving automorphism in ${\rm Aut}M$ and in ${\rm
Aut}_{\frac{1}{2}}G$. Then }

$$o_{max}(M,g) \ \leq \ o_{max}(G,g),$$

\no{and the equality hold for at least one map underlying the
graph $G$.}

\vskip 3mm

The proof of the Theorem $2.4.3$ will be delayed to the next
chapter after we prove the Theorem $3.1.1$. By this theorem, we
get the following interesting results.

\vskip 4mm

\no{\bf Corollary $2.4.2$} {\it The maximum order of an
orientation-preserving automorphism of a complete map ${\mathcal
K}_n, n\geq 3$, is at most $n$.}

\vskip 4mm

\no{\bf Corollary $2.4.3$} {\it The maximum order of an
orientation-preserving automorphism of a plane tree ${\mathcal T}$
is at most $|{\mathcal T}|-1$ and attains the upper bound only if
the underlying tree is the star.}

\newpage

%%%%%%%%%%%%%%%Headings%%%%%%%%%%%%%%%%%%%%%%%%%%%%%%%%%%%%%%%%%
\thispagestyle{empty} \pagestyle{myheadings} \topmargin 5mm
\headheight 8mm \headsep 10mm

\markright {\scriptsize Chapter $3$\quad On the Automorphisms of a
Graph on Surfaces}
%%%%%%%%%%%%%%%%%%%%%%%%%%%%%%%%%%%%%%%%%%%%%%%%%%%%%%%%%%%%%%%%

\vskip 35mm

\no{\bf\large Chapter $3$\quad On the Automorphisms of a Graph on
Surfaces}

\vskip 15mm

For determining the automorphisms of a map, an alternate approach
is to consider the action of the semi-arc automorphism group of
its underlying graph on the quadricells and to distinguish which
is an automorphism of the map and which is not. This approach is
first appeared in the reference $[43]$ as an initial step for the
enumeration of the non-equivalent embeddings of a graph on
surfaces, and also important for enumeration unrooted maps
underlying a graph on surfaces used in Chapter $4$.

\vskip 10mm

{\bf\S 1. A necessary and sufficient condition for a group of a
graph being that of a map}

\vskip 8mm

Let $\Gamma =(V,E)$ be a connected graph. Its automorphism is
denoted by $Aut\Gamma$. Choose the base set $X=E(\Gamma )$. Then
its quadricells ${\mathcal X}_{\alpha ,\beta }$ is defined to be:

$${\mathcal X}_{\alpha ,\beta }=\bigcup\limits_{x\in X}\{x,\alpha x,
\beta x,\beta\alpha\beta x\},$$

\no{where, $K=\{1,\alpha ,\beta ,\alpha\beta\}$ is the Klein $4$-
elements group.}

For $\forall g\in{\rm Aut}\Gamma$, define an {\it induced action}
$g|^{{\mathcal X}_{\alpha,\beta}}$ of $g$ on ${\mathcal
X}_{\alpha,\beta}$ as follows.

\vskip 2mm

{\it For $\forall x\in {\mathcal X}_{\alpha,\beta}$, if $x^g =y$,
then define $(\alpha x)^g =\alpha y$, $(\beta x)^g =\beta y$ and
$(\alpha\beta x)^g =\alpha\beta y$.}

Let $M=({\mathcal X}_{\alpha ,\beta },{\mathcal P})$ be a map.
According to the Theorem $1.2.5$, for an automorphism $g\in {\rm
Aut}M$ and $g|_{V(M)}:u\rightarrow v, \ u,v\in V(M)$, if $u^g= v$,
then call $g$ an {\it orientation-preserving automorphism}. if
$u^g= v^{-1}$, then call $g$ an {\it orientation-reversing
automorphism}. For any $g\in{\rm Aut}M$, it is obvious that {\it
$g|_{\Gamma}$ is orientation-preserving or orientation-reversing
and the product of two orientation-preserving automorphisms or
orientation-reversing automorphisms is orientation-preserving, the
product of an orientation-preserving automorphism with an
orientation-reversing automorphism is orientation-reversing}.

For a subgroup $G\preceq {\rm Aut}M$, define $G^+ \ \preceq \ G$
being the orientation-preserving subgroup of $G$. Then the index
of $G^+$ in $G$ is $2$. Assume the vertex $v$ to be
$v=(x_1,x_2,\cdots,x_{\rho (v)})(\alpha x_{\rho (v)},\cdots,\alpha
x_2,\alpha x_1)$. Denote by $<v>$ the cyclic group generated by
$v$. Then we get the following property for the automorphisms of a
map.

\vskip 4mm

\no{\bf Lemma $3.1.1$} {\it Let $G\preceq {\rm Aut}M$ be an
automorphism group of a map $M$. Then $\forall v\in V(M)$,

($i$) if $\forall g\in G$,$g$ is orientation-preserving, then
$G_v\preceq \ <v>$, is a cyclic group;

($ii$) $G_v\preceq \ <v>\times <\alpha >$.}

\vskip 3mm

{\it Proof}  ($i$)Let $M=({\mathcal X}_{\alpha ,\beta },{\mathcal
P})$. since for any $\forall g\in G$, $g$ is
orientation-preserving, we know for $\forall v\in V(M),h\in G_v$,
$v^h= v$. Assume the vertex

$$v=(x_1,x_2,\cdots,x_{\rho (v)})(\alpha x_{\rho (v)},\alpha
x_{\rho (v)-1},\cdots,\alpha x_1).$$

\no{Then}

$$
[(x_1,x_2,\cdots,x_{\rho (v)})(\alpha x_{\rho (v)},\cdots,\alpha
x_2,\alpha x_1)]^h=(x_1,x_2,\cdots,x_{\rho (v)})(\alpha x_{\rho
(v)},\cdots,\alpha x_2,\alpha x_1).
$$

Therefore, if $h(x_1)=x_{k+1},1\leq k\leq\rho (v)$, then

$$
h=[(x_1,x_2,\cdots,x_{\rho (v)})(\alpha x_{\rho (v)},\alpha
x_{\rho (v)-1},\cdots,\alpha x_1)]^k = v^k.
$$

\no{If $h(x_1)=\alpha x_{\rho (v)-k+1},1\leq k\leq\rho (v)$, then}

$$
h=[(x_1,x_2,\cdots,x_{\rho (v)})(\alpha x_{\rho (v)},\alpha
x_{\rho (v)-1},\cdots,\alpha x_1)]^k \alpha = v^k\alpha .
$$

But if $h=v^k\alpha$, then we know that $v^h=v^{\alpha}=v^{-1}$,
i.e., $h$ is not orientation-preserving. Whence, $h=v^k, 1\leq
k\leq\rho (v)$, i.e., each element in $G_v$ is the power of $v$.
Assume $\xi$ is the least power of elements in $G_v$. Then
$G_v=<v^{\xi}>\preceq \ <v>$ is a cyclic group generated by
$v^{\xi}$.

($ii$)For $\forall g\in G_v$, $v^g=v$, i.e.,

$$
[(x_1,x_2,\cdots,x_{\rho})(\alpha x_{\rho},\alpha
x_{\rho-1},\cdots,\alpha x_1)]^g=(x_1,x_2,\cdots,x_{\rho})(\alpha
x_{\rho},\alpha x_{\rho-1},\cdots,\alpha x_1).
$$

Similar to the proof of $(i)$, we know there exists an integer
$s,1\leq s\leq\rho$, such that $g=v^s$ or $g=v^s\alpha$. Whence,
$g\in <v>$ or $g\in <v>\alpha$, i.e.,

$$G_v\preceq \ <v>\times <\alpha>.\quad\quad \natural$$

\vskip 4mm

\no{\bf Lemma $3.1.2$} {\it Let $\Gamma$ be a connected graph. If
$G\preceq {\rm Aut}\Gamma$, and $\forall v\in
V(\Gamma)$,$G_v\preceq \ <v>\times <\alpha>$, then the action of
$G$ on ${\mathcal X}_{\alpha,\beta}$ is fixed-free.}

\vskip 3mm

{\it Proof} Choose a quadricell $x\in {\mathcal X}_{\alpha
,\beta}$. We prove that $G_x=\{{\bf 1}_G\}$. In fact, if $g\in
G_x$, then $x^g=x$. Particularly, the incident vertex $u$ is
stable under the action of $g$, i.e., $u^g=u$. assume

$$u=(x,y_1,\cdots,y_{\rho (u)-1})(\alpha x,\alpha
y_{\rho (u)-1},\cdots,\alpha y_1),$$

\no{then since $G_u\preceq \ <u>\times <\alpha>$, we get that}

$$x^g=x,y_1^g=y_1,\cdots, y_{\rho (u)-1}^g=y_{\rho (u)-1}$$

\no{and}

$$(\alpha x)^g=\alpha x,(\alpha y_1)^g=\alpha y_1,\cdots, (\alpha
y_{\rho (u)-1})^g=\alpha y_{\rho (u)-1},$$

\no{that is, for any quadricell $e_u$ incident with the vertex
 $u$, $e_u^g=e_u$. According to the definition of the induced action ${\rm
Aut}\Gamma$ on ${\mathcal X}_{\alpha,\beta}$, we know that}

$$(\beta x)^g=\beta x,(\beta y_1)^g=\beta y_1,\cdots,
(\beta y_{\rho (u)-1})^g=\beta y_{\rho (u)-1}$$

\no{and}

$$(\alpha\beta x)^g=\alpha\beta x,(\alpha\beta y_1)^g=\alpha\beta y_1,\cdots,
(\alpha\beta y_{\rho (u)-1})^g=\alpha\beta y_{\rho (u)-1}.$$

Whence, for any quadricell $y\in{\mathcal X}_{\alpha,\beta}$,
assume the incident vertex of $y$ is $w$, then by the connectivity
of the graph $\Gamma$, we know that there is a path
$P(u,w)=uv_1v_2\cdots v_sw$ in $\Gamma$ connecting the vertex $u$
and$w$. Not loss of generality, we assume that $\beta y_k$ is
incident with the vertex $v_1$.  Since $(\beta y_k)^g=\beta y_k$
and $G_{v_1}\preceq \ <v_1>\times <\alpha>$, we know that for any
quadricell $e_{v_1}$ incident with the vertex $v_1$, $
e_{v_1}^g=e_{v_1}$.

Similarly, if a quadricell $e_{v_i}$ incident with the vertex
$v_i$ is stable under the action of $g$, i.e.,
$(e_{v_i})^g=e_{v_i}$, then we can prove that any quadricell
$e_{v_{i+1}}$ incident with the vertex $v_{i+1}$ is stable under
the action of $g$. This process can be well done until we arrive
the vertex $w$. Therefore, we can get that any quadricell $e_{w}$
incident with the vertex $w$ is stable under the action of $g$.
Particularly, we have that $y^g=y$.

Therefore, we get that $g={\bf 1}_G$. Whence,$G_x=\{{\bf
1}_G\}$.\quad\quad $\natural$

Now we prove a necessary and sufficient condition for a subgroup
of a graph being an automorphism group of a map underlying this
graph.

\vskip 4mm

\no{\bf Theorem $3.1.1$}  {\it Let $\Gamma$ be a connected graph.
If $G \ \preceq \ {\rm Aut}\Gamma$, then $G$ is an automorphism
group of a map underlying the graph $\Gamma$ iff for $\forall v\in
V(\Gamma)$, the stabler $G_v\preceq \ <v>\times <\alpha>$.}

\vskip 3mm

{\it Proof} According to the Lemma $3.1.1(ii)$, the condition of
the Theorem $3.1.1$ is necessary. Now we prove its sufficiency.

By the Lemma $3.1.2$, we know that the action of $G$ on ${\mathcal
X}_{\alpha ,\beta}$ is fixed-free, i.e., for$\forall x\in
{\mathcal X}_{\alpha ,\beta}$, $|G_x|=1$. Whence, the length of
orbit of $x$ under the action $G$ is $|x^G|=|G_x||x^G|=|G|$, i.e.,
for $\forall x\in {\mathcal X}_{\alpha ,\beta}$, the length of $x$
under the action of $G$ is $|G|$.

Assume that there are $s$ orbits $O_1,O_2,\cdots,O_s$ of $G$
action on ÔÚ$V(\Gamma )$, where, $O_1=\{u_1,u_2,\cdots,u_k\}$,
$O_2=\{v_1,v_2,\cdots,v_l\}$,$\cdots$,$O_s=\{w_1,w_2,\cdots,w_t\}$.
We construct the conjugatcy permutation pair for every vertex in
the graph $\Gamma $ such that they product ${\mathcal P}$ is
stable under the action of $G$.

Notice that for $\forall u\in V(\Gamma )$, since $|G|=|G_u||u^G|$,
we know that $[k,l,\cdots,t]| \ |G|$.

In the first, we determine the conjugate permutation pairs for
each vertex in the orbit $O_1$. Choose any vertex $u_1\in O_1$,
assume that the stabler $G_{u_1}$ is $\{1_G,g_1,g_2g_1,\cdots,
\prod\limits_{i=1}^{m-1}g_{m-i}\}$, where, $m=|G_{u_1}|$ and the
quadricells incident with vertex $u_1$ is $\widetilde{N(u_1)}$ in
the graph $\Gamma$ . We arrange the elements in
$\widetilde{N(u_1)}$ as follows.

Choose a quadricell $u_1^a\in \widetilde{N(u_1)}$. We use
$G_{u_1}$ action on $u_1^a$ and $\alpha u_1^a$, respectively. Then
we get the quadricell set $A_1=\{u_1^a,g_1(u_1^a),\cdots
,\prod\limits_{i=1}^{m-1}g_{m-i}(u_1^a)\}$ and $\alpha
A_1=\{\alpha u_1^a,\alpha g_1(u_1^a),\cdots
,\alpha\prod\limits_{i=1}^{m-1}g_{m-i}(u_1^a)\}$. By the
definition of the action of an automorphism of a graph on its
quadricells we know that $A_1\bigcap\alpha A_1=\emptyset$. Arrange
the elements in $A_1$ as $\overrightarrow{A_1} =
u_1^a,g_1(u_1^a),\cdots ,\prod\limits_{i=1}^{m-1}g_{m-i}(u_1^a)$.

If $\widetilde{N(u_1)}\setminus A_1\bigcup\alpha A_1=\emptyset$,
then the arrangement of elements in $\widetilde{N(u_1)}$ is
$\overrightarrow{A_1}$. If $\widetilde{N(u_1)}\setminus
A_1\bigcup\alpha A_1\not=\emptyset$, choose a quadricell $u_1^b\in
\widetilde{N(u_1)}\setminus A_1\bigcup \alpha A_1$. Similarly,
using the group $G_{u_1}$ acts on $u_1^b$, we get that
$A_2=\{u_1^b,g_1(u_1^b),\cdots
,\prod\limits_{i=1}^{m-1}g_{m-i}(u_1^b)\}$ and $\alpha
A_2=\{\alpha u_1^b,\alpha g_1(u_1^b),\cdots
,\alpha\prod\limits_{i=1}^{m-1}g_{m-i}(u_1^b)\}$. Arrange the
elements in $A_1\bigcup A_2$ as

$$\overrightarrow{A_1\bigcup A_2} = u_1^a,g_1(u_1^a),\cdots
,\prod\limits_{i=1}^{m-1}g_{m-i}(u_1^a); u_1^b,g_1(u_1^b),\cdots
,\prod\limits_{i=1}^{m-1}g_{m-i}(u_1^b).$$

If $\widetilde{N(u_1)}\setminus (A_1\bigcup A_2\bigcup\alpha
A_1\bigcup \alpha A_2) =\emptyset$, then the arrangement of
elements in $A_1\bigcup A_2$ is $\overrightarrow{A_1\bigcup A_2}$.
Otherwise, $\widetilde{N(u_1)}\setminus (A_1\bigcup
A_2\bigcup\alpha A_1\bigcup \alpha A_2) \not=\emptyset$. We can
choose another quadricell $u_1^c\in\widetilde{N(u_1)}\setminus
(A_1\bigcup A_2\bigcup$ $\alpha A_1\bigcup \alpha A_2)$.
Generally, If we have gotten the quadricell sets $A_1,A_2,\cdots
,A_r, 1\leq r\leq 2k$, and the arrangement of element in them is
$\overrightarrow{A_1\bigcup A_2\bigcup\cdots\bigcup A_r}$, if
$\widetilde{N(u_1)}\setminus £¨A_1\bigcup A_2\bigcup\cdots\bigcup
A_r\bigcup\alpha A_1\bigcup \alpha A_2\bigcup$ $\cdots\bigcup
\alpha A_r) \not=\emptyset$, then we can choose an element
$u_1^d\in\widetilde{N(u_1)}\setminus (A_1\bigcup
A_2\bigcup\cdots\bigcup A_r\bigcup\alpha A_1$ $\bigcup \alpha
A_2\bigcup\cdots \bigcup \alpha A_r)$ and define the quadricell
set

$$A_{r+1}=\{u_1^d,g_1(u_1^d),\cdots
,\prod\limits_{i=1}^{m-1}g_{m-i}(u_1^d)\}$$

$$\alpha A_{r+1}=\{\alpha u_1^d,\alpha g_1(u_1^d),\cdots
,\alpha\prod\limits_{i=1}^{m-1}g_{m-i}(u_1^d)\}$$

\no{and the arrangement of elements in $A_{r+1}$ is}

$$\overrightarrow{A_{r+1}}=u_1^d,g_1(u_1^d),\cdots
,\prod\limits_{i=1}^{m-1}g_{m-i}(u_1^d).$$

\no{Define the arrangement of elements in
$\bigcup\limits_{j=1}^{r+1}A_j$ to be}

$$\overrightarrow{\bigcup\limits_{j=1}^{r+1}A_j}=
\overrightarrow{\bigcup\limits_{i=1}^{r}A_i};\overrightarrow{A_{r+1}}
.$$

\no{Whence,}

$$
\widetilde{N(u_1)}=(\bigcup\limits_{j=1}^kA_j)\bigcup
(\alpha\bigcup\limits_{j=1}^kA_j)
$$

\no{and $A_k$ is obtained by the action of the stabler $G_{u_1}$
on $u_1^e$. At the same time, the arrangement of elements in the
subset $\bigcup\limits_{j=1}^kA_j$ of $\widetilde{N(u_1)}$ to be
$\overrightarrow{\bigcup\limits_{j=1}^kA_j}$.}

Define the conjugatcy permutation pair of the vertex $u_1$ to be

$$\varrho_{u_1}=(C)(\alpha C^{-1}\alpha),$$

\no{where£¬}

$$
C=(u_1^a,u_1^b,\cdots,u_1^e;g_1(u_1^a),g_1(u_1^b),\cdots,g_1(u_1^e),\cdots,
\prod\limits_{i=1}^{m-1}(u_1^a),\prod\limits_{i=1}^{m-1}(u_1^b),\cdots
,\prod\limits_{i=1}^{m-1}(u_1^e)).
$$

For any vertex $u_i\in O_1,1\leq i\leq k$, assume that
$h(u_1)=u_i$, where $h\in G$, then we define the conjugatcy
permutation pair $\varrho_{u_i}$ of the vertex $u_i$ to be

$$\varrho_{u_i}= \varrho_{u_1}^h = (C^h)(\alpha C^{-1}\alpha^{-1}).$$

Since $O_1$ is an orbit of the action $G$ on $V(\Gamma )$, then we
have that

$$
(\prod\limits_{i=1}^k \varrho_{u_i})^G \ = \ \prod\limits_{i=1}^k
\varrho_{u_i}.
$$

Similarly, we can define the conjugatcy permutation pairs
$\varrho_{v_1},\varrho_{v_2},\cdots,\varrho_{v_l},\cdots,
\varrho_{w_1},$ $\varrho_{w_2},\cdots,\varrho_{w_t}$ of vertices
in the orbits $O_2,\cdots,O_s$. We also have that

$$
(\prod\limits_{i=1}^l \varrho_{v_i})^G \ = \ \prod\limits_{i=1}^l
\varrho_{v_i}.
$$

$$
\cdots\cdots\cdots\cdots\cdots\cdots\cdots
$$

$$
(\prod\limits_{i=1}^t \varrho_{w_i})^G \ = \ \prod\limits_{i=1}^t
\varrho_{w_i}.
$$

Now define the permutation

$$
{\mathcal P}=(\prod\limits_{i=1}^k \varrho_{u_i})\times
(\prod\limits_{i=1}^l \varrho_{v_i})\times\cdots\times
(\prod\limits_{i=1}^t \varrho_{w_i}).
$$

Then since $O_1,O_2,\cdots,O_s$ are the orbits of $V(\Gamma )$
under the action of $G$, we get that

\begin{eqnarray*}
{\mathcal P}^G &=& (\prod\limits_{i=1}^k \varrho_{u_i})^G\times
(\prod\limits_{i=1}^l \varrho_{v_i})^G\times\cdots\times
(\prod\limits_{i=1}^t \varrho_{w_i})^G\\
&=& (\prod\limits_{i=1}^k \varrho_{u_i})\times
(\prod\limits_{i=1}^l \varrho_{v_i})\times\cdots\times
(\prod\limits_{i=1}^t \varrho_{w_i})= {\mathcal P}.
\end{eqnarray*}

Whence, if we define the map $M=({\mathcal X}_{\alpha
,\beta},{\mathcal P})$£¬then $G$ is an automorphism of the map
$M$.\quad\quad $\natural$

For the orientation-preserving automorphism, we have the following
result.

\vskip 4mm

\no{\bf Theorem $3.1.2$}  {\it Let $\Gamma$ be a connected graph.
If \ $G \ \preceq \ {\rm Aut}\Gamma$, then $G$ is an
orientation-preserving automorphism group of a map underlying the
graph $\Gamma$ iff for $\forall v\in V(\Gamma)$, the stabler $G_v
\ \preceq \ <v>$ is a cyclic group.}

\vskip 3mm

{\it Proof}   According to the Lemma $3.1.1(i)$£¬we know the
necessary. Notice the approach of construction the conjugatcy
permutation pair in the proof of the Theorem $3.1.1$. We know that
$G$ is also an orientation-preserving automorphism group of the
map $M$.\quad\quad $\natural$

According to the Theorem $3.1.2$, we can prove the Theorem $2.4.3$
now.

\vskip 4mm

\no{\bf\it The Proof of the Theorem  $2.4.3$}

\vskip 3mm

Since every subgroup of a cyclic group is also a cyclic group, we
know that any cyclic orientation-preserving automorphism group of
the graph $\Gamma$ is an orientation-preserving automorphism group
of a map underlying the graph $\Gamma$ by the Theorem $3.1.2$.
Whence, we get that

$$o_{max}(M,g) \ \leq \ o_{max}(G,g).\quad\quad \natural$$

\vskip 4mm

\no{\bf Corollary $3.1.1$}   {\it For any positive integer $n$,
there exists a vertex transitive map $M$ underlying a circultant
such that $Z_n$ is an orientation-preserving automorphism group of
the map $M$.}

\vskip 3mm

\no{\bf Remark $3.1.1$} Gardiner et al proved in $[20]$ that if
add an additional condition , i.e, $G$ is transitive on the
vertices in $\Gamma$, then there is a regular map underlying the
graph $\Gamma$.

\vskip 10mm

{\bf\S 2. The automorphisms of a complete graph on surfaces}

\vskip 8mm

A map is called a {\it complete map} if its underlying graph is a
complete graph. For a connected graph $\Gamma$, the notations
${\cal E }^{O}(\Gamma ), {\cal E}^{N}(\Gamma )$ and ${\cal
E}^{L}(\Gamma )$ denote the embeddings of $\Gamma$ on the
orientable surfaces, non-orientable surfaces and locally surfaces,
respectively. For $\forall e=(u,v)\in E(\Gamma)$, its quadricell
$Ke=\{e,\alpha e,\beta e, \alpha\beta e\}$ can be represented by
$Ke=\{u^{v+},u^{v-},v^{u+},v^{u-}\}$.

Let $K_n$ be a complete graph of order $n$. Label its vertices by
integers $1,2,...,n$. Then its edge set is $\{ij| 1\leq i,j\leq n
,i\not=j\, ij=ji\}$, and

$${\cal X}_{\alpha,\beta} (K_n)=\{ i^{j+} : 1\leq i,j\leq n ,i\not=j \}
\bigcup\{ i^{j-} : 1\leq i,j\leq n ,i\not=j \},$$

$$\alpha = \prod\limits_{1\leq i,j\leq n ,i\not=j}  ( i^{j+} , i^{j-}),$$

$$\beta = \prod\limits_{1\leq i,j\leq n ,i\not=j}  ( i^{j+} , i^{j+})( i^{j-} , i^{j-}).$$

We determine all the automorphisms of complete maps of order $n$
and give they concrete representation in this section.

First, we need some useful lemmas for an automorphism of a map
induced by an automorphism of its underlying graph.

\vskip 4mm

\no{{\bf Lemma $3.2.1$} {\it Let $\Gamma$ be a connected graph and
$g\in{\rm Aut \Gamma}$. If there is a map $M\in {\cal
E}^{L}(\Gamma)$ such that the induced action $g^* \in{\rm AutM}$,
then for $\forall (u,v),(x,y)\in{E(\Gamma)}$,}}

$$[l^g(u),l^g(v)]=[l^g(x),l^g(y)]= constant,$$

\no{\it where, $l^g(w)$ denotes the length of the cycle containing
the vertex $w$ in the cycle decomposition of $g$.}

\vskip 3mm

 {\it Proof}  According to the Lemma $2.2.1$, we know that the length of
a  quadricell $u^{v+}$ or $u^{v-}$ under the action $g^*$ is
$[l^g(u),l^g(v)]$. Since $g^*$ is an automorphism of map,
therefore, $g^*$ is semi-regular. Whence, we get that
 $$ [l^g(u),l^g(v)]=[l^g(x),l^g(y)]= constant.\quad\quad\natural$$

Now we consider conditions for an induced automorphism of a map by
an automorphism of graph to be an  orientation-reversing
automorphism of a map.

\vskip 4mm

\no{\bf Lemma $3.2.2$} {\it If $\xi\alpha$ is an automorphism of a
map, then $\xi\alpha =\alpha\xi .$}

\vskip 3mm

{\it Proof} Since $\xi\alpha$ is an automorphism of a map, we know
that

$$(\xi\alpha)\alpha = \alpha (\xi\alpha).$$

\no{That is, $\xi\alpha =\alpha\xi .\quad\quad\natural$}

\vskip 4mm

\no{\bf Lemma $3.2.3$} {\it If $\xi$ is an automorphism of map
$M=({\cal X}_{\alpha,\beta},{\cal P})$, then $\xi\alpha $ is
semi-regular on ${\cal X}_{\alpha,\beta}$ with order $o(\xi)$ if
$o(\xi)\equiv 0(mod 2)$ and $2o(\xi)$ if $o(\xi)\equiv 1(mod 2)$.}

\vskip 3mm

{\it Proof} Since $\xi$ is an automorphism of map by the Lemma
$3.2.2$, we know that the cycle decomposition of $\xi$ can be
represented by

$$\xi=\prod_k (x_1,x_2,\cdots,x_k)(\alpha x_1,\alpha x_2,\cdots, \alpha x_k),$$

\no{where, $\prod_k$ denotes the product of disjoint cycles with
length $k=o(\xi)$.}

Therefore,  if $k\equiv 0 (mod 2)$, we get that

$$\xi\alpha = \prod_k (x_1,\alpha x_2,x_3,\cdots,\alpha x_k)$$

\no{and if $k\equiv 1 (mod 2)$, we get that}

$$\xi\alpha = \prod_{2k} (x_1,\alpha x_2,x_3,\cdots, x_k,\alpha x_1,x_2,\alpha x_3,\cdots,\alpha x_k).$$

\no{Whence, $\xi$ is semi-regular acting on ${\cal
X}_{\alpha,\beta}$.\quad\quad $\natural$}

Now we can prove the following result for orientation-reversing
automorphisms of a map.

\vskip 4mm

\no{\bf Lemma $3.2.4$}  {\it For a connected graph $\Gamma$, let
${\cal K}$  be all automorphisms in $Aut\Gamma$ whose extending
action on ${\cal X}_{\alpha,\beta}$, $ X=E(\Gamma)$, are
automorphisms of maps underlying the graph $\Gamma$. Then for
$\forall\xi\in {\cal K}$, $o(\xi^*)\geq 2$, $\xi^*\alpha \in {\cal
K}$ if and only if $o(\xi^*)\equiv 0(mod 2).$}

 \vskip 3mm

{\it Proof} Notice that by the Lemma $3.2.3$, if $\xi^*$ is an
automorphism of map with underlying graph $\Gamma$, then
$\xi^*\alpha$ is semi-regular acting on ${\cal X}_{\alpha,\beta}$.

Assume $\xi^*$ is an automorphism of map $M=({\cal
X}_{\alpha,\beta},{\cal P})$. Without loss of generality, we
assume that

$$
{\cal P}=C_1C_2\cdots C_k,
$$
\no{where,$C_i =(x_{i1},x_{i2},\cdots,x_{ij_i})$ is a cycle in the
decomposition of $\xi |_{V(\Gamma)}$ and
$x_{it}=\{(e^{i1},e^{i2},\cdots,e^{it_i})(\alpha e^{i1},\alpha
e^{it_i},\cdots,\alpha e^{i2})\}$ and.}

$$\xi |_{E(\Gamma)}=(e_{11},e_{12},\cdots,e_{s_1})(e_{21},e_{22},\cdots,e_{2s_2})\cdots
(e_{l1},e_{l2},\cdots,e_{ls_l}).$$

\no{and}
$$
\xi^*=C (\alpha C^{-1}\alpha ),
$$

\no{where,
$C=(e_{11},e_{12},\cdots,e_{s_1})(e_{21},e_{22},\cdots,e_{2s_2})\cdots
(e_{l1},e_{l2},\cdots,e_{ls_l})$. Now since $\xi^*$ is an
automorphism of map, we get that $s_1=s_2=\cdots =s_l
=o(\xi^*)=s.$}

If $o(\xi^*)\equiv 0(mod 2)$, define a map $M^*=({\cal
X}_{\alpha,\beta},{\cal P}^*)$, where,

$$
{\cal P}^*=C_1^*C_2^*\cdots C_k^*,
$$
\no{where, $C_i^* =(x_{i1}^*,x_{i2}^*,\cdots,x_{ij_i}^*)$,
$x_{it}^* =\{(e_{i1}^*,e_{i2}^*,\cdots,e_{it_i}^*)(\alpha
e_{i1}^*,\alpha e_{it_i}^*,\cdots, e_{i2}^*)\}$ and
$e_{ij}^*=e_{pq}$. Take $e_{ij}^*=e_{pq}$ if  $ q\equiv 1(mod 2)$
and $e_{ij}^*=\alpha e_{pq}$ if $q\equiv 0(mod 2)$.  Then we get
that $M^{\xi\alpha}= M$.}

Now if $o(\xi^*)\equiv 1(mod 2)$, by the Lemma $3.2.3$,
$o(\xi^*\alpha)=2o(\xi^*).$ Therefore, any chosen quadricells
$(e^{i1},e^{i2},\cdots,e^{it_i})$ adjacent to the vertex $x_{i1}$
for $i=1,2,\cdots,n$, where, $n=|\Gamma|$, the resultant map $M$
is unstable under the action of $\xi\alpha$. Whence, $\xi\alpha$
is not an automorphism of one map with underlying graph $\Gamma$.
\quad\quad $\natural$

Now we determine all automorphisms of a complete map underlying a
graph $K_n$ by applying the previous results. Recall that the
automorphism group of the graph $K_n$ is just the symmetry group
of degree $n$, that is, ${\rm AutK}_n=S_n$.

\vskip 4mm

\no{\bf Theorem $3.2.1$} {\it All orientation-preserving
automorphisms of non-orientable complete maps of order$\geq 4$ are
extended actions of elements in}

$${\cal E}_{[s^{\frac{n}{s}}]}, \quad {\cal E}_{[1,s^{\frac{n-1}{s}}]},$$

\no{\it and all orientation-reversing automorphisms of
non-orientable complete maps of order$\geq 4$ are extended actions
of elements in}

$$
\alpha {\cal E}_{[(2s)^{\frac{n}{2s}}]},\quad \alpha {\cal
E}_{[(2s)^{\frac{4}{2s}}]}, \quad\alpha {\cal E}_{[1,1,2]},
$$
\no{\it where, ${\cal E}_{\theta}$ denotes the conjugatcy class
containing element $\theta$ in the symmetry group $S_n$}

\vskip 3mm

{\it Proof} First, we prove that the induced permutation $\xi^*$
on a complete map of order $n$ by an element $\xi\in S_n$ is an
cyclic order-preserving automorphism of non-orientable map, if,
and only if,

$$\xi\in {\cal E}_{s^{\frac{n}{s}}}\bigcup {\cal E}_{[1,s^{\frac{n-1}{s}}]}$$

Assume the cycle index of $\xi$ is
$[1^{k_1},2^{k_2},...,n^{k_n}]$. If there exist two integers
$k_i,k_j \not=0$, and $i,j\geq 2, i\not=j$, then in the cycle
decomposition of $\xi$ , there are two cycles

$$(u_1,u_2,...,u_i)\quad  {\rm and}\quad  (v_1,v_2,...,v_j).$$

Since
$$[l^{\xi}(u_1), l^{\xi}(u_2)]=i\quad  {\rm and}\quad  [l^{\xi}(v_1), l^{\xi}(v_2)]=j$$

\no{and $i\not= j$, we know that $\xi^*$ is not an automorphism of
embedding by the Lemma $2.5$. Whence, the cycle index of $\xi$
must be the form of $[1^k,s^l]$.}

Now if $k\geq 2$, let $(u),(v)$ be two cycles of length $1$ in the
cycle decomposition of $\xi$. By the Lemma $2.5$, we know that

$$[l^{\xi}(u), l^{\xi}(v)] = 1.$$

If there is a cycle $(w,...)$ in the cycle decomposition of $\xi$
whose length greater or equal to two, we get that

$$[l^{\xi}(u), l^{\xi}(w)]=[1, l^{\xi}(w)]= l^{\xi}(w).$$

According to the Lemma $3.2.1$, we get that $l^{\xi}(w)=1$, a
contradiction. Therefore, the cycle index of $\xi$ must be the
forms of $[s^l]$ or $[1, s^l]$. Whence,$ sl=n$ or $sl+1=n$.
Calculation shows that $l=\frac{n}{s}$ or $l=\frac{n-1}{s}$. That
is, the cycle index of $\xi$ is one of the following three types
$[1^n]$,  $[1,s^{\frac{n-1}{s}}]$  and   $[s^{\frac{n}{s}}]$  for
some integer $s$ .

Now we only need to prove that for each element $\xi$ in ${\cal
E}_{[1,s^{\frac{n-1}{s}}]}$  and ${\cal E}_{[s^{\frac{n}{s}}]}$ ,
there exists an non-orientable complete map $M$ of order $n$ with
the induced permutation $\xi^*$ being its cyclic order-preserving
automorphism of surface. The discussion are divided into two
cases.

\vskip 2mm

{\bf Case 1}   \hspace{40mm}    $ \xi\in{ {\cal
E}_{[s^{\frac{n}{s}}]}}$

\vskip 2mm

Assume the cycle decomposition of $\xi$ being $\xi =
(a,b,\cdots,c)\cdots(x,y,\cdots,z)\cdots(u,v,$ $\cdots,w)$, where,
the length of each cycle is $k$, and $1\leq
a,b,\cdots,c,x,y,\cdots,z,u,v,\cdots,w\leq n$ . In this case, we
can construct a non-orientable complete map $M_{1} = ({\cal
X}_{\alpha,\beta}^{1},{\cal{P}}_{1})$ as follows.

$${\cal X}_{\alpha,\beta}^{1}=\{ i^{j+} : 1\leq i,j\leq n ,i(j \}\bigcup
\{ i^{j-} : 1\leq i,j\leq n ,i\not=j \},$$

$${\cal{P}}_{1}= \prod\limits_{x\in\{a,b,\cdots,c,\cdots,x,y,\cdots,z,u,v,\cdots,w\}}
(C(x))(\alpha C(x)^{-1}\alpha),$$

\no{where,}

$$C(x) = (x^{a+},\cdots ,x^{x*},\cdots,x^{u+},x^{b+},x^{y+},\cdots,
\cdots,x^{v+},x^{c+},\cdots,x^{z+},\cdots ,x^{w+}),$$

\no{where $x^{x*}$ denotes an empty position and }

$$\alpha C(x)^{-1}\alpha =(x^{a-},x^{w-},\cdots,x^{z-},\cdots ,x^{c-},x^{v-},\cdots,
x^{b-},x^{u-},\cdots,x^{y-},\cdots ).$$

It is clear that $M_1^{\xi^*} = M_1$. Therefore, $\xi^*$ is an
cyclic order-preserving automorphism of the map $M_1$.

\vskip 3mm {\bf Case 2}  \hspace{40mm}  $\xi \in{ {\cal
E}_{[1,s^{\frac{n-1}{s}}]}}$ \vskip 3mm We assume the cycle
decomposition of $\xi$ being

$$\xi = (a,b,...,c)...(x,y,...,z)...(u,v,...,w)(t),$$

\no{where, the length of each cycle is $k$ beside the final cycle,
and $1\leq a,b...c, x,y..., z,$ $u,v,...,w,t\leq n$ . In this
case, we construct a non-orientable complete map $M_{2} = ({\cal
X}_{\alpha,\beta}^{2},{\cal{P}}_{2})$
 as follows.}

$${\cal X}_{\alpha,\beta}^{2}=\{ i^{j+} : 1\leq i,j\leq n ,i\not=j \}\bigcup
\{ i^{j-} : 1\leq i,j\leq n ,i\not=j \},$$

$${\cal{P}}_{2} = (A)( \alpha A^{-1})\prod\limits_{x\in\{a,b,...,c,...,
x,y,...z,u,v,...,w\}} (C(x))(\alpha C(x)^{-1}\alpha),$$

\no{where,}

$$A= (t^{a+},t^{x+},...t^{u+},t^{b+},t^{y+},...,t^{v+},...,
t^{c+},t^{z+},...,t^{w+})$$

\no{and}

$$\alpha A^{-1}\alpha=(t^{a-},t^{w-},...t^{z-},t^{c-},t^{v-},...,t^{y-},...,
t^{b-},t^{u-},...,t^{x-})$$

\no{and}
$$C(x) = (x^{a+},...,x^{x*},...,x^{u+},x^{b+},...,x^{y+},...,x^{v+},...,x^{c+},...,
x^{z+},...,x^{w+})$$

\no{and}

$$\alpha C(x)^{-1}\alpha =(x^{a-},x^{w-},..,x^{z-},...,x^{c-},...,x^{v-},...,
x^{y-},...,x^{b-},x^{u-},...).$$

It is also clear that $M_2^{\xi^*} = M_2$. Therefore, $\xi^*$ is
an automorphism of a map $M_2$ .

Now we consider the case of orientation-reversing automorphism of
a complete map. According to the Lemma $3.2.4$, we know that an
element $\xi\alpha$, where, $\xi\in S_n$, is an
orientation-reversing automorphism of a complete map, only if,

$$\xi\in {\cal E}_{[k^{\frac{n_1}{k}},(2k)^{\frac{n-n_1}{2k}}]}.$$

Our discussion is divided into two parts.

\vskip 3mm

\no{\bf Case $3$}\quad\quad\quad\quad   $n_1 = n$

\vskip 3mm

Without loss of generality, we can assume the cycle decomposition
of $\xi$ has the following form in this case.

$$
\xi
=(1,2,\cdots,k)(k+1,k+2,\cdots,2k)\cdots(n-k+1,n-k+2,\cdots,n).
$$

\vskip 3mm

\no{\bf Subcase $3.1$}\quad\quad\quad $k\equiv 1(mod 2)$ and $k >
1$

\vskip 3mm

According to the Lemma $3.2.4$, we know that $\xi^*\alpha$ is not
an automorphism of map since $o(\xi^*)=k\equiv 1(mod 2)$.

\vskip 3mm

\no{\bf Subcase $3.2$} \quad\quad\quad\quad  $k\equiv 0(mod 2)$

\vskip 3mm

Construct a non-orientable map $M_3=({\cal
X}_{\alpha,\beta}^3,{\cal P}_3)$, where $ X^3= E(K_n)$ and

$$
{\cal P}_3=\prod\limits_{i\in \{1,2,\cdots,n\}}(C(i))(\alpha
C(i)^{-1}\alpha),
$$

\no{where, if $i\equiv 1(mod 2)$, then}

$$C(i)=(i^{1+},i^{k+1+},\cdots,i^{n-k+1+},i^{2+},\cdots,i^{n-k+2+},\cdots,i^{i*},\cdots,i^{k+},i^{2k+},\cdots,i^{n+}),$$

$$\alpha C(i)^{-1}\alpha=(i^{1-},i^{n-},\cdots,i^{2k-},i^{k-},\cdots,i^{k+1-})$$

\no{and if $i\equiv 0(mod 2)$, then}

$$C(i)=(i^{1-},i^{k+1-},\cdots,i^{n-k+1-},i^{2-},\cdots,i^{n-k+2-},\cdots,i^{i*},\cdots,i^{k-},i^{2k-},\cdots,i^{n-}),$$

$$\alpha C(i)^{-1}\alpha=(i^{1+},i^{n+},\cdots,i^{2k+},i^{k+},\cdots,i^{k+1+}).$$

Where, $i^{i*}$ denotes the empty position, for example,
$(2^1,2^{2*},2^3,2^4,2^5)=(2^1,2^3,2^4,2^5)$. It is clear that
${\cal P}_3^{\xi\alpha}={\cal P}_3$, that is, $\xi\alpha$ is an
automorphism of map $M_3$.

\vskip 3mm

\no{\bf Case 4}\quad\quad\quad\quad $n_1\not= n$

\vskip 3mm

Without loss of generality, we can assume that

\begin{eqnarray*}
\xi &=& (1,2,\cdots,k)(k+1,k+2,\cdots,n_1)\cdots (n_1-k+1,n_1-k+2,\cdots,n_1)\\
&\times&
(n_1+1,n_1+2,\cdots,n_1+2k)(n_1+2k+1,\cdots,n_1+4k)\cdots(n-2k+1,\cdots,n)
\end{eqnarray*}

\vskip 3mm

\no{\bf Subcase $4.1$} \quad\quad\quad\quad $k\equiv 0(mod 2)$

\vskip 3mm

Consider the orbits of $1^{2+}$ and $n_1+2k+1^{1+}$ under the
action of $<\xi\alpha >$, we get that

$$|orb((1^{2+})^{<\xi\alpha >})| = k$$

\no{and}

$$|orb(((n_1+2k+1)^{1+})^{< \xi\alpha >})| = 2k.$$

\no{Contradicts to the Lemma $3.2.1$.}

\vskip 3mm

\no{\bf Subcase $4.2$} \quad\quad\quad\quad $k\equiv 1(mod 2)$

\vskip 3mm

In this case, if $k\not= 1$, then $k\geq 3$. Similar to the
discussion of the Subcase $3.1$, we know that $\xi\alpha$ is not
an automorphism of complete map. Whence, $k=1$ and

$$\xi\in {\cal E}_{[1^{n_1},2^{n_2}]}.$$

\no{Without loss of generality, assume that}

$$\xi =(1)(2)\cdots(n_1)(n_1+1,n_1+2)(n_1+3,n_1+4)\cdots(n_1+n_2-1,n_1+n_2).$$

If $n_2\geq 2$, and there exists a map $M=({\cal
X}_{\alpha,\beta},{\cal P})$, assume the vertex $v_1$ in $M$ being

$$
v_1=(1^{l_{12}+},1^{l_{13}+},\cdots,1^{l_{1n}+})(1^{l_{12}-},1^{l_{1n}-},\cdots,1^{l_{13}-})
$$
\no{where, $l_{1i}\in\{+2,-2,+3,-3,\cdots,+n,-n\}$ and
$l_{1i}\not= l_{1j}$ if $i\not= j$.}

Then we get that

$$
(v_1)^{\xi\alpha}=(1^{l_{12}-},1^{l_{13}-},\cdots,1^{l_{1n}-})(1^{l_{12}+},1^{l_{1n}+},\cdots,1^{l_{13}+})\not=v_1.
$$

\no{Whence, $\xi\alpha$ is not an automorphism of map $M$. A
contradiction.

Therefore, $n_2 =1$. Similarly, we can also get that $n_1 =2$.
Whence,  $\xi =(1)(2)(34)$ and $n=4$. We construct a stable
non-orientable map $M_4$ under the action of $\xi\alpha$ as
follows.

$$M_4=({\cal X}_{\alpha,\beta}^4,{\cal P}_4),$$

\no{where,}

\begin{eqnarray*}
{\cal P}_4  &=& (1^{2+},1^{3+},1^{4+})(2^{1+},2^{3+},2^{4+})(3^{1+},3^{2+},3^{4+})(4^{1+},4^{2+},4^{3+})\\
& \times &
(1^{2-},1^{4-},1^{3-})(2^{1-},2^{4-},2^{3-})(3^{1-},3^{4-},3^{2-})(4^{1-},4^{3-},4^{2-}).
\end{eqnarray*}

Therefore, all orientation-preserving automorphisms of
non-orientable complete maps are extended actions of elements in

$${\cal E}_{[s^{\frac{n}{s}}]}, \quad {\cal E}_{[1,s^{\frac{n-1}{s}}]}$$

\no{and all orientation-reversing automorphisms of non-orientable
complete maps are extended actions of elements in}

$$
 \alpha {\cal E}_{[(2s)^{\frac{n}{2s}}]},\quad \alpha {\cal E}_{[(2s)^{\frac{4}{2s}}]}
\quad  \alpha {\cal E}_{[1,1,2]}.
$$
\no{This completes the proof.\quad\quad $\natural$}

According to the Rotation Embedding Scheme for orientable
embedding of a graph, First presented by Heffter in $1891$ and
formalized by Edmonds in $[17]$, each orientable complete map is
just the case of eliminating the sign + and - in our
representation for complete map. Whence,we have the following
result for an automorphism of orientable surfaces, which is
similar to the Theorem $3.2.1$.

\vskip 4mm

\no{\bf Theorem $3.2.2$} {\it All orientation-preserving
automorphisms of orientable complete maps of order$\geq 4$ are
extended actions of elements in}

$${\cal E}_{[s^{\frac{n}{s}}]}, \quad {\cal E}_{[1,s^{\frac{n-1}{s}}]} $$

\no{\it and all orientation-reversing automorphisms of orientable
complete maps of order$\geq 4$ are extended actions of elements in
}
$$
\alpha {\cal E}_{[(2s)^{\frac{n}{2s}}]},\quad \alpha {\cal
E}_{[(2s)^{\frac{4}{2s}}]}, \quad\alpha {\cal E}_{[1,1,2]},
$$
\no{\it where,${\cal E}_{\theta}$ denotes the conjugatcy class
containing $\theta$ in $S_n$.}

\vskip 3mm

{\it Proof}  The proof is similar to that of the Theorem $3.2.1$.
For completion, we only need to construct orientable maps $M_i^O ,
i=1,2,3,4$, to replace the non-orientable maps $M_i, i=1,2,3,4$ in
the proof of the Theorem $3.2.1$.

In fact, for orientation-preserving case, we only need to take
$M_1^O$, $M_2^O$ to be the resultant maps eliminating the sign +
and - in $M_1$, $M_2$ constructed in the proof of the Theorem
$3.2.1$.

For  the orientation-reversing case, we take  $M_3^O=
(E(K_n)_{\alpha ,\beta},{\cal P}_3^O)$ with

$$
{\cal P}_3=\prod\limits_{i\in \{1,2,\cdots,n\}}(C(i)),
$$

\no{where, if $i\equiv 1(mod 2)$, then}

$$C(i)=(i^{1},i^{k+1},\cdots,i^{n-k+1},i^{2},\cdots,i^{n-k+2},\cdots,i^{i*},\cdots,i^{k},i^{2k},\cdots,i^{n}),$$

\no{and if $i\equiv 0(mod 2)$, then}

$$C(i)=(i^{1},i^{k+1},\cdots,i^{n-k+1},i^{2},\cdots,i^{n-k+2},\cdots,i^{i*},\cdots,i^{k},i^{2k},\cdots,i^{n})^{-1},$$

\no{where $i^{i*}$ denotes the empty position and $M_4^O
=(E(K_4)_{\alpha ,\beta}, {\cal P}_4)$ with}

$$
{\cal P}_4= (1^2,1^3,1^4)(2^1,2^3,2^4)(3^1,3^4,3^2)(4^1,4^2,4^3).
$$

It can be shown that $(M_i^O)^{\xi^*\alpha}=M_i^O$ for $i=1,2,3$
and $4$. \quad\quad $\natural$

\vskip 10mm

{\bf\S 3. The automorphisms of a semi-regular graph on surfaces}

\vskip 8mm

A  graph is called a {\it semi-regular graph} if it is simple and
its automorphism group action on its ordered pair of adjacent
vertices is fixed-free, which is considered in $[43], [50]$ for
enumeration its non-equivalent embeddings on surfaces. A map
underlying a semi-regular graph is called a {\it semi-regular
map}. We determine all automorphisms of maps underlying a
semi-regular graph in this section.

Comparing with the Theorem $3.1.2$, we get a necessary and
sufficient condition for an automorphism of a graph being that of
a map.

\vskip 4mm

\no{{\bf Theorem $3.3.1$} {\it For a connected graph $\Gamma$, an
automorphism $\xi\in {\rm Aut} \Gamma$ is an
orientation-preserving automorphism of a non-orientable map
underlying the graph $\Gamma$ iff $\xi$ is semi-regular acting on
its ordered pairs of adjacent vertices. }}

\vskip 3mm

{\it Proof} According to the Lemma $2.2.1$, if $\xi\in{\rm Aut}
\Gamma$ is an orientation-preserving automorphism of a map $M$
underlying graph $\Gamma$, then $\xi$ is semi-regular acting on
its ordered pairs of adjacent vertices.

Now assume that $\xi\in {\rm Aut} \Gamma$ is semi-regular acting
on its ordered pairs of adjacent vertices. Denote by
$\xi|_{V(\Gamma)}$, $\xi|_{{E(\Gamma)}_{\beta}}$ the action of
$\xi$ on $V(\Gamma)$ and on its ordered pairs of adjacent
vertices,respectively. By the given condition, we can assume that

$$
\xi|_{V(\Gamma)}=(a,b,\cdots,c)\cdots (g,h,\cdots,k)\cdots
(x,y,\cdots,z)
$$

\no{and}

$$\xi|_{{E(\Gamma)}_{\beta}}= C_1\cdots C_i\cdots C_m,$$

\no{where,\ $s_a|C(a)| = \cdots = s_g|C(g)| = \cdots =
s_x|C(x)|$,\ and C(g) denotes the cycle containing $g$ in
$\xi|_{V(\Gamma)}$ and}

\vskip 2mm
$C_1=(a^1,b^1,\cdots,c^1,a^2,b^2,\cdots,c^2,\cdots,a^{s_a},b^{s_a},\cdots,c^{s_a}),$
\vskip 2mm
$\cdots\cdots\cdots\cdots\cdots\cdots\cdots\cdots\cdots\cdots\cdots\cdots\cdots\cdots\cdots\cdots\cdots,$
\vskip 2mm
$C_i=(g^1,h^1,\cdots,k^1,g^2,k^2,\cdots,k^2,\cdots,g^{s_g},h^{s_g},\cdots,k^{s_g}),$
\vskip 2mm
$\cdots\cdots\cdots\cdots\cdots\cdots\cdots\cdots\cdots\cdots\cdots\cdots\cdots\cdots\cdots\cdots\cdots,$
\vskip 2mm
$C_m=(x^1,y^1,\cdots,z^1,\cdots,x^2,y^2,\cdots,z^2,\cdots,x^{s_x},y^{s_x},\cdots,z^{s_x}).$

\vskip 2mm

Now for $\forall\xi, \xi\in {\rm Aut} \Gamma$. We construct a
stable map $M = ({\mathcal X}_{\alpha,\beta},{\mathcal P})$ under
the action of $\xi$ as follows.

 $$  X =E(\Gamma)$$

\no{and}

$$
{\mathcal P}=\prod\limits_{g\in T_{\xi}^V}\prod\limits_{x\in C(g)}
(C_x)(\alpha C_x^{-1}).
$$

Assume that $u=\xi^f (g)$, and

$$N_{\Gamma} (g)=\{g^{z_1},g^{z_2},\cdots,g^{z_l}\}.$$

Obviously, all degrees of vertices in $C(g)$ are same. Notices
that $\xi|_{N_{\Gamma}(g)}$ is circular acting on $N_{\Gamma} (g)$
by the Theorem $3.1.2$.  Whence, it is semi-regular acting on
$N_{\Gamma} (g)$. Without loss of generality, we can assume that

$$
\xi|_{N_{\Gamma}(g)}= (g^{z_1},g^{z_2},\cdots,g^{z_s})
(g^{z_{s+1}},g^{z_{s+2}},\cdots,g^{z_{2s}})\cdots
(g^{z_{(k-1)s+1}},g^{z_{(k-1)s+2}},\cdots,g^{z_{ks}}),
$$

\no{where, $l=ks$. Choose}

$$
C_g
=(g^{z_1+},g^{z_{s+1}+},\cdots,g^{z_{(k-1)s+1}+},g^{z_2+},g^{z_{s+2}+},\cdots,
g^{z_s+},g^{z_{2s}},\cdots,g^{z_{ks}+}).
$$

\no{Then,}

$$
C_x
=(x^{z_1+},x^{z_{s+1}+},\cdots,x^{z_{(k-1)s+1}+},x^{z_2+},x^{z_{s+2}+},\cdots,
x^{z_s+},x^{z_{2s}},\cdots,x^{z_{ks}+}),
$$

\no{where,}

$$x^{z_i+}=\xi^f (g^{z_i+}),$$

\no{for $i=1,2,\cdots,ks.$ and}

$$
\alpha C_x^{-1}=(\alpha x^{z_1+},\alpha x^{z_{s+1}+},\cdots,\alpha
x^{z_{(k-1)s+1}+}, \alpha x^{z_s+},\alpha x^{z_{2s}},\cdots,\alpha
x^{z_{ks}+}).
$$

Immediately, we get that $M^{\xi}=\xi M \xi^{-1} =M$ by this
construction. Whence, $\xi$ is an orientation-preserving
automorphism of the map $M$. \quad\quad $\natural$

By the Rotation Embedding Scheme, eliminating $\alpha$ on each
quadricell in Tutte's representation of an embedding induces an
orientable embedding with the same underlying graph. Since an
automorphism of an embedding is commutative with $\alpha$ and
$\beta$, we get the following result for the orientable-preserving
automorphisms of the orientable maps underlying a semi-regular
graph.

\vskip 4mm

\no{{\bf Theorem $3.3.2$} {\it If $\Gamma$ is a connected
semi-regular graph, then for $\forall\xi\in {\rm Aut} \Gamma$,
$\xi$ is an orientation-preserving automorphism of orientable maps
underlying the graph $\Gamma$.}}

\vskip 3mm

According to the Theorem $3.3.1$ and $3.3.2$, if $\Gamma$ is
semi-regular, i.e., each automorphism acting on the ordered pairs
of adjacent vertices in $\Gamma$ is fixed-free, then every
automorphism of the graph $\Gamma$ is an orientation-preserving
automorphism of orientable maps and non-orientable maps underlying
the graph $\Gamma$. We restated this result as the following.

\vskip 4mm

\no{\bf Theorem $3.3.3$}\quad {\it If\ $\Gamma$ is a connected
semi-regular graph, then for $\forall\xi\in {\rm Aut} \Gamma$,
$\xi$ is an orientation-preserving automorphism of orientable maps
and non-orientable maps underlying the graph $\Gamma$.}

\vskip 3mm

Notice that if $\varsigma^*$ is an orientation-reversing
automorphisms of a map, then $\varsigma^*\alpha$ is an
orientation-preserving automorphism of the same map. By the Lemma
$3.2.4$, if $\tau$ is an automorphism of maps underlying a graph
$\Gamma$, then $\tau\alpha$ is an automorphism of maps underlying
this graph if and only if $\tau\equiv 0(mod 2)$. Whence, we have
the following result for the automorphisms of maps underlying a
semi-regular graph

\vskip 4mm

\no{\bf Theorem $3.3.4$} {\it Let $\Gamma$ be a semi-regular
graph. Then all the automorphisms of orientable maps underlying
the graph $\Gamma$ are

$$g|^{{\mathcal X}_{\alpha ,\beta}} \ and \ \alpha h|^{{\mathcal X}_{\alpha ,\beta}},
g,h\in {\rm Aut}\Gamma \ with \ h\equiv 0(mod2).$$

\no{and all the automorphisms of non-orientable maps underlying
the graph $\Gamma$ are also

$$g|^{{\mathcal X}_{\alpha ,\beta}} \ and \ \alpha h|^{{\mathcal X}_{\alpha ,\beta}},
g,h\in {\rm Aut}\Gamma \ with \ h\equiv 0(mod2).$$}}

Theorem $3.3.4$ will be used in the Chapter $4$ for the
enumeration of unrooted maps on surfaces underlying a semi-regular
graph.

\vskip 10mm

{\bf\S 4. The automorphisms of one face maps}

\vskip 8mm

There is a well-know result for the automorphisms of a map and its
dual in topological graph theory, i.e., the automorphism group of
a map is same as its dual. Therefore, for determining the
automorphisms of one face maps, we can determine them by the
automorphisms of a bouquet $B_n$ on surfaces.

A map underlying a graph $B_n, n\geq 1$ has the form ${\mathcal
B}_n=({\cal X}_{\alpha,\beta },{\mathcal P}_n)$ with
$X=E(B_n)=\{e_1,e_2,\cdots,e_n \}$ and

$${\mathcal P}_n=(x_1,x_2,\cdots,x_{2n})(\alpha x_1,\alpha x_{2n},\cdots,x_2)$$

\no{where, $x_i\in X,\beta X$ or $\alpha\beta X$ and satisfying
the condition $(Ci)$ and $(Cii)$ in the Section $2.2$ of Chapter
$1$. For a given bouquet $B_n$ with $n$ edges, its semi-arc
automorphism group is}

$$Aut_{\frac{1}{2}}B_n = S_n[S_2].$$

\no{Form group theory, we know that each element in $S_n[S_2]$ can
be represented by $(g;h_1,h_2,\cdots, h_n)$ with $g\in S_n$ and
$h_i\in S_2=\{1,\alpha\beta\}$ for $i=1,2,\cdots,n$. The action of
$(g;h_1,h_2,\cdots, h_n)$ on a map ${\cal B}_n$ underlying a graph
$B_n$ by the following rule:}

\vskip 2mm {\it if $x\in \{e_i,\alpha e_i,\beta e_i,\alpha\beta
e_i\}$, then $(g;h_1,h_2,\cdots, h_n)(x)=g(h_i(x))$}. \vskip 2mm

\no{For example, if $h_1=\alpha\beta$, then,$(g;h_1,h_2,\cdots,
h_n)(e_1)=\alpha\beta g(e_1)$, $(g;h_1,h_2,\cdots, h_n)(\alpha
e_1)$ $=\beta g(e_1)$, $(g;h_1,h_2,\cdots, h_n)(\beta e_1)=\alpha
g(e_1)$ and $(g;h_1,h_2,\cdots, h_n)(\alpha\beta e_1)= g(e_1)$.}

The following result for automorphisms of a map underlying a graph
$B_n$ is obvious.

\vskip 4mm

\no{{\bf Lemma $3.4.1$} {\it Let $(g;h_1,h_2,\cdots, h_n)$ be an
automorphism of a map ${\mathcal B}_n$ underlying a graph $B_n$.
Then

$$
(g;h_1,h_2,\cdots, h_n)=(x_1,x_2,...,x_{2n})^k,
$$

\no{and if $(g;h_1,h_2,\cdots, h_n)\alpha$ is an automorphism of a
map ${\mathcal B}_n$, then }

$$
(g;h_1,h_2,\cdots, h_n)\alpha =(x_1,x_2,\cdots,x_{2n})^k
$$

\no{for some integer $k,1\leq k\leq n$. Where,
$x_i\in\{e_1,e_2,\cdots,e_n\}, i=1,2,\cdots ,2n$ and $x_i\not=x_j$
if $i\not=j.$}}

\vskip 3mm

Analyzing the structure of elements in the group $S_n[S_2]$, we
get the automorphisms of maps underlying a graph $B_n$ by the
Theorem $3.3.1$ and $3.3.2$ as follows.

\vskip 4mm

\no{\bf Theorem $3.4.1$} {\it Let $B_n$ be a bouquet with $n$
edges $1,2,\cdots ,n$. Then the automorphisms $(g;h_1,h_2,\cdots
,h_n)$ of orientable maps underlying a $B_n,n\geq 1$, are
respective }

($O1$)\quad $g\in {\mathcal E}_{[k^{\frac{n}{k}}]}, \ h_i=1,
i=1,2,\cdots ,n; $

($O2$)\quad $g\in  {\mathcal E}_{[k^{\frac{n}{k}}]} \ and \ if \
g=\prod\limits_{i=1}^{n/k}(i_1,i_2,\cdots i_k),$ $ \ where \
i_j\in\{1,2,\cdots ,n\},  n/k\equiv 0(mod 2), \ then \ h_{i_1}=
(1,\alpha\beta), i=1,2,\cdots ,\frac{n}{k} \ and \ h_{i_j}=1 \ for
\ j\geq 2;$

($O3$)\quad $g\in {\mathcal E}_{[k^{2s},(2k)^{\frac{n-2ks}{2k}}]}
\ and \ if \ g=\prod\limits_{i=1}^{2s}(i_1,i_2,\cdots
i_k)\prod\limits_{j=1}^{(n-2ks)/2k}(e_{j_1},e_{j_2},\cdots
,e_{j_{2k}}),$ $\ where \ i_j,e_{j_t}\in\{1,2,\cdots ,n\}, \ then
\ h_{i_1}= (1,\alpha\beta), \ i=1,2,\cdots ,s, h_{i_l}=1 \ for \
l\geq 2 \ and \ h_{j_t}=1 \ for \ t=1,2, \cdots , 2k$

\no{\it and the automorphisms $(g;h_1,h_2,\cdots ,h_n)$ of
non-orientable maps underlying a $B_n,n$ $\geq 1$, are respective
}

($N1$)\quad $g\in {\mathcal E}_{[k^{\frac{n}{k}}]}, \ h_i=1,
i=1,2,\cdots ,n; $

($N2$)\quad $g\in  {\mathcal E}_{[k^{\frac{n}{k}}]} \ and \ if \
g=\prod\limits_{i=1}^{n/k}(i_1,i_2,\cdots i_k),$ $ \ where \
i_j\in\{1,2,\cdots ,n\}, n/k\equiv 0(mod 2), \ then \ h_{i_1} $ $=
(1,\alpha\beta), (1,\beta ) \ with \ at \ least \ one \
h_{{i_0}_1}(1,\beta ), i=1,2,\cdots ,\frac{n}{k} \ and \ h_{i_j}=1
\ for \ j\geq 2;$

($N3$)\quad $g\in {\mathcal E}_{[k^{2s},(2k)^{\frac{n-2ks}{2k}}]}
\ and \ if \ g=\prod\limits_{i=1}^{2s}(i_1,i_2,\cdots
i_k)\prod\limits_{j=1}^{(n-2ks)/2k}(e_{j_1},e_{j_2},\cdots
,e_{j_{2k}}),$ $ \ where \ i_j,e_{j_t}\in\{1,2,\cdots ,n\}, \ then
\ h_{i_1}= (1,\alpha\beta),(1,\beta ) \ with \ at \ least \ one \
h_{{i_0}_1}=(1,\beta ), i=1,2,\cdots ,s, h_{i_l}=1 \ for \ l\geq 2
\ and \ h_{j_t}=1 \ for \ t=1,2, \cdots , 2k,$

\no{\it where, ${\mathcal E}_{\theta}$ denotes the conjugacy class
in the symmetry group $S_n$ containing the element $\theta$.}

\vskip 3mm

{\it Proof} By the structure of the group $S_n[S_2]$, it is clear
that the elements in the cases $(i), (ii)$ and $(iii)$ are its all
semi-regular elements. We only need to construct an orientable or
non-orientable map ${\mathcal B}_n=({\cal X}_{\alpha,\beta
},{\mathcal P}_n)$ underlying $B_n$ stable under the action of an
element in each case.

\vskip 3mm

$(i)$\quad $g=\prod\limits_{i=1}^{n/k}(i_1,i_2,\cdots i_k)$ and
$h_i=1, i=1,2,\cdots ,n$, where  $i_j\in\{1,2,\cdots ,n\}$.

\vskip 2mm

Choose

$${\mathcal X}_{\alpha
,\beta}^1=\bigcup\limits_{i=1}^{n/k}K\{i_1,i_2,\cdots ,i_k\},$$

\no{where $K=\{1,\alpha ,\beta ,\alpha\beta\}$ and}

$${\mathcal P}_n^1=\ C_1 (\alpha C_1^{-1}\alpha^{-1})$$

\no{with}

\begin{eqnarray*}
C_1=(&1_1,& 2_1,\cdots ,(\frac{n}{k})_1,\alpha\beta
1_1,\alpha\beta 2_1,\cdots ,\alpha\beta (\frac{n}{k})_1, 1_2,2_2,
\cdots ,(\frac{n}{k})_2,\\
&\alpha\beta 1_2,& \alpha\beta 2_2, \cdots ,\alpha\beta
(\frac{n}{k})_2, \cdots ,1_k,2_k,\cdots
,(\frac{n}{k})_k,\alpha\beta 1_k,\alpha\beta 1_k,\cdots
,\alpha\beta (\frac{n}{k})_k).
\end{eqnarray*}

\no{Then the map ${\mathcal B}_n^1=({\cal X}_{\alpha,\beta
}^1,{\mathcal P}_n^1)$ is an orientable map underlying the graph
$B_n$ and stable under the action of $(g;h_1,h_2,\cdots ,h_n)$.}

For the non-orientable case, we can chose

\begin{eqnarray*}
C_1 = (&1_1,& 2_1,\cdots ,(\frac{n}{k})_1,\beta 1_1,\beta
2_1,\cdots ,\beta (\frac{n}{k})_1, 1_2,2_2,\cdots
,(\frac{n}{k})_2,\\
&\beta 1_2,& \beta 2_2, \cdots ,\beta (\frac{n}{k})_2, \cdots
,1_k,2_k,\cdots ,(\frac{n}{k})_k,\beta 1_k,\beta 1_k,\cdots ,\beta
(\frac{n}{k})_k).
\end{eqnarray*}

\no{Then the map ${\mathcal B}_n^1=({\cal X}_{\alpha,\beta
}^1,{\mathcal P}_n^1)$ is a non-orientable map underlying the
graph $B_n$ and stable under the action of $(g;h_1,h_2,\cdots
,h_n)$.}

\vskip 3mm

($ii$)\quad $g=\prod\limits_{i=1}^{n/k}(i_1,i_2,\cdots i_k)$,
$h_i=(1,\beta ) \ or \ (1,\alpha\beta), i=1,2,\cdots ,n, \
\frac{n}{k}\equiv 0(mod 2)$, where $i_j\in\{1,2,\cdots ,n\}$.

\vskip 2mm

If $h_{i_1}=(1,\alpha\beta)$ for $i=1,2,\cdots ,\frac{n}{k}$ and
$h_{i_t}=1$ for $t\geq 2$. Then

$$(g;h_1,h_2,\cdots ,h_n)=\prod\limits_{i=1}^{n/k}(i_1,\alpha\beta i_2,\cdots
\alpha\beta i_k,\alpha\beta i_1,i_2,\cdots ,i_k).$$

Similar to the case $(i)$, we take ${\mathcal X}_{\alpha
,\beta}^2={\mathcal X}_{\alpha ,\beta}^1$ and

$${\mathcal P}_n^2=\ C_2 (\alpha C_2^{-1}\alpha^{-1})$$

\no{with}

\begin{eqnarray*}
C_2=(&1_1,& 2_1,\cdots ,(\frac{n}{k})_1,\alpha\beta
1_2,\alpha\beta 2_2,\cdots ,\alpha\beta
(\frac{n}{k})_2,\alpha\beta 1_k,\alpha\beta 2_k,\\
&\cdots& ,\alpha\beta (\frac{n}{k})_k, \alpha\beta 1_1,\alpha\beta
2_1, \cdots ,\alpha\beta (\frac{n}{k})_1, 1_2,2_2,\cdots
,(\frac{n}{k})_2, \cdots
 ,1_k,2_k,\cdots
,(\frac{n}{k})_k).
\end{eqnarray*}

\no{Then the map ${\mathcal B}_n^2=({\cal X}_{\alpha,\beta
}^2,{\mathcal P}_n^2)$ is an orientable map underlying the graph
$B_n$ and stable under the action of $(g;h_1,h_2,\cdots ,h_n)$.
For the non-orientable case, the construction is similar. It only
need to replace each element $\alpha\beta i_j$ by $\beta i_j$ in
the construction of the orientable case if $h_{i_j}=(1,\beta )$.}

\vskip 3mm

($iii$)\quad  $ g=\prod\limits_{i=1}^{2s}(i_1,i_2,\cdots
i_k)\prod\limits_{j=1}^{(n-2ks)/2k}(e_{j_1},e_{j_2},\cdots
,e_{j_{2k}})$ and $ h_{i_1}= (1,\alpha\beta), \ i=1,2,\cdots ,s,
h_{i_l}=1 \ for \ l\geq 2 \ and \ h_{j_t}=1 \ for \ t=1,2, \cdots
, 2k$

\vskip 2mm

In this case, we know that

$$(g;h_1,h_2,\cdots ,h_n)=\prod\limits_{i=1}^{s}(i_1,\alpha\beta i_2,\cdots
\alpha\beta i_k,\alpha\beta i_1,i_2,\cdots
,i_k)\prod\limits_{j=1}^{(n-2ks)/2k}(e_{j_1},e_{j_2},\cdots
,e_{j_{2k}}).$$

Denote by $p$ the number $(n-2ks)/2k$.  We construct an orientable
map ${\mathcal B}_n^3=({\cal X}_{\alpha,\beta }^3,{\mathcal
P}_n^3)$ underlying $B_n$ stable under the action of
$(g;h_1,h_2,\cdots ,h_n)$ as follows.

Take ${\mathcal X}_{\alpha ,\beta}^3={\mathcal X}_{\alpha
,\beta}^1$ and

$${\mathcal P}_n^3=\ C_3 (\alpha C_3^{-1}\alpha^{-1})$$

\no{with}

\begin{eqnarray*}
C_3= (&1_1&, 2_1,\cdots ,s_1, e_{1_1},e_{2_1},\cdots ,e_{p_1},
\alpha\beta 1_2,\alpha\beta 2_2,\cdots ,\alpha\beta s_2,\\
&e_{1_2},& e_{2_2},\cdots ,e_{p_2},\cdots ,\alpha\beta
1_k,\alpha\beta 2_k,\cdots ,\alpha\beta s_k,e_{1_k},e_{2_k},\cdots
,\\
&e_{p_k},& \alpha\beta 1_1, \alpha\beta 2_1,\cdots ,\alpha\beta
s_1,e_{1_{k+1}},e_{2_{k+1}},\cdots ,e_{p_{k+1}},1_2,2_2,\cdots ,\\
&s_2,& e_{1_{k+2}}, e_{2_{k+2}},\cdots ,e_{p_{k+2}},\cdots
,1_k,2_k,\cdots ,s_k,e_{1_{2k}},e_{2_{2k}},\cdots ,e_{p_{2k}}).
\end{eqnarray*}

\no{Then the map ${\mathcal B}_n^3=({\cal X}_{\alpha,\beta
}^3,{\mathcal P}_n^3)$ is an orientable map underlying the graph
$B_n$ and stable under the action of $(g;h_1,h_2,\cdots ,h_n)$.

Similarly, replacing each element $\alpha\beta i_j$ by $\beta i_j$
in the construction of the orientable case if $h_{i_j}=(1,\beta
)$, a non-orientable map underlying the graph $B_n$ and stable
under the action of $(g;h_1,h_2,\cdots ,h_n)$ can be also
constructed. This completes the proof.\quad\quad $\natural$

We use the Lemma $3.4.1$ for the enumeration of unrooted one face
maps on surfaces in the next chapter.

\newpage

%%%%%%%%%%%%%%%Headings%%%%%%%%%%%%%%%%%%%%%%%%%%%%%%%%%%%%%%%%%
\thispagestyle{empty} \pagestyle{myheadings} \topmargin 5mm
\headheight 8mm \headsep 10mm

\markright {\scriptsize  Chapter $4$\quad Application to the
Enumeration of Unrooted Maps and $s$-Manifolds}
%%%%%%%%%%%%%%%%%%%%%%%%%%%%%%%%%%%%%%%%%%%%%%%%%%%%%%%%%%%%%%%%

\vskip 30mm

\no{\bf\large Chapter $4$\quad Application to the Enumeration of
Unrooted Maps and $s$-manifolds}

\vskip 10mm

All the results gotten in the Chapter $3$ is more useful for the
enumeration of unrooted maps on surfaces underlying a graph.
According to the theory in Chapter $1$, the enumeration of
unrooted maps on surfaces underlying a graph can be carried out by
the following programming:

\vskip 3mm

STEP 1. Determine all automorphisms of maps underlying a graph;

\vskip 2mm

STEP 2. Calculate the the fixing set $Fix(\varsigma)$ for each
automorphism $\varsigma$ of a map;

\vskip 2mm

STEP 3. Enumerate the unrooted maps on surfaces underlying a graph
by Corollary $1.3.1$.

\vskip 2mm

\no{The advantage of this approach is its independence of the
orientability, which enables us to enumerate orientable or
non-orientable maps on surafces. We present the enumeration
results by this programming in this chapter.}

\vskip 10mm

{\bf\S 1. The enumeration of unrooted complete maps on surfaces}

\vskip 8mm

We first consider a permutation and its stabilizer . A permutation
with the following form $(x_1,x_2,...,x_n)( \alpha x_n, \alpha
x_2,...,\alpha x_1)$ is called a {\it permutation pair}. The
following result is obvious.

\vskip 4mm

\no{{\bf Lemma $4.1.1$} {\it Let $g$ be a permutation on the set
$\Omega=\{x_1,x_2,...,x_n\}$ such that  $g\alpha=\alpha g $. If}}

$$g(x_1,x_2,...,x_n)( \alpha x_n, \alpha x_{n-1},..., \alpha x_1)g^{-1}=
(x_1,x_2,...,x_n)( \alpha x_n, \alpha x_{n-1},..., \alpha x_1),$$

\no{\it then}

$$
g=(x_1,x_2,...,x_n)^k
$$

\no{\it and if}

$$g\alpha (x_1,x_2,...,x_n)( \alpha x_n, \alpha x_{n-1},..., \alpha x_1)(g\alpha )^{-1}=
(x_1,x_2,...,x_n)( \alpha x_n, \alpha x_{n-1},..., \alpha x_1),$$

\no{\it then}

$$
g\alpha =(\alpha x_n, \alpha x_{n-1},..., \alpha x_1)^k
$$

\no{for some integer $k,1\leq k\leq n$.}

\vskip 4mm

\no{{\bf Lemma $4.1.2$} {\it For each permutation $g,g\in
{\mathcal E}_{[k^{\frac{n}{k}}]}$ satisfying $g\alpha=\alpha g$ on
the set $\Omega=\{x_1,x_2,...,x_n\}$, the number of stable
permutation pairs in $\Omega$ under the action of $g$ or $g\alpha$
is}}

$$
\frac{2\phi (k)(n-1)!}{|{\mathcal E}_{[k^{\frac{n}{k}}]}|},
$$

\no{\it  where  $\phi (k)$ denotes the Euler function.}

\vskip 3mm

{\it Proof}  Denote the number of stable pair permutations under
the action of $g$ or $g\alpha$ by $n(g)$ and $\mathcal C$ the set
of pair permutations. Define the set $A=\{(g,C)| g\in {\mathcal
E}_{[k^{\frac{n}{k}}]}, C\in{\mathcal C}\quad and \quad C^g=C \ or
\ $ $C^{g\alpha}=C\}$. Clearly, for $\forall g_1,g_2\in{ {\mathcal
E}_{[k^{\frac{n}{k}}]}}$, we have $n(g_1)=n(g_2)$. Whence, we get
that

$$
|A| = | {\mathcal E}_{[k^{\frac{n}{k}}]}| n(g) . \hspace{30mm}
(4.1.1)
$$

On the other hand, by the Lemma $4.1.1$, for any permutation pair
$C=(x_1,x_2,...,x_n)$ $( \alpha x_n, \alpha x_{n-1},..., \alpha
x_1)$, since $C$ is stable under the action of $g$, there must be
$g=(x_1,x_2,...,x_n)^l $ or $g\alpha =(\alpha x_n, \alpha
x_{n-1},..., \alpha x_1)^l$, where $l=s\frac{n}{k}, 1\leq s\leq k$
and $(s,k)=1$. Therefore, there are $2\phi(k)$ permutations in
${\mathcal E}_{[k^{\frac{n}{k}}]}$ acting on it stable. Whence, we
also have

$$
|A| = 2\phi (k) |{\mathcal C}|. \hspace{30mm}         (4.1.2)
$$

Combining ($4.1.1$) with ($4.1.2$), we get that

$$
n(g)=\frac{2\phi (k)|{\mathcal C}|}{|{\mathcal
E}_{[k^{\frac{n}{k}}]}|}= \frac{2\phi (k)(n-1)!}{ |{\mathcal
E}_{[k^{\frac{n}{k}}]}| }.\hspace{20mm}\natural
$$

Now we can enumerate the unrooted complete maps on surfaces.

\vskip 4mm

\no{{\bf Theorem $4.1.1$} {\it The number $n^{L}(K_n)$ of complete
maps of order $n\geq 5$ on surfaces is}}

$$
n^{L}(K_n)=\frac{1}{2}(\sum\limits_{k|n}+\sum\limits_{k|n, k\equiv
0(mod 2)})
 \frac{2^{\alpha (n,k)}(n-2)!^{\frac{n}{k}}}
{k^{\frac{n}{k}}(\frac{n}{k})!}+ \sum\limits_{k|(n-1),k\not=1}
\frac{\phi (k)2^{\beta (n,k)}(n-2)!^{\frac{n-1}{k}}}{n-1},
$$
\no{where, }

\[
\alpha (n,k)=\left\{\begin{array}{cc}
\frac{n(n-3)}{2k}, & {\rm if}\quad k\equiv 1(mod 2);\\
\frac{n(n-2)}{2k}, & {\rm if}\quad k\equiv 0(mod 2),
\end{array}
\right.
\]

\no{and}

\[
\beta (n,k)=\left\{\begin{array}{cc}
\frac{(n-1)(n-2)}{2k}, & {\rm if}\quad k\equiv 1(mod 2);\\
\frac{(n-1)(n-3)}{2k}, & {\rm if}\quad k\equiv 0(mod 2).
\end{array}
\right.
\]
\no{and $n^L(K_4)=11.$}

\vskip 3mm

{\it Proof}   According to ($1.3.3$) in the Corollary $1.3.1$ and
the Theorem $3.2.1$ for $n\geq 5$, we know that

\begin{eqnarray*}
n^{L}(K_n) &=& \frac{1}{2|{\rm AutK}_n|}
\times(\sum\limits_{g_1\in{{\mathcal E}_{[k^{\frac{n}{k}}]}}}|\Phi
(g_1)| + \sum\limits_{g_2\in{{\mathcal
E}_{[(2s)^{\frac{n}{2s}}]}}}|\Phi (g_2\alpha)|\\
&+& \sum\limits_{h\in{{\mathcal E}_{[1,k^{\frac{n-1}{k}}]}}}|\Phi (h)|) \\
&=& \frac{1}{2n!}\times (\sum\limits_{k|n}|{\mathcal
E}_{[k^{\frac{n}{k}}]}| |\Phi (g_1)|+ \sum\limits_{l|n,l\equiv
0(mod2)}|{\mathcal E}_{[l^{\frac{n}{l}}]}| |\Phi (g_2\alpha)|\\
&+& \sum\limits_{l|(n-1)}|{\mathcal E}_{[1,l^{\frac{n-1}{l}}]}|
|\Phi (h)|),
\end{eqnarray*}

\no{where, $g_1\in{{\mathcal
E}_{[k^{\frac{n}{k}}]}},g_2\in{{\mathcal E}_{[l^{\frac{n}{l}}]}}$
and $h\in{{\mathcal E}_{[1,k^{\frac{n-1}{k}}]}}$ are three chosen
elements. }

Without loss of generality, we assume that an element
$g,g\in{{\mathcal E}_{[k^{\frac{n}{k}}]}}$ has the following cycle
decomposition.

$$g = (1,2,...,k)(k+1,k+2,...,2k)...((\frac{n}{k}-1)k+1,(\frac{n}{k}-1)k+2,...,n)$$

and

$$
{\mathcal P}={\prod}_1\times{\prod}_2,
$$

where,

$$
{\prod}_1=(1^{i_{21}},1^{i_{31}},...,1^{i_{n1}})(2^{i_{12}},2^{i_{32}},...,2^{i_{n2}})...
(n^{i_{1n}},n^{i_{2n}},...,n^{i_{(n-1)n}}),
$$

and

$${\prod}_2=\alpha ({{\prod}_1}^{-1})\alpha^{-1},$$

\no{being a complete map which is stable under the action of $g$,
where $s_{ij}\in \{k+,k-| k= 1,2,...n \}$.}

Notice that the quadricells adjacent to the vertex "$1$" can make
$2^{n-2}(n-2)! $different pair permutations and for each chosen
pair permutation, the pair permutations adjacent to the vertices
$2,3,...,k$ are uniquely determined since $\mathcal P$  is stable
under the action of $g$.

Similarly, for each given pair permutation adjacent to the vertex
$k+1,2k+1,..., (\frac{n}{k}-1)k+1$, the pair permutations adjacent
to $k+2,k+3,...,2k$ and $2k+2,2k+3,...,3k$ and,...,and
$(\frac{n}{k}-1)k+2,(\frac{n}{k}-1)k+3,...n$ are also uniquely
determined because $\mathcal P$ is stable under the action of $g$.

Now for an orientable embedding $M_1$ of $K_n$, all the induced
embeddings by exchanging two sides of some edges and retaining the
others unchanged in $M_1$ are the same as $M_1$ by the definition
of maps. Whence, the number of different stable embeddings under
the action of $g$ gotten by exchanging $x$ and $\alpha x$ in $M_1$
for $x\in U, U\subset{\mathcal X}_{\beta}$, where ${\mathcal
X}_{\beta}= \bigcup\limits_{x\in E(K_n)} \{x,\beta x \}$ , is
$2^{g(\varepsilon)-\frac{n}{k}}$, where $g(\varepsilon)$ is the
number of orbits of $E(K_n)$ under the action of $g$ and we
substract $\frac{n}{k}$ because we can chosen
$1^{2+},k+1^{1+},2k+1^{1+},\cdots,n-k+1^{1+}$ first in our
enumeration.

Notice that the length of each orbit under the action of $g$ is
$k$ for $\forall x\in E(K_n)$ if $k$ is odd and is $\frac{k}{2}$
for $x=i^{i+\frac{k}{2}}, i=1,k+1,\cdots, n-k+1$, or $k$ for all
other edges if $k$ is even. Therefore, we get that

\[
g(\varepsilon)=\left\{\begin{array}{cc}
\frac{\varepsilon (K_n)}{k}, & {\rm if}\quad k\equiv 1(mod 2);\\
\frac{\varepsilon (K_n)-\frac{n}{2}}{k}, & {\rm if}\quad k\equiv
0(mod 2).
\end{array}
\right.
\]

\no{Whence, we have that}

\[
\alpha (n,k)=g(\varepsilon) - \frac{n}{k}=\left\{\begin{array}{cc}
\frac{n(n-3)}{2k}, & {\rm if}\quad k\equiv 1(mod 2);\\
\frac{n(n-2)}{2k}, & {\rm if}\quad k\equiv 0(mod 2),
\end{array}
\right.
\]

\no{and}

$$
|\Phi (g)|=2^{\alpha (n,k)}(n-2)!^{\frac{n}{k}},\hspace{30mm}
(4.1.3)
$$

Similarly, if $k\equiv 0(mod2)$, we get also that

$$
|\Phi (g\alpha)|=2^{\alpha (n,k)}(n-2)!^{\frac{n}{k}}
\hspace{30mm} (4.1.4)
$$

\no{for an chosen element $g$,$g\in {\mathcal
E}_{[k^{\frac{n}{k}}]}.$}

Now for $\forall h\in{{\mathcal E}_{[1,k^{\frac{n-1}{k}}]}} $,
without loss of generality, we assume that $h = (1,2,...,k)$
$(k+1,k+2,...,2k)...
((\frac{n-1}{k}-1)k+1,(\frac{n-1}{k}-1)k+2,...,(n-1))(n)$. Then
the above statement is also true for the complete graph $K_{n-1}$
with the vertices $1,2,\cdots,n-1$. Notice that the quadricells
$n^{1+},n^{2+},\cdots,n^{n-1+}$ can be chosen first in our
enumeration and they are not belong to the graph $K_{n-1}$.
According to the Lemma $4.1.2$, we get that

$$
|\Phi (h)|=2^{\beta (n,k)}(n-2)!^{\frac{n-1}{k}}\times\frac{2\phi
(k)(n-2)!} {| {\mathcal E}_{[1,k^{\frac{n-1}{k}}]}|},\hspace{10mm}
(4.1.5)
$$

\no{Where}

\[
\beta (n,k)=h(\varepsilon)=\left\{\begin{array}{cc}
\frac{\varepsilon
(K_{n-1})}{k}-\frac{n-1}{k}=\frac{(n-1)(n-4)}{2k},
& {\rm if}\quad k\equiv 1(mod 2);\\
\frac{\varepsilon
(K_{n-1})}{k}-\frac{n-1}{k}=\frac{(n-1)(n-3)}{2k}, & {\rm if}\quad
k\equiv 0(mod 2).
\end{array}
\right.
\]

Combining $(4.1.3)-(4.1.5)$, we get that

\begin{eqnarray*}
n^{L}(K_n) &=&  \frac{1}{2n!}\times (\sum\limits_{k|n}|{\mathcal
E}_{[k^{\frac{n}{k}}]}| |\Phi (g_0)|+\sum\limits_{l|n,l\equiv
0(mod2)}|{\mathcal E}_{[l^{\frac{n}{l}}]}| |\Phi (g_1\alpha)|\\
&+& \sum\limits_{l|(n-1)}|{\mathcal E}_{[1,l^{\frac{n-1}{l}}]}|
|\Phi
(h)|)\\
 &=& \frac{1}{2n!}\times ( \sum\limits_{k|n} \frac{n! 2^{\alpha
(n,k)}(n-2)!^{\frac{n}{k}}}{k^{\frac{n}{k}}(\frac{n}{k})!}+
\sum\limits_{k|n,k\equiv 0(mod2)}
\frac{n! 2^{\alpha (n,k)}(n-2)!^{\frac{n}{k}}}{k^{\frac{n}{k}}(\frac{n}{k})!}\\
&+&\sum\limits_{k|(n-1),k\not =1}
\frac{n!}{k^{\frac{n-1}{k}}(\frac{n-1}{k})!} \times\frac{2\phi
(k)(n-2)! 2^{\beta (n,k)}(n-2)!^{\frac{n-1}{k}}}
{\frac{(n-1)!}{k^{\frac{n-1}{k}}(\frac{n-1}{k})!}})\\
&=& \frac{1}{2}(\sum\limits_{k|n}+\sum\limits_{k|n,k\equiv
0(mod2)}) \frac{2^{\alpha (n,k)}(n-2)!^{\frac{n}{k}}}
{k^{\frac{n}{k}}(\frac{n}{k})!}+ \sum\limits_{k|(n-1),k\not=1}
\frac{\phi (k)2^{\beta (n,k)}(n-2)!^{\frac{n-1}{k}}}{n-1}.
\end{eqnarray*}

 For $n=4$, similar calculation shows that $n^L(K_4)=11$
by consider the fixing set of permutations in ${\mathcal
E}_{[s^{\frac{4}{s}}]}$,${\mathcal E}_{[1,s^{\frac{3}{s}}]}$,
${\mathcal E}_{[(2s)^{\frac{4}{2s}}]}$,$\alpha{\mathcal
E}_{[(2s)^{\frac{4}{2s}}]}$ and $\alpha {\mathcal E}_{[1,1,2]}$.
\quad\quad $\natural$

For the orientable case, we get the number $n^{O}(K_n)$ of
orientable complete maps of order $n$ as follows.

\vskip 4mm

\no{{\bf Theorem $4.1.2$} {\it The number $n^O((K_n)$ of complete
maps of order $n\geq 5$ on orientable surfaces is}}

$$
n^{O}(K_n)=\frac{1}{2}(\sum\limits_{k|n}+\sum\limits_{k|n, k\equiv
0(mod 2)} ) \frac{(n-2)!^{\frac{n}{k}}}
{k^{\frac{n}{k}}(\frac{n}{k})!}+ \sum\limits_{k|(n-1),k\not=1}
\frac{\phi (k)(n-2)!^{\frac{n-1}{k}}}{n-1} .
$$

\no{and $n(K_4)=3.$}

\vskip 3mm

{\it Proof}  According to Tutte's algebraic representation of map,
a map $M=({\mathcal X}_{\alpha,\beta},{\mathcal P})$ is orientable
if and only if for $\forall x\in {\mathcal X}_{\alpha ,\beta}$,
$x$ and $\alpha\beta x$ are in a same orbit of ${\mathcal
X}_{\alpha ,\beta}$ under the action of the group
$\Psi_I=<\alpha\beta,{\mathcal P}>$. Now applying ($1.3.1$) in the
Corollary $1.3.1$ and the Theorem $3.2.1$, Similar to the proof of
the Theorem $4.1.1$, we get the number $n^O(K_n)$ for $n\geq 5$ as
follows}

$$
n^{O}(K_n)=\frac{1}{2}(\sum\limits_{k|n}+\sum\limits_{k|n, k\equiv
0(mod 2)} ) \frac{(n-2)!^{\frac{n}{k}}}
{k^{\frac{n}{k}}(\frac{n}{k})!}+ \sum\limits_{k|(n-1),k\not=1}
\frac{\phi (k)(n-2)!^{\frac{n-1}{k}}}{n-1} .
$$

\no{and for the complete graph $K_4$, calculation shows that
$n(K_4)=3. \quad\quad \natural$}

Notice that $n^O(K_n)+n^N(K_n)=n^L(K_n)$. Therefore, we can also
get the number $n^N(K_n)$ of unrooted complete maps of order $n$
on non-orientable surfaces by the Theorem $4.1.1$ and the Theorem
$4.1.2$ as follows.

\vskip 4mm

\no{{\bf Theorem $4.1.3$} {\it The number $n^N(K_n)$ of unrooted
complete maps of order $n, n\geq 5$ on non-orientable surfaces
is}}

\begin{eqnarray*}
n^{N}(K_n)&=& \frac{1}{2}(\sum\limits_{k|n}+\sum\limits_{k|n,
k\equiv 0(mod 2)} ) \frac{(2^{\alpha
(n,k)}-1)(n-2)!^{\frac{n}{k}}}
{k^{\frac{n}{k}}(\frac{n}{k})!}\\
&+& \sum\limits_{k|(n-1),k\not=1} \frac{\phi (k)(2^{\beta
(n,k)}-1)(n-2)!^{\frac{n-1}{k}}}{n-1},
\end{eqnarray*}

\no{and $n^N(K_4)=8$. Where, $\alpha (n,k)$ and $\beta(n,k)$ are
same as in Theorem $4.1.1$.}\vskip 3mm

\includegraphics[bb=10 10 200 400]{7mg2.eps}

\vskip 3mm \c{Fig.$4.1$}

For $n=5$, calculation shows that $n^L(K_5)=1080$ and $n^{O}(K_5)
=45$ based on the Theorem $4.1.1$ and $4.1.2$. For $n=4$, there
are 3 unrooted orientable maps and 8 non-orientable maps shown in
the Fig.$4.1$. All the 11 maps of $K_4$ on surfaces are
non-isomorphic.

Now consider the action of orientation-preserving automorphisms of
complete maps, determined in the Theorem $3.2.1$ on all orientable
embeddings of a complete graph of order $n$. Similar to the proof
of the Theorem $4.1.2$, we can get the number of non-equivalent
embeddings of a complete graph of order $n$, which has been done
in $[43]$ and the result gotten is same as the result of Mull et
al in $[54]$.

\vskip 10mm

{\bf\S 2. The enumeration of a semi-regular graph on surfaces}

\vskip 8mm

The non-equivalent embeddings of a semi-regular graph on surfaces
are enumerated in the reference $[50]$. Here, we enumerate the
unrooted maps underlying a semi-regular graph on orientable or
non-orientable surfaces.

The following lemma is implied in the proof of the Theorem $4.1$
in $[50]$.

\vskip 4mm

\no{\bf Lemma $4.2.1$} {\it Let $\Gamma =(V,E)$ be a semi-regular
graph. Then for $\xi\in {\rm Aut}\Gamma$}

$$
|\Phi^{O} (\xi)|=\prod\limits_{x\in T_{\xi}^V}
(\frac{d(x)}{o(\xi|_{N_{\Gamma}(x)})}-1)!
$$

\no{\it and}

$$
|\Phi^{L} (\xi)|= 2^{|T_{\xi}^E| -|T_{\xi}^V|} \prod\limits_{x\in
T_{\xi}^V} (\frac{d(x)}{o(\xi|_{N_{\Gamma}(x)})}-1)!,
$$

\no{\it where, $T_{\xi}^V, T_{\xi}^E$ are the representations of
orbits of $\xi$ acting on $v(\Gamma)$ and $E(\Gamma)$
,respectively and $\xi_{N_{\Gamma}(x)}$ the restriction of $\xi$
to $N_{\Gamma}(x)$.}

\vskip 3mm

According to the Corollary $1.3.1$, we get the following
enumeration results.

\vskip 4mm

\no{\bf Theorem $4.2.1$} {\it Let $\Gamma$ be a semi-regular
graph. Then the numbers of unrooted maps on orientable and
non-orientable surfaces underlying the graph $\Gamma$ are}

$$
n^O(\Gamma )=\frac{1}{|{\rm Aut}\Gamma|}(\sum\limits_{\xi \in {\rm
Aut}\Gamma}\lambda (\xi )\prod\limits_{x\in T_{\xi}^V}
(\frac{d(x)}{o(\xi|_{N_{\Gamma}(x)})}-1)!
$$

\no{and}

$$
n^N(\Gamma )=\frac{1}{|Aut \Gamma|}\times \sum\limits_{\xi\in Aut
\Gamma} (2^{|T_{\xi}^E| -|T_{\xi}^V|}-1)\lambda (\xi )
\prod\limits_{x\in T_{\xi}^V}
(\frac{d(x)}{o(\xi|_{N_{\Gamma}(x)})}-1)!,
$$

\no{\it where $\lambda (\xi )=1$ if $o(\xi )\equiv 0(mod 2)$ and
$\frac{1}{2}$, otherwise.}

\vskip 3mm

{\it Proof}  By the Corollary $1.3.1$, we know that

$$
n^{O}(\Gamma)=\frac{1}{2|\rm
Aut_{\frac{1}{2}}\Gamma|}\sum\limits_{g \in {\rm
Aut_{\frac{1}{2}}\Gamma}} |\Phi_1^{O} (g)|
$$

\no{and}

$$
n^{L}(\Gamma)=\frac{1}{2|\rm
Aut_{\frac{1}{2}}\Gamma|}\sum\limits_{g \in {\rm
Aut_{\frac{1}{2}}\Gamma}} |\Phi_1^{T} (g)|.
$$

According to the Theorem $3.3.4$, all the automorphisms of
orientable maps underlying the graph $\Gamma$ are

$$g|^{{\mathcal X}_{\alpha ,\beta}} \ {\rm and} \ \alpha h|^{{\mathcal X}_{\alpha ,\beta}},
g,h\in {\rm Aut}\Gamma \ {\rm with}  \ h\equiv 0(mod2).$$

\no{and all the automorphisms of non-orientable maps underlying
the graph $\Gamma$ are also

$$g|^{{\mathcal X}_{\alpha ,\beta}} \ {\rm and} \ \alpha h|^{{\mathcal X}_{\alpha ,\beta}},
g,h\in {\rm Aut}\Gamma \ {\rm with}  \ h\equiv 0(mod2).$$

Whence, we get the number of unrooted orientable maps by the Lemma
$4.2.1$ as follows.

\begin{eqnarray*}
n^{O}(\Gamma) &=& \frac{1}{2|{\rm Aut}\Gamma|}\sum\limits_{g \in
{\rm Aut}\Gamma} |\Phi_1^{O} (g)|\\
&=& \frac{1}{2|{\rm Aut}\Gamma|}\{(\sum\limits_{\xi \in {\rm
Aut}\Gamma}\prod\limits_{x\in T_{\xi}^V}
(\frac{d(x)}{o(\xi|_{N_{\Gamma}(x)})}-1)!\\
&+& \sum\limits_{\varsigma \in {\rm Aut}\Gamma ,o(\varsigma
)\equiv 0(mod 2)}\prod\limits_{x\in T_{\varsigma}^V}
(\frac{d(x)}{o(\varsigma|_{N_{\Gamma}(x)})}-1)!)\\
&=& \frac{1}{|{\rm Aut}\Gamma|}(\sum\limits_{\xi\in {\rm
Aut}\Gamma}\lambda (\xi )\prod\limits_{x\in T_{\xi}^V}
(\frac{d(x)}{o(\xi|_{N_{\Gamma}(x)})}-1)!.
\end{eqnarray*}

Similarly, we enumerate the unrooted maps on locally orientable
surface underlying the graph $\Gamma$ by ($1.3.3$) in the
Corollary $1.3.1$ as follows.

\begin{eqnarray*}
n^{L}(\Gamma) &=&\frac{1}{2|{\rm Aut}\Gamma|}\sum\limits_{g \in
{\rm Aut}\Gamma} |\Phi_1^{T} (g)|\\
&=& \frac{1}{2|{\rm Aut}\Gamma|}(\sum\limits_{\xi \in {\rm
Aut}\Gamma}^{|T_{\xi}^E| -|T_{\xi}^V|}\prod\limits_{x\in
T_{\xi}^V}
(\frac{d(x)}{o(\xi|_{N_{\Gamma}(x)})}-1)!\\
&+& \sum\limits_{\varsigma \in {\rm Aut}\Gamma ,o(\varsigma
)\equiv 0(mod 2)}2^{|T_{\varsigma}^E|
-|T_{\varsigma}^V|}\prod\limits_{x\in T_{\varsigma}^V}
(\frac{d(x)}{o(\varsigma|_{N_{\Gamma}(x)})}-1)!)\\
&=& \frac{1}{|{\rm Aut}\Gamma|}\sum\limits_{\xi \in {\rm
Aut}\Gamma}\lambda (\xi )2^{|T_{\xi}^E|
-|T_{\xi}^V|}\prod\limits_{x\in T_{\xi}^V}
(\frac{d(x)}{o(\xi|_{N_{\Gamma}(x)})}-1)!.
\end{eqnarray*}

Notice that $n^{O}(\Gamma)+n^{N}(\Gamma)=n^{L}(\Gamma)$, we get
the number of unrooted maps on non-orientable surfaces underlying
the graph $\Gamma$ as follows.

\begin{eqnarray*}
n^{N}(\Gamma) &=& n^{L}(\Gamma)-n^{O}(\Gamma)\\
&=& \frac{1}{|Aut \Gamma|}\times \sum\limits_{\xi\in Aut \Gamma}
(2^{|T_{\xi}^E| -|T_{\xi}^V|}-1)\lambda (\xi ) \prod\limits_{x\in
T_{\xi}^V} (\frac{d(x)}{o(\xi|_{N_{\Gamma}(x)})}-1)!
\end{eqnarray*}

This completes the proof.\quad\quad $\natural$

If $\Gamma$ is also a $k$-regular graph, we can get a simple
result for the numbers of unrooted maps on orientable or
non-orientable surfaces.

\vskip 4mm

\no{\bf Corollary $4.2.1$} {\it Let $\Gamma$ be a $k$-regular
semi-regular graph. Then the numbers of unrooted maps on
orientable or non-orientable surfaces underlying the graph
$\Gamma$ are respective}

$$
n^O(\Gamma )=\frac{1}{|{\rm Aut} \Gamma|}\times\sum\limits_{g\in
{\rm Aut \Gamma}}\lambda (g) (k-1)!^{|T_g^V|}
$$

\no{\it and}

$$
n^N(\Gamma )=\frac{1}{|{\rm Aut} \Gamma|}\times\sum\limits_{g\in
{\rm Aut \Gamma}}\lambda (g)
(2^{|T_g^E|-|T_g^V|}-1)(k-1)!^{|T_g^V|},
$$

\no{\it where, $\lambda (\xi )=1$ if $o(\xi )\equiv 0(mod 2)$ and
 $\frac{1}{2}$, otherwise.}

\vskip 3mm

{\it Proof} In the proof of the Theorem $4.2$ in $[50]$, it has
been proved that

$$\forall x\in V(\Gamma),\quad\quad o(\xi_{N_{\Gamma}(x)})=1.$$

\no{Whence, we get $n^O(\Gamma )$ and $n^N(\Gamma )$ by the
Theorem $4.2.1$.\quad \quad $\natural$}

Similarly, if $\Gamma ={\rm Cay}(Z_p:S)$ for a prime $p$, we can
also get closed formulas for the number of unrooted maps
underlying the graph $\Gamma$.

\vskip 4mm

\no{\bf Corollary $4.2.2$} {\it Let $\Gamma ={\rm Cay}(Z_p:S)$ be
connected graph of prime order $p$ with $(p-1,|S|)=2$. Then}

$$
n^O(\Gamma,{\cal M})=\frac{(|S|-1)!^p+
2p(|S|-1)!^{\frac{p+1}{2}}+(p-1)(|S|-1)!}{4p}
$$

\no{\it and}

\begin{eqnarray*}
n^N(\Gamma,{\cal M}) &=& \frac{(2^{\frac{p|S|}{2}-p}-1)(|S|-1)!^p+
2(2^{\frac{p|S|-2p-2)}{4}}-1)p(|S|-1)!^{\frac{p+1}{2}}}{2p}\\
&+& \frac{(2^{\frac{|S|-2}{2}}-1)(p-1)(|S|-1)!}{4p}.
\end{eqnarray*}

\vskip 3mm

{\it Proof} By the proof of the Corollary $4.1$ in $[50]$, we have
known that

\[
|T_{g}^V|=\left\{\begin{array}{cc}
1, & {\rm if}\quad o(g)=p\\
\frac{p+1}{2}, & {\rm if}\quad o(g)=2\\
p, & {\rm if}\quad o(g)=1
\end{array}
\right.
\]

\no{and}

\[
|T_{g}^E|=\left\{\begin{array}{cc}
\frac{|S|}{2}, & {\rm if}\quad o(g)=p\\
\frac{p|S|}{4}, & {\rm if}\quad o(g)=2\\
\frac{p|S|}{2}, & {\rm if}\quad o(g)=1
\end{array}
\right.
\]

Notice that ${\rm Aut}\Gamma = D_p$ (see $[3][12]$) and there are
$p$ elements order $2$, one order $1$ and $p-1$ order $p$. Whence,
we have

$$
n^O(\Gamma,{\cal M})=\frac{(|S|-1)!^p+
2p(|S|-1)!^{\frac{p+1}{2}}+(p-1)(|S|-1)!}{4p}
$$

\no{\it and}

\begin{eqnarray*}
n^N(\Gamma,{\cal M}) &=& \frac{(2^{\frac{p|S|}{2}-p}-1)(|S|-1)!^p+
2(2^{\frac{p|S|-2p-2)}{4}}-1)p(|S|-1)!^{\frac{p+1}{2}}}{2p}\\
&+& \frac{(2^{\frac{|S|-2}{2}}-1)(p-1)(|S|-1)!}{4p}.
\end{eqnarray*}

\no{By the Corollary $4.2.1$.\quad\quad $\natural$}

\vskip 10mm

{\bf\S 3. The enumeration of a bouquet on surfaces}

\vskip 8mm

For any integer $k,k|2n$, let ${\mathcal J}_k$ be the conjugacy
class in $S_n[S_2]$ with each cycle in the decomposition of a
permutation in ${\mathcal J}_k$ being $k$-cycle. According to the
Corollary $1.3.1$, we need to determine the numbers $|\Phi^O
(\xi)|$ and $|\Phi^L (\xi)|$ for each automorphism of a map
underlying a graph $B_n$.

\vskip 4mm

\no{\bf Lemma $4.3.1$} {\it Let
$\xi=\prod\limits_{i=1}^{2n/k}(C(i))(\alpha C(i)\alpha^{-1} )\in
{\mathcal J}_k$, where, $C(i)=(x_{i1},x_{i2},\cdots,x_{ik})$ is a
$k$-cycle, be a cycle decomposition. Then

$(i)$\quad If $k\not=2n$, then

 $$|\Phi^O
(\xi)|=k^{\frac{2n}{k}}(\frac{2n}{k}-1)!$$

\no{and if k=2n, then $|\Phi^O (\xi)|=\phi (2n)$.}

$(ii)$\quad If $k\geq 3$ and $k\not= 2n$, then

$$|\Phi^L (\xi)|=(2k)^{\frac{2n}{k}-1}(\frac{2n}{k}-1)!,$$

\no{and }

$$|\Phi^L (\xi)| \ = \ 2^n(2n-1)!$$

\no{if $\xi =(x_1)(x_2)\cdots (x_n)(\alpha x_1)(\alpha x_2)\cdots
(\alpha x_n)(\beta x_1)(\beta x_2)\cdots (\beta x_n)(\alpha\beta
x_1)(\alpha\beta x_2)$ $\cdots (\alpha\beta x_n)$, and}

$$|\Phi^L (\xi)| \ = \ 1$$

\no{if $\xi =(x_1,\alpha\beta x_1)(x_2,\alpha\beta x_2)\cdots
(x_n,\alpha\beta x_n)(\alpha x_1,\beta x_1)(\alpha x_2,\beta
x_2)\cdots(\alpha x_n,\beta x_n)$, and}

$$|\Phi^L (\xi)|= \frac{n!}{(n-2s)!s!}$$

\no{if $\xi ={\zeta ;\varepsilon_1,\varepsilon_2,\cdots
,\varepsilon_n}$ and $\zeta\in{\mathcal E}_{[1^{n-2s},2^{s}]}$ for
some integer $s$, $\varepsilon_i=(1,\alpha\beta )$ for $1\leq
i\leq s$ and $\varepsilon_j=1$ for $s+1\leq j\leq n$,where,
${\mathcal E}_{[1^{n-2s},2^{s}]}$ denotes the conjugate class with
the type $[1^{n-2s},2^{s}]$ in the symmetry group $S_n$, and}

$$|\Phi^L (\xi)|=\phi (2n)$$

\no{if $\xi ={\theta ;\varepsilon_1,\varepsilon_2,\cdots
,\varepsilon_n}$ and $\theta\in {\mathcal E}_{[n^1]}$ and
$\varepsilon_i=1$ for $1\leq i\leq n-1$,
$\varepsilon_n=(1,\alpha\beta)$, where, $\phi (t)$ is the Euler
function.}}

\vskip 3mm

{\it Proof}\quad $(i)$ Notice that for a given representation of
$C(i)$, $i=1,2,\cdots,\frac{2n}{k}$, since $<{\mathcal
P}_n,\alpha\beta
>$ is not transitive on ${\cal X}_{\alpha,\beta }$, there is one
and only one stable orientable map ${\mathcal B}_n=({\cal
X}_{\alpha,\beta },{\mathcal P}_n)$ with $X=E(B_n)$ and ${\mathcal
P}_n=C(\alpha C^{-1}\alpha^{-})$, where,

$$
C=(x_{11},x_{21},\cdots,x_{\frac{2n}{k}1},x_{21},x_{22},\cdots,
x_{\frac{2n}{k}2},x_{1k},x_{2k},\cdots,x_{\frac{2n}{k}k}).
$$

\no{Counting ways of each possible order for $C(i),
i=1,2,\cdots,\frac{2n}{k}$, and different representations for
$C(i)$, we know that}

 $$|\Phi^O
(\xi)|=k^{\frac{2n}{k}}(\frac{2n}{k}-1)!$$

\no{for $k\not=2n$.}

 Now if $k=2n$, then the permutation itself is
a map with the underlying graph $B_n$. Whence, it is also an
automorphism of the map with the permutation is its power.
Therefore, we get that

$$|\Phi^O (\xi)|=\phi (2n)$$.

$(ii)$ For $k\geq 3$ and $k\not= 2n$, since the group $<{\mathcal
P}_n,\alpha\beta>$ is transitive on ${\mathcal X}_{\alpha,\beta }$
or not, we can interchange $C(i)$ by $\alpha C(i)^{-1}\alpha^{-1}$
for each cycle not containing the quadricell $x_{11}$. Notice that
we only get the same map if the two sides of some edges are
interchanged altogether or not. Whence, we get that

$$|\Phi^L (\xi)|=2^{\frac{2n}{k}-1}k^{\frac{2n}{k}-1}(\frac{2n}{k}-1)!
=(2k)^{\frac{2n}{k}-1}(\frac{2n}{k}-1)!.$$

Now if $\xi =(x_1,\alpha\beta x_1)(x_2,\alpha\beta x_2)\cdots
(x_n,\alpha\beta x_n)(\alpha x_1,\beta x_1)(\alpha x_2,\beta
x_2)\cdots(\alpha x_n,\beta x_n)$, there is one and only one
stable map $({\cal X}_{\alpha,\beta },{\mathcal P}^1_n)$ under the
action of $\xi$, where

$${\mathcal P}^1_n=(x_1,x_2,\cdots,x_n,\alpha\beta x_1,\alpha\beta x_2,
\cdots,\alpha\beta x_n) (\alpha x_1,\beta x_n,\cdots,\beta
x_1,\alpha x_n,\cdots, \alpha x_1).$$

\no{Which is an orientable map. Whence, $|\Phi^L (\xi)|=|\Phi^O
(\xi)|=1$. }

If $\xi =(x_1)(x_2)\cdots (x_n)(\alpha x_1)(\alpha x_2)\cdots
(\alpha x_n)(\beta x_1)(\beta x_2)\cdots (\beta x_n)(\alpha\beta
x_1)(\alpha\beta x_2)$\\ $\cdots (\alpha\beta x_n)$, we can
interchange $(\alpha\beta x_i)$ with $(\beta x_i)$ and obtain
different embeddings of $B_n$ on surfaces. Whence, we know that

$$|\Phi^L (\xi)|=2^n(2n-1)!.$$

Now if $\xi =(\zeta ;\varepsilon_1,\varepsilon_2,\cdots
,\varepsilon_n)$ and $\zeta\in{\mathcal E}_{[1^{n-2s},2^{s}]}$ for
some integer $s$, $\varepsilon_i=(1,\alpha\beta )$ for $1\leq
i\leq s$ and $\varepsilon_j=1$ for $s+1\leq j\leq n$, we can not
interchange $(x_i,\alpha\beta x_i)$ with $(\alpha x_i,\beta x_i)$
to get different embeddings of $B_n$ for it just interchanging the
two sides of one edge. Whence, we get that

$$|\Phi^L (\xi)|= \frac{n!}{1^{n-2s}(n-2s)!2^ss!}\times 2^s =\frac{n!}{(n-2s)!s!}. $$

For $\xi =(\theta ;\varepsilon_1,\varepsilon_2,\cdots
,\varepsilon_n)$, $\theta\in {\mathcal E}_{[n^1]}$ and
$\varepsilon_i=1$ for $1\leq i\leq n-1$,
$\varepsilon_n=(1,\alpha\beta)$, we can not get different
embeddings of $B_n$  by interchanging the two conjugate cycles,
whence, we get that

$$|\Phi^L (\xi)|=|\Phi^O (\xi)|=\phi (2n).$$

\no{This completes the proof. \quad\quad $\natural$}}

Recall that the {\it cycle index of a group} $G$, denoted by
$Z(G;s_1,s_2,\cdots,s_n)$, is defined by ($[30]$)

$$
Z(G;s_1,s_2,\cdots,s_n)=\frac{1}{|G|}\sum\limits_{g\in
G}s_1^{\lambda_1 (g)}s_2^{\lambda_2 (g)}\cdots s_n^{\lambda_n
(g)},
$$

\no{where, $\lambda_i (g)$ is the number of $i$-cycles in the
cycle decomposition of $g$. For the symmetric group $S_n$, its
cycle index is known as follows:}

$$
Z(S_n;s_1,s_2,\cdots,s_n)=\sum\limits_{\lambda_1+2\lambda_2+\cdots+k\lambda_k=n}
\frac{s_1^{\lambda_1} s_2^{\lambda_2 }\cdots
s_k^{\lambda_k}}{1^{\lambda_1}\lambda_1!2^{\lambda_2}\lambda_2!
\cdots k^{\lambda_k}\lambda_k!}.
$$

\no{For example, we have that  $Z(S_2)=\frac{s_1^2+s_2}{2}$. By a
result of Polya ($[56]$), we know the cycle index of $S_n[S_2]$ as
follows:}

$$
Z(S_n[S_2];s_1,s_2,\cdots,s_{2n})=\frac{1}{2^nn!}\sum\limits_{\lambda_1+2\lambda_2+\cdots+k\lambda_k=n}
\frac{(\frac{s_1^2+s_2}{2})^{\lambda_1}
(\frac{s_2^2+s_4}{2})^{\lambda_2 }\cdots
(\frac{s_k^2+s_{2k}}{2})^{\lambda_k}}{1^{\lambda_1}\lambda_1!2^{\lambda_2}\lambda_2!
\cdots k^{\lambda_k}\lambda_k!}
$$

Now we can count maps on surfaces with an underlying graph $B_n$.

\vskip 4mm

\no{\bf Theorem $4.3.1$} {\it The number $n^O(B_n)$ of
non-isomorphic maps on orientable surfaces underlying a graph
$B_n$ is }

\begin{eqnarray*}
n^O(B_n) &=& \sum\limits_{k|2n,k\not= 2n}
k^{\frac{2n}{k}-1}(\frac{2n}{k}-1)! \frac{1}{(\frac{2n}{k})!}
\frac{\partial^{\frac{2n}{k}}(Z(S_n[S_2]))}{\partial
s_k^{\frac{2n}{k}}}|_{s_k=0}\\
&+& \phi (2n)\frac{\partial (Z(S_n[S_2]))}{\partial
s_{2n}}|_{s_{2n}=0}
\end{eqnarray*}

{\it Proof}\quad According to the formula (1.3.1) in the Corollary
$1.3.1$, we know that

$$
n^{O}(B_n)=\frac{1}{2\times 2^nn!}\sum\limits_{\xi\in
{S_n[S_2]\times\prec\alpha\succ}} |\Phi^{T} (\xi)|.
$$

Now since for $\forall \xi_1,\xi_2\in S_n[S_2]$, if there exists
an element $\theta\in S_n[S_2]$, such that
$\xi_2=\theta\xi_1\theta^{-1}$, then we have that
$|\Phi^O(\xi_1)|=|\Phi^O(\xi_2)|$ and
$|\Phi^O(\xi)|=|\Phi^O(\xi\alpha)|$.  Notice that $|\Phi^O(\xi)|$
has been gotten by the Lemma $4.3.1$. Using the Lemma $4.3.1 (i)$
and the cycle index $Z(S_n[S_2])$, we get that

\begin{eqnarray*}
n^{O}(B_n)&=& \frac{1}{2\times 2^nn!}(\sum\limits_{k|2n, k\not=
2n}k^{\frac{2n}{k}-1}(\frac{2n}{k}-1)!|{\mathcal J}_{k}|+\phi
(2n)|{\mathcal J}_{2n}|)\\
&=& \sum\limits_{k|2n,k\not= 2n}
k^{\frac{2n}{k}-1}(\frac{2n}{k}-1)! \frac{1}{(\frac{2n}{k})!}
\frac{\partial^{\frac{2n}{k}}(Z(S_n[S_2]))}{\partial
s_k^{\frac{2n}{k}}}|_{s_k=0}\\
&+& \phi (2n)\frac{\partial (Z(S_n[S_2]))}{\partial
s_{2n}}|_{s_{2n}=0} \quad\natural
\end{eqnarray*}

Now we consider maps on the non-orientable surfaces with an
underlying graph $B_n$. Similar to the discussion of the Theorem
$4.1$, we get the following enumeration result for the
non-isomorphic maps on the non-orientable surfaces.

\vskip 4mm

\no{\bf Theorem $4.3.2$} {\it The number $n^N(B_n)$ of
non-isomorphic maps on the non-orientable surfaces with an
underlying graph $B_n$ is }

\begin{eqnarray*}
n^N(B_n)&=& \frac{(2n-1)!}{n!}+\sum\limits_{k|2n,3\leq
k<2n}(2k)^{\frac{2n}{k}-1}(\frac{2n}{k}-1)!
\frac{\partial^{\frac{2n}{k}}(Z(S_n[S_2]))}{\partial
s_k^{\frac{2n}{k}}}|_{s_k=0}\\
 &+& \frac{1}{2^nn!}(\sum\limits_{s\geq
1}\frac{n!}{(n-2s)!s!}+4^n(n-1)!(\frac{\partial^n(Z(S_n[S_2]))}{\partial
s_2^n}|_{s_2=0}-\lfloor\frac{n}{2}\rfloor)) .
\end{eqnarray*}

\vskip 3mm

{\it Proof}\quad Similar to the proof of the Theorem $4.3.1$,
applying the formula ($1.3.3$) in the Corollary $1.3.1$ and the
Lemma $4.3.1(ii)$, we get that

\begin{eqnarray*}
n^L(B_n)&=& \frac{(2n-1)!}{n!}+ \phi (2n)\frac{\partial^n
(Z(S_n[S_2]))}{\partial s_{2n}^n}|_{s_{2n}=0}\\
&+& \frac{1}{2^nn!}(\sum\limits_{s\geq
0}\frac{n!}{(n-2s)!s!}+4^n(n-1)!(\frac{\partial^n(Z(S_n[S_2]))}{\partial
s_2^n}|_{s_2=0}-\lfloor\frac{n}{2}\rfloor))\\
&+& \sum\limits_{k|2n,3\leq
k<2n}(2k)^{\frac{2n}{k}-1}(\frac{2n}{k}-1)!
\frac{\partial^{\frac{2n}{k}}(Z(S_n[S_2]))}{\partial
s_k^{\frac{2n}{k}}}|_{s_k=0}.
\end{eqnarray*}

Notice that $n^O(B_n)+n^N(B_n)=n^L(B_n)$. Applying the result in
the Theorem $4.3.1$, we know that

\begin{eqnarray*}
n^N(B_n)&=& \frac{(2n-1)!}{n!}+\sum\limits_{k|2n,3\leq
k<2n}(2k)^{\frac{2n}{k}-1}(\frac{2n}{k}-1)!
\frac{\partial^{\frac{2n}{k}}(Z(S_n[S_2]))}{\partial
s_k^{\frac{2n}{k}}}|_{s_k=0}\\
&+& \frac{1}{2^nn!}(\sum\limits_{s\geq
1}\frac{n!}{(n-2s)!s!}+4^n(n-1)!(\frac{\partial^n(Z(S_n[S_2]))}{\partial
s_2^n}|_{s_2=0}-\lfloor\frac{n}{2}\rfloor)).
\end{eqnarray*}

\no{This completes the proof. \quad\quad $\natural$}

\includegraphics[bb=10 10 200 230]{12mg2.eps}
\vskip 2mm \c{\bf Fig. $4.2$}

Calculation shows that

$$
Z(S_1[S_2])=\frac{s_1^2+s_2}{2}
$$

\no{and}

$$
Z(S_2[S_2])=\frac{s_1^4+2s_1^2s_2+3s_2^2+2s_4}{8},
$$

For $n=2$, calculation shows that there are $1$ map on the plane,
$2$ maps on the projective plane, $1$ map on the torus and $2$
maps on the Klein bottle. All of those maps are non-isomorphic and
same as the results gotten by the Theorem $4.3.1$ and $4.3.2$,
which are shown in the Fig. $2$.}

\vskip 5mm

{\bf \S $4.$ A classification of the closed $s$-manifolds}

\vskip 3mm

According to the Theorem $1.2.8$, We can classify the closed
$s$-manifolds by triangular maps with valency in $\{5,6,7\}$ as
follows:\vskip 2mm

{\bf Classical Type}:\vskip 2mm

$(1)$  $\Delta_1=\{5-regular \ triangular \ maps\}$ ({\it
elliptic});

$(2)$ $\Delta_2=\{6-regular  \ triangular \ maps\}$({\it euclid});

$(3)$ $\Delta_3=\{7-regular  \ triangular \ maps\}$({\it
hyperbolic}).\vskip 2mm

{\bf Smarandache Type}:\vskip 2mm

$(4)$ $\Delta_4=\{triangular \ maps \ with \ vertex \ valency \ 5
\ and \ 6\}$ ({\it euclid-elliptic});

$(5)$ $\Delta_5=\{triangular \ maps \ with \ vertex \ valency \ 5
\ and \ 7\}$ ({\it elliptic-hyperbolic});

$(6)$ $\Delta_6=\{triangular \ maps \ with \ vertex \ valency \ 6
\ and \ 7\}$ ({\it euclid-hyperbolic});

$(7)$ $\Delta_7=\{triangular \ maps \ with \ vertex \ valency \ 5
, 6 \ and \ 7\}$ ({\it mixed}).

We prove each type is not empty in this section.

\vskip 3mm

\no{\bf Theorem $4.4.1$} {\it For classical types
$\Delta_1\sim\Delta_3$, we have that

$(1)$ $\Delta_1=\{O_{20}, P_{10}\};$

$(2)$ $\Delta_2=\{T_i,K_j, 1 \leq i,j\leq +\infty\};$

$(3)$ $\Delta_3=\{H_i, 1\leq i\leq +\infty\},$

\no{where, $O_{20},P_{10}$ are shown in the Fig.$4.3$, $T_3, K_3$
are shown in the Fig. $4.4$ and $H_i$ is the Hurwitz maps, i.e.,
triangular maps of valency $7$.}}

\includegraphics[bb=10 10 180 180]{sg2.eps}
\vskip 2mm \c{\bf Fig. $4.3$}

\includegraphics[bb=10 10 180 180]{sg3.eps}
\vskip 2mm \c{\bf Fig. $4.4$}

\vskip 3mm

{\it Proof}  If $M$ is a $k$-regular triangulation, we get that
$2\varepsilon (M)=3\phi (M)=k\nu (M)$. Whence, we have

$$\varepsilon (M)=\frac{3\phi (M)}{2} \ {\rm and} \ \nu (M)=\frac{3\varepsilon (M)}{k}.$$

By the Euler-Poincare formula, we know that

$$\chi (M)=\nu (M)-\varepsilon (M)+\phi (M)=(\frac{3}{k}-\frac{1}{2})\phi (M).$$

If $M$ is elliptic, then $k=5$. Whence, $\chi (M)=\frac{\phi
(M)}{10} \ >0$. Therefore, if $M$ is orientable, then $\chi
(M)=2$, Whence, $\phi (M)=20, \nu (M)=12$ and $\varepsilon
(M)=30$, which is the map $O_{20}$. If if $M$ is non-orientable,
then $\chi (M)=1$, Whence, $\phi (M)=10, \nu (M)=6$ and
$\varepsilon (M)=15$, which is the map $P_{10}$.

If $M$ is euclid, then $k=6$. Whence, $\chi (M)=0$, i.e., $M$ is a
$6$-regular triangulation $T_i$ or $K_j$ for some integer $i$ or
$j$ on the torus or Klein bottle, which is infinite.

If $M$ is hyperbolic, then $k=7$. Whence, $\chi (M) \ <0$. $M$ is
a $7$-regular triangulation, i.e., the Hurwitz map. According to
the results in $[65]$, there are infinite Hurwitz maps on
surfaces. This completes the proof. \quad\quad $\natural$

For the Smarandache Types, the situation is complex. But we can
also obtain the enumeration results for each of the types
$\Delta_4\sim\Delta_7$. First, we prove a condition for the
numbers of vertex valency $5$ with $7$.

\vskip 4mm

\no{\bf Lemma $4.4.1$} {\it Let Let ${\mathcal C}(T,n)$ be an
$s$-manifold. Then}

$$v_7\geq v_5+2$$

\no{\it if $\chi ({\mathcal C}(T,n))\leq -1$  and}

$$v_7\leq v_5-2$$

\no{\it if  $\chi ({\mathcal C}(T,n))\geq 1$. where $v_i$ denotes
the number of vertices of valency $i$ in ${\mathcal C}(T,n)$.}

\vskip 3mm

{\it Proof} Notice that we have know

$$\varepsilon ({\mathcal C}(T,n) )= \frac{-\chi ({\mathcal C}(T,n))}{\frac{1}{3}-\frac{2}{k}},$$

\no{where $k$ is the average valency of vertices in ${\mathcal
C}(T,n)$. Since}

$$k = \frac{5v_5+6v_6+7v_7}{v_5+v_6+v_7}$$

\no{and $\varepsilon ({\mathcal C}(T,n) )\geq 3$. Therefore, we
get that}

$(i)$ If $\chi ({\mathcal C}(T,n))\leq -1$, then

$$\frac{1}{3}-\frac{2v_5+2v_6+2v_7}{5v_5+6v_6+7v_7} \ > 0, $$

\no{i.e., $v_7\geq v_5+1$. Now if $v_7=v_5+1$, then}

$$5v_5+6v_6+7v_7=12v_5+6v_6+7\equiv 1(mod 2).$$

\no{Contradicts to the fact that $\sum_{v\in
V(\Gamma)}\rho_{\Gamma}(v)=2\varepsilon (\Gamma)\equiv 0(mod 2)$
for a graph $\Gamma$. Whence we get that}

$$v_7\geq v_5+2.$$

$(ii)$ If $\chi ({\mathcal C}(T,n))\geq 1$, then

$$\frac{1}{3}-\frac{2v_5+2v_6+2v_7}{5v_5+6v_6+7v_7} \ < 0, $$

\no{i.e., $v_7\leq v_5-1$. Now if $v_7=v_5-1$, then}

$$5v_5+6v_6+7v_7=12v_5+6v_6-7\equiv 1(mod 2).$$

\no{Also contradicts to the fact that $\sum_{v\in
V(\Gamma)}\rho_{\Gamma}(v)=2\varepsilon (\Gamma)\equiv 0(mod 2)$
for a graph $\Gamma$. Whence, we have that}

$$v_7\leq v_5-2. \quad\quad \natural$$

\vskip 4mm

\no{\bf Corollary $4.4.1$} {\it There is no an $s$-manifold
${\mathcal C}(T,n)$ such that}

$$|v_7 - v_5|\leq 1,$$

\no{where $v_i$ denotes the number of vertices of valency $i$ in
${\mathcal C}(T,n)$.}

Define an {\it operator $\Xi : M\rightarrow M^*$} on a
triangulation $M$ of a surface as follows:

{\it Choose each midpoint on each edge in $M$ and connect the
midpoint in each triangle as shown in the Fig.$4.5$. Then the
resultant $M^*$ is a triangulation of the same surface and the
valency of each new vertex is $6$. }

\c{\includegraphics[bb=0 0 500 200]{sg4.eps}}\c{\bf Fig. $4.5$}

Then we get the following result.

\vskip 3mm

\no{\bf Theorem $4.4.2$} {\it For the Smarandache Types
$\Delta_4\sim\Delta_7$, we have }

$(i)$ $|\Delta_5|\geq 2$;

$(ii)$ Each of $|\Delta_4|, |\Delta_6|$ and $|\Delta_7|$ is
infinite.}

\vskip 3mm

{\it Proof} For $M\in \Delta_4$, let $k$ be the average valency of
vertices in $M$, since

$$k=\frac{5v_5+6v_6}{v_5+v_6} \ <  6 \ {\rm and} \  \varepsilon (M)=
\frac{-\chi (M)}{\frac{1}{3}-\frac{2}{k}},$$

\no{we have that $\chi (M)\geq 1$. Calculation shows that $v_5=6$
if $\chi (M)=1$ and $v_5=12$ if $\chi (M)=2$. We can construct a
triangulation with vertex valency $5,6$ on the plane and the
projective plane in the Fig. $4.6$.}

\c{\includegraphics[bb=0 0 500 200]{sg5.eps}}  \c{\bf Fig. $4.6$}

Now let $M$ be a map in the Fig. $4.6$. Then $M^{\Theta}$ is also
a triangulation of the same surface with vertex valency $5,6$ and
$M^{\Theta}\not=M$. Whence, $|\Delta_4|$ is infinite.

For $M\in \Delta_5$, by the Lemma $4.4.1$, we know that $v_7\leq
v_5-2$ if $\chi (M)\geq 1$ and $v_7\geq v_5+2$ if $\chi(M)\leq
-1$. We construct a triangulation on the plane and on the
projective plane in the Fig.$4.6$.

\c{\includegraphics[bb=0 0 500 230]{sg6.eps}}  \c{\bf Fig. $4.7$}

For $M\in \Delta_6$, we know that $k=\frac{6v_6+7v_7}{v_6+v_7} \
>6.$ Whence, $\chi (M)\leq -1$. Since $3\phi (M)=6v_6+7v_7=2\varepsilon
(M)$, we get that

$$
v_6+v_7-\frac{6v_6+7v_7}{2}+\frac{6v_6+7v_7}{3}=\chi (M).
$$

\no{Therefore, we have $v_7=-\chi (M)$. Since there are infinite
Hurwitz maps $M$ on surfaces. Then the resultant triangular map
$M^*$ is a triangulation with vertex valency $6,7$ and
$M^*\not=M$. Whence, $|\Delta_6|$ is infinite.}

For $M\in \Delta_7$, we construct a triangulation with vertex
valency $5,6,7$ in the Fig. $4.8$.

Let $M$ be one of the maps in the Fig.$4.8$. Then the action of
$\Theta$ on $M$ results infinite triangulations of valency $5,6$
or $7$. This completes the proof. \quad\quad $\natural$

\vskip 4mm

\no{\bf Conjecture $4.4.1$}\quad {\it The number $| \Delta_5|$  is
infinite.}

\includegraphics[bb=10 10 180 180]{sg7.eps}
\vskip 2mm \c{\bf Fig. $4.8$}

\newpage

%%%%%%%%%%%%%%%Headings%%%%%%%%%%%%%%%%%%%%%%%%%%%%%%%%%%%%%%%%%
\thispagestyle{empty} \pagestyle{myheadings} \topmargin 5mm
\headheight 8mm \headsep 10mm

\markright {\scriptsize  Chapter $5$\quad Open Problems for
Combinatorial Maps}
%%%%%%%%%%%%%%%%%%%%%%%%%%%%%%%%%%%%%%%%%%%%%%%%%%%%%%%%%%%%%%%%

\vskip 35mm

\no{\bf\large Chapter $5$\quad Open Problems for Combinatorial
Maps}

\vskip 25mm

As a well kind of decomposition of a surface, maps are more
concerned by mathematician in the last century, especially in the
enumeration of maps ($[33]-[35]$) and graphs embedding on a
surface ($[22],[35],[53],[70]$). This has its own right,also its
contribution to other branch of mathematics and sciences. Among
those map application to other branch mathematics, one typical
example is its contribution to the topology for the classification
of compact surfaces by one face, or its dual, one vertex maps,
i.e., a bouquet $B_n$ on surfaces. Another is its contribution to
the group theory for finding the Higman-Sims group (a sporadic
simple group ([$6$])) and an one sentence proof of the
Nielsen-Schreier theorem, i.e., {\it every subgroup of a free
group is also free} ($[63]-[64]$). All those applications are to
the algebra, a branch of mathematics without measures. From this
view, the topics discussed in the graph theory can be seen just
the algebraic questions. But our world is full of measures. For
applying combinatorics to other branch of mathematics,  a good
idea is pullback measures on combinatorial objects again, ignored
by the classical combinatorics and reconstructed or make
combinatorial generalization for the classical mathematics, such
as, the algebra, differential geometry, Riemann geometry,
Smarandache geometries, $\cdots$ and the mechanics, theoretical
physics, $\cdots$. For doing this, the most possible way is,
perhaps by the combinatorial maps. Here is a collection of open
problems concerned maps with the Riemann geometry and Smarandache
geometries.

\vskip 8mm

{\bf $5.1$ The uniformization theorem for simple connected Riemann
surfaces}

\vskip 5mm

The {\it uniformization theorem} for simple connected Riemann
surfaces is one of those beautiful results in the Riemann surface
theory, which is stated as follows([$18$]).

\vskip 3mm

{\it If ${\mathcal S}$ is a simple connected Riemann surface, then
${\mathcal S}$ is conformally equivalent to one and only one of
the following three:

$(a)$\quad ${\mathcal C}\bigcup {\infty}$;

$(b)$\quad ${\mathcal C}$;

$(c)$\quad $\triangle =\{z\in {\mathcal C}| |z| \ < 1\}.$}\vskip
2mm

\no{We have proved in the Chapter $2$ that any automorphism of a
map is conformal. Therefore, we can also introduced the conformal
mapping between maps. Then, {\it how can we define the conformal
equivalence for maps enabling us to get the uniformization theorem
of maps? What is the correspondence class maps with the three type
$(a)-(c)$ Riemann surfaces?}}

\vskip 8mm

{\bf $5.2$ The Riemann-Roch theorem}

\vskip 5mm

Let ${\mathcal S}$ be a Riemann surface. A {\it divisor} on
${\mathcal S}$ is a formal symbol

$$
{\mathcal U}=\prod\limits_{i=1}^kP_i^{\alpha (P_i)}
$$

\no{with $P_i\in {\mathcal S}$, $\alpha (P_i)\in {\bf Z}$. Denote
by $Div({\mathcal S})$ the free commutative group on the points in
${\mathcal S}$ and define}

$$deg{\mathcal U}=\sum\limits_{i=1}^k\alpha (P_i).$$

Denote by ${\mathcal H}({\mathcal S})$ the field of meromorphic
function on ${\mathcal S}$. Then for $\forall f\in {\mathcal
H}({\mathcal S})\setminus \{0\}$, $f$ determines a divisor $(f)\in
Div({\mathcal S})$ by

$$
(f)=\prod\limits_{P\in {\mathcal S}}P^{ord_Pf},
$$

\no{where, if we write $f(z)=z^ng(z)$ with $g$ holomorphic and
non-zero at $z=P$, then the $ord_Pf=n$. For ${\mathcal
U}_1=\prod\limits_{P\in {\mathcal S}}P^{\alpha_1 (P)},{\mathcal
U}_2=\prod\limits_{P\in {\mathcal S}}P^{\alpha_2 (P)}, \in
Div({\mathcal S})$, call ${\mathcal U}_1\geq {\mathcal U}_2$ if
$\alpha_1(P)\geq\alpha_2(P)$. Now we define a vector space}

$$L({\mathcal U})=\{f\in {\mathcal H}({\mathcal S})|(f)\geq
{\mathcal U}, {\mathcal U}\in Div({\mathcal S})\}$$

$$\Omega ({\mathcal U})=\{\omega| \omega \ is \ an \ abelian \ differential \ with \ (\omega)\geq {\mathcal U}\}.$$

The Riemann-Roch theorem says that($[71]$)

$$dim(L({\mathcal U}^{-1}))=deg{\mathcal U}-g({\mathcal S})+1+dim\Omega ({\mathcal S}).$$

Comparing with the divisors and their vector space, there ia also
cycle space and cocycle space in graphical space theory ([$35$]).
{\it What is their relation? Whether can we rebuilt the
Riemann-Roch theorem by map theory?}

\vskip 8mm

{\bf $5.3$ Combinatorial construction of an algebraic curve of
genus}

\vskip 6mm

A {\it complex plane algebraic curve ${\mathcal C}_l$} is a
homogeneous equation $f(x,y,z)=0$ in $P_2{\mathcal
C}=(C^2\setminus (0,0,0))/\sim$, where $f(x,y,z)$ is a polynomial
in $x,y$ and $z$ with coefficients in ${\mathcal C}$. The degree
of $f(x,y,z)$ is said the {\it degree of the curve} ${\mathcal
C}_l$. For a Riemann surface $S$, a well-known result is
($[71]$){\it there is a holomorphic mapping $\varphi :
S\rightarrow P_2{\mathcal C}$ such that $\varphi (S)$ is a complex
plane algebraic curve and}

$$g(S)=\frac{(d(\varphi (S))-1)(d(\varphi (S))-2)}{2}.$$

By map theory, we know a combinatorial map also is on a surface
with genus. Then {\it whether we can get an algebraic curve by all
edges in a map or by make operations on the vertices or edges of
the map to get plane algebraic curve with given $k$-multiple
points?} and {\it how do we find the equation $f(x,y,z)=0$?}

\vskip 8mm

{\bf $5.4$ Classification of $s$-manifolds by maps}

\vskip 6mm

We have classified the closed $s$-manifolds by maps in the Section
$4$ of Chapter $4$. For the general $s$-manifolds, their
correspondence combinatorial model is the maps on surfaces with
boundary, founded by Bryant and Singerman in $1985$ ([8]). The
later is also related to the modular groups of spaces and need to
investigate further itself. The questions are\vskip 3mm

$(i)$ {\it how can we combinatorially classify the general
$s$-manifolds by maps with boundary?}

$(ii)$ {\it how can we find the automorphism group of an
$s$-manifold?}

$(iii)$ {\it how can we know the numbers of non-isomorphic
$s$-manifolds, with or without root?}

$(iv)$ {\it find rulers for drawing an $s$-manifold on a surface,
such as, the torus, the projective plane or Klein bottle, not the
plane.}\vskip 2mm

The Smarandache manifolds only using the triangulations of
surfaces with vertex valency in $\{5,6,7\}$. Then {\it what are
the geometrical mean of the other maps, such as, the $4$-regular
maps on surfaces.} It is already known that the later is related
to the Gauss cross problem of curves($[35]$). May be we can get a
geometry even more general than that of the Smarandache
geometries.

\vskip 8mm

{\bf $5.5$ Gauss mapping among surfaces}

\vskip 6mm

In the classical differential geometry, a {\it Gauss mapping}
among surfaces is defined as follows([42]):\vskip 3mm

{\it Let ${\mathcal S}\subset R^3$ be a surface with an
orientation {$\bf N$}. The mapping $ N: {\mathcal S}\rightarrow
R^3$ takes its value in the unit sphere}

$$S^2=\{(x,y,z)\in R^3|x^2+y^2+z^2=1\}$$

\no{\it along the orientation {$\bf N$}. The map $ N: {\mathcal
S}\rightarrow S^2$, thus defined, is called the Gauss
mapping.}\vskip 3mm

we know that for a point $P\in {\mathcal S}$ such that the
Gaussian curvature $K(P)\not=0$ and $V$ a connected neighborhood
of $P$ with $K$ does not change sign,

$$K(P)=\lim_{A\rightarrow 0}\frac{N(A)}{A},$$

\no{where $A$ is the area of a region $B\subset V$ and $ N(A)$ is
the area of the image of $B$ by the Gauss mapping $ N: {\mathcal
S}\rightarrow S^2$. The questions are}

{\it $(i)$ what is its combinatorial meaning of the Gauss mapping?
How to realizes it by maps?

$(ii)$ how can we define various curvatures for maps and rebuilt
the results in the classical differential geometry?}

\vskip 8mm

{\bf $5.6$ The Gauss-Bonnet theorem}

\vskip 6mm

{\it Let ${\mathcal S}$ be a compact orientable surface. Then}

$$\int\int_{\mathcal S}Kd\sigma =2\pi\chi ({\mathcal S}),$$

\no{\it where $K$ is Gaussian curvature on ${\mathcal S}$.}

 This is the famous {\it Gauss-Bonnet theorem} for compact surface
($[14],[71]-[72])$. This theorem should has a combinatorial form.
The questions are

{\it $(i)$  how can we define the area of a map? }(Notice that we
give a definition of non-Euclid area of maps in Chapter $2$.)

{\it ($ii$) can we rebuilt the Gauss-Bonnet theorem by maps?}

\vskip 8mm

{\bf $5.7$ Riemann manifolds}

\vskip 6mm

A Riemann surface is just a Riemann $2$-manifold. A {\it Riemann
$n$-manifold} $(M,g)$ is a $n$-manifold $M$ with a Riemann metric
$g$. Many important results in Riemann surfaces are generalized to
Riemann manifolds with higher dimension ($[14],[71]-[72]$). For
example, let ${\mathcal M}$ be a complete, simple-connected
Riemann $n$-manifold with constant sectional curvature $c$, then
we know that {\it ${\mathcal M}$ is isometric to one of the model
spaces ${\mathcal R}^n, S_{{\mathcal R}^n}$ or $H_{{\mathcal
R}^n}$}. There is also a combinatorial map theory for higher
dimension manifolds (see $[67]-[68]$). {\it Whether can we
systematically rebuilt the Riemann manifold theory by
combinatorial maps? or  can we make a combinatorial generalization
of results in the Riemann geometry, for example, the
Chern-Gauss-Bonnet theorem ($[14],[37],[71]$)? } If we can, a new
system for the Einstein's relative theory will be found.

\newpage

%%%%%%%%%%%%%%%Headings%%%%%%%%%%%%%%%%%%%%%%%%%%%%%%%%%%%%%%%%%
\thispagestyle{empty} \pagestyle{myheadings} \topmargin 5mm
\headheight 8mm \headsep 10mm

\markright {\scriptsize  References}
%%%%%%%%%%%%%%%%%%%%%%%%%%%%%%%%%%%%%%%%%%%%%%%%%%%%%%%%%%%%%%%%

\vskip 35mm

\no{\bf\large References}\vskip 5mm

\re{[1]}R.D.M.Accola, On the number of automorphisms of a closed
Riemann surface, {Tran.Amer.Math.Soc.},vol.131,398-408(1968).

\re{[2]}N.L.Alling, {\it Foundation of the theory of Klein
surfaces}, Lect.notes in Math.,219, Springer-Verlag,
Berlin,etc.(1971).

\re{[3]}B.Alspach, Point-symmetric graphs and digraphs of prime
order and transitive permutation groups of prime degree,
J.Combin.Theory,Ser.B, 15(1973),12-17.

\re{[4]}C.Ashbacher, Smarandache geometries, {\it Smarandache
Notions Journal}, Vol.8, No. 1-2-3(1997),212-215.

\re{[5]}L.Babai,Automorphism groups, isomorphism, reconstruction,
in R. Graham, M.Grotschel and L.Lovasz ed: {\it Handbook of
Combinatorics},Elsevier Science B.V, (1995), 1447-1540.

\re{[6]} N.L.Biggs and A.T.White, {\it Permutation Groups and
Combinatoric Structure}, Cambridge University Press (1979).

\re{[7]}L.Brankovic, M.Miller et al, A note on constructing large
Cayley graphs of given degree and diameter by voltage
assignments,{\it The Electronic J. Combinatorics}, 5(1998),\#R9.

\re{[8]}R.P.Bryant and D.Singerman, Foundations of the theory of
maps on surfaces with boundary,{\it
Quart.J.Math.Oxford}(2),36(1985), 17-41.

\re{[9]}E.Bujalance, Cyclic group automorphisms of compact
non-orientable Klein surfaces without boundary,{\it Pacific J.
Mathematics}, vol.109, No.2,279-289(1983).

\re{[10]}E. Bujalance,J.J.Etayo and J.M.Gamboa, Hyperelliptic
Klein surfaces, {\it Qurt. J.Math.Oxford}(2),36(1985), 141-157.

\re{[11]}E.Bujalance,J.J.Etayo,J.M.Gamboa and G.Gromadzki,{\it
Automorphism groups of compact bordered Klein surfaces},
Lect.notes in Math.,1439, Springer-Verlag, Berlin,etc.(1990).

\re{[12]}C.Y.Chao, On the classification of symmetrical graphs
with a prime number of vertices, Trans.Amer.Math.Soc. 158(1971),
247-256.

\re{[13]}J.Chern, J.L.Gross and R.G.Rieper, Overlap matrices and
total imbedding distributions, {\it Discrete
Math.}128(1994),73-94.

\re{[14]}S.S.Chern  and W.H.Chern, {\it Lectures in Differential
Geometry}£¬Peking University Press, 2001.

\re{[15]}B.P.Chetia and K.Patra, On metabelian groups of
automorphisms of compact Riemann surfaces, {\it J.London
Math.Soc}, Vol.33,No.2, 467-472(1986).

\re{[16]}S.Chimienti and M.Bencze, Smarandache paradoxist
geometry, {\it Bulletin of Pure and Applied Sciences}, Delhi,
India, Vol.17E, No.1(1998), 123-1124.

\re{[17]}J.Edmonds, A combinatorial representation for polyhedral
surfaces, {\it Notices Amer. Math. Soc} 7 (1960)

\re{[18]}H. M.Farkas and I. Kra, {\it Riemann Surfaces},
Springer-Verlag New York inc(1980).

\re{[19]}M.L.Furst, J.L.Gross and R.Statman, Genus distributions
for two class of graphs, {\it J.Combin.Theory, Ser B},
46(1989),22-36.

\re{[20]}A.Gardiner,R.Nedela,J.\v{S}ir\'{a}\v{n} and
M.\v{S}kovera, characterization of graphs which underlie regular
maps on closed surfaces, {\it J.London Math.
Soc.}(2)59(1999),100-108.

\re{[21]}G.Gromadzki, Abelian Groups of automorphisms of compact
non-orientable Klein surfaces without boundary, {\it
Commentationes Mathematicae}, 28(1989),197-217.

\re{[22]}J.L.Gross and T.W.Tucker, {\it Topological graph theory},
John Wiley \& Sons,1987.

\re{[23]} J.L.Gross and M.L.Furst, Hierarchy for
imbedding-distribution invariants of a graph, {\it J.Graph
Theory},11(1987), 205-220.

\re{[24]}J.L.Gross and M.L.Furst, Genus distribution for bouquets
of circles, {\it J.Combin. Theory, Ser B}, 47(1989), 292-306.

\re{[25]}F.Harary and W.T.Tutte, On the order of the group of a
planar map, {\it J.Combin. Theory}, 1(1966), 194-195.

\re{[26]}W.J.Harwey,Cyclic groups of automorphisms of compact
Riemann surfaces, {\it Quart J.Math. Oxford}, vol.17, No.2,
86-97(1966).

\re{[27]}H.Iseri, {\it Smarandache Manifolds}, American Research
Press, Rehoboth, NM, 2002.

\re{[28]}H.Iseri, {\it Partially Paradoxist Smarandache
Geometries}, http://www.gallup.unm.
edu/\~smarandache/Howard-Iseri-paper.htm.

\re{[29]}G.A.Jones and D.Singerman, Theory of maps on orientable
surfaces, {\it Proc.London Math. Soc.}(3),37(1978),273-307.

\re{[30]}V.Krishnamurthy, {\it Combinatorics: Theory and
Application}, Ellis Horwood Limited,1986.

\re{[31]}J.H.Kwak and S.H.Shim, Total embedding distribution for
bouquets of circles, {\it Discrete Math.}248(2002),93-108.

\re{[32]}L.Kuciuk and M.Antholy, An Introduction to Smarandache
Geometries, {\it Mathematics Magazine, Aurora, Canada},
Vol.12(2003), and online:

\c{http://www.mathematicsmagazine.com/1-2004/Sm\_Geom\_1\_2004.htm;}

also at New Zealand Mathematics Colloquium, Massey University,
Palmerston

North, New Zealand, December 3-6,2001

\c{http://atlas-conferences.com/c/a/h/f/09.htm;}

also at the International Congress of Mathematicians (ICM2002),
Beijing, China,

20-28, August, 2002,

\c{http://www.icm2002.org.cn/B/Schedule.Section04.htm.}

\re{[33]}Y.P.Liu, {\it Enumerative Theory of Maps}, Kluwer
Academic Publishers , Dordrecht/ Boston/ London,1999.

\re{[34]}Y.P.Liu, {\it Advances in Combinatorial Maps}, Northern
Jiaotong University Publisher, Beijing (2003).

\re{[35]}Y.P.Liu, {\it Embeddability in Graphs}, Kluwer Academic
Publisher, Dordrecht / Boston / London (1995).

\re{[36]} Y.P. Liu, The maximum non-orientable genus of a graph
(in Chinese), {\it{Scientia Sinica (Special Issue on
Math)}},{\bf{I}}(1979),191-201.

\re{[37]}J.M.Lee, {\it Riemann Manifolds}, Springer-Verlag New
York,Inc(1997).

\re{[38]}Lv Y.N. and Zhang X.L.{\it Riemann Surfaces}£¬Sciences
Publisher Press£¬Beijing, 1991.

\re{[39]}C.Maclachlan, Abelian groups of automorphisms of compact
Riemann surfaces, {\it Proc.London Math.Soc.},
vol.15,No.3,699-712(1965).

\re{[40]}A.Malnic, Group action,coverings and lifts of
automorphisms, {\it Discrete Math}, 182(1998),203-218.

\re{[41]}A.Malinic,R.Nedela and M.$\check{S}$koviera, Lifting
graph automorphisms by voltage assignment, {\it Europ.
J.Combinatorics}, 21(2000),927-947.

\re{[42]}Mantredo P.de Carmao, {\it Differential Geometry of
Curves and Surfaces}, Pearson Education asia Ltd (2004).

\re{[43]}L.F.Mao, {\it A census of maps on surface with given
underlying graph}, A doctor thesis in Northern Jiaotong
University, Beijing, 2002.

\re{[44]}L.F.Mao and Y.P.Liu, New automorphism groups identity of
trees, {\it Chinese Advance in Math.},113-117,5(2002).

\re{[45]}L.F.Mao and Y.P.Liu, On the roots on orientable
embeddings of graph, {\it Acta Math. Scientia},
23A(3),287-293(2003).

\re{[46]}L.F.Mao and Y.P.Liu, Group action approach for
enumerating maps on surfaces,{\it J.Applied Math. \& Computing},
vol.13, No.1-2,201-215.

\re{[47]}L.F.Mao and Y.P.Liu, A new approach for enumerating maps
on orientable surfaces, {\it Australasian J. Combinatorics},
vol.30(2004), 247-259.

\re{[48]} L.F.Mao and Y.P.Liu, The semi-arc automorphism group of
a graph with application to map enumeration, {\it Graphs and
Combinatorics}(accepted and to appear).

\re{[49]}L.F.Mao and Y.P.Liu, Automorphisms of maps with
underlying Cayley graph and their application to enumeration,
Submitted to {\it Discrete Math}.

\re{[50]} L.F.Mao, Y.P.Liu and F.Tian, Automorphisms of maps with
a given underlying graph and their application to enumeration,
{\it Acta Mathematica Sinica}, Vol.21, 2(2005),225-236.

\re{[51]}L.F.Mao and F.Tian, On oriented $2$-factorable graphs,
{\it J.Appl.Math. \& Computing}, Vol.17(2005), 25-38.

\re{[52]}W.S.Massey, {\it Algebraic topology: an introduction},
Springer-Verlag,New York, etc.(1977).

\re{[53]}B.Mohar and C.Thomassen, {\it Graphs on Surfaces}, The
Johns Hopkins University Press, London, 2001.

\re{[54]}B.P.Mull,R.G.Rieper and A.T.White, Enumeration 2-cell
imbeddings of connected graphs,{\it Proc.Amer.Math.Soc.},
103(1988), 321~330.

\re{[55]}B.P.Mull, Enumerating the orientable $2$-cell imbeddings
of complete bipartite graphs, {\it J.Graph Theory}, vol 30,
2(1999),77-90.

\re{[56]}R.Nedela and M $\check{S}$koviera, Regular embeddings of
canonical double coverings of graphs, {\it J.combinatorial
Theory},Ser.B 67, 249-277(1996).

\re{[57]}R.Nedela and M.$\check{S}$koviera, Regular maps from
voltage assignments and exponent groups, {\it
Europ.J.Combinatorics}, 18(1997),807-823.

\re{[58]}PlanetMath, {\it Smarandache Geometries},
http://planetmath.org/encyclopedia/ SmarandacheGeometries.htm1.

\re{[59]}K.Polthier and M.Schmies, Straightest geodesics on
polyhedral surfaces, in {\it Mathematical Visualization} (ed. by
H.C.Hege and K.Polthier), Springer-Verlag, Berlin, 1998.

\re{[60]}D.Singerman, Automorphisms of compact non--orientable
Riemann surface, {\it Glasgow J. math.,} 12(1971), 50-59.

\re{[61]}F.Smarandache ,Paradoxist mathematics, {\it Collected
Papers},Vol.,II, 5-28, University of Kishinev Press, 1997.

\re{[62]} S.Stahl, Generalized embedding schemes, {\it J.Graph
Theory} 2,41-52,1978.

\re{[63]}T.R.Stallings, Graphical theory of automorphisms of free
groups, in {\it Combinatorial group theory and topology}(ed. by
S.M.Gersten and J.R.Stallings), Princeton University Press,1987.

\re{[64]}J.Stillwell, {\it Classical topology and combinatorial
group theory}, Springer-Verlag New York Inc., (1980).

\re{[65]}D.B.Surowski, Lifting map automorphisms and MacBeath's
theorem, {\it J. Combin. Theory}, Ser B,50(1990),135-149.

\re{[66]}W.T.Tutte, What is a maps? in New Directions in the
Theory of Graphs (ed.by F.Harary), Academic Press (1973), 309~325.

\re{[67]}A.Vince, Combinatorial maps,{\it J. Combin. Theory}, Ser
B 34 (1983), 1-21.

\re{[68]}A.Vince, Regular combinatorial maps,{\it J. Combin.
Theory}, Ser B 35 (1983),256-277.

\re{[69]}J.K.Weeks, {\it The shape of Space}, New York, Marcel
Dekkler, Inc, 1985.

\re{[70]}A.T.White, {\it Graphs of Group on Surfaces- interactions
and models}, Elsevier Science B.V. (2001).

\re{[71]}H.X.Wu, Y.N.Lv and Z.H.Chern£¬{\it Introduction to
Compact Riemann Surfaces}, Science Publisher Press, Beijing, 1999.

\re{[72]}H.X.Wu, C.L.Shen and Y.L.Yu,{\it Elementary Riemann
Geometry}, Peking University Press, 1989.

\re{[73]}M.Y.Xu, {\it Introduction to the Finite Groups ({\rm\bf
I)}}, Science Publisher Press, Beijing, 1999.

\end{document}